\documentclass[runningheads]{llncs}
\usepackage{graphicx}
\usepackage{float}
\usepackage{fullpage}
\usepackage{mathtools}
\usepackage{amssymb, amsmath, latexsym, bm}
\usepackage{nicefrac}
\usepackage{hhline}
\usepackage{multirow}
\usepackage{pifont}% http://ctan.org/pkg/pifont

\usepackage[colorinlistoftodos,bordercolor=orange,backgroundcolor=orange!20,linecolor=orange,textsize=scriptsize]{todonotes}
 \usepackage{algorithm}
 \usepackage{algpseudocode}
\usepackage{caption}
\usepackage{subcaption}
\usepackage{enumitem}
\newcommand{\argmin}{\mathop{\arg\!\min}}

\newcommand{\EndProof}{\begin{flushright}$\square$\end{flushright}}

  \newcommand{\ab}[1]{{\color{black} #1}} %Aleksandr's new text
  \newcommand{\abc}[1]{{\color{black} #1}} %Aleksandr's new text

\usepackage{chngcntr}

\usepackage[utf8]{inputenc}
\usepackage[T1]{fontenc} 
\usepackage{makecell}

\usepackage[colorlinks=true,linkcolor=red, citecolor=blue, urlcolor=blue]{hyperref}%

\newcommand{\cA}{{\cal A}}

\newcommand{\cM}{{\cal M}}

\newcommand{\EE}{\mathbf{E}}

\def\R{\mathbb{R}}

\def\R{\mathbb R}

\def\EE{\mathbb E}

\def\e{\varepsilon}
\def\la{\langle}
\def\ra{\rangle}

\def\y{\mathbf{y}}
\def\a{\mathbf{a}}
\def\b{\mathbf{b}}
\def\x{\mathbf{x}}
\def\w{\mathbf{w}}
\def\z{\mathbf{z}}
\def\u{\mathbf{u}}
\def\g{\mathbf{g}}
\def\c{\mathbf{c}}

\def\lm{\lambda}

\newcommand{\eigWmax}{\lambda_{\max}(W)}
\newcommand{\eigWmin}{\lambda_{\min}^+(W)}
\newcommand{\eigWhatmax}{\lambda_{\max}(\hat W)}
\newcommand{\eigWhatmin}{\lambda_{\min}^+(\hat W)}

\newtheorem{assumption}{Assumption}

\usepackage{booktabs}

\usepackage{caption}
\usepackage{float}
\usepackage{capt-of}

\usepackage{array}
\usepackage{arydshln}
\usepackage{rotating}

\begin{document}
\title{Decentralized Personalized Federated Learning: Lower Bounds and Optimal Algorithm for All Personalization Modes
}
\titlerunning{Decentralized Personalized Federated Learning}
\author{Abdurakhmon Sadiev\inst{1, 2} \and
Ekaterina Borodich\inst{1,3} \and
Aleksandr Beznosikov\inst{1, 2, 3} \and
Darina Dvinskikh\inst{1,2,4}\and
\\
Saveliy Chezhegov\inst{1} \and
Rachael Tappenden \inst{5} \and
 Martin Tak\'a\v{c}\inst{2} \and
Alexander Gasnikov \inst{1,2,4}}
\authorrunning{A. Sadiev et al.}
\institute{Moscow Institute of Physics and Technology, Russia \and
Mohamed bin Zayed University of Artificial Intelligence (MBZUAI),
UAE \and
HSE University, Russia \and
Institute for Information Transmission Problems RAS, Russia
\and University of Canterbury, New Zealand
}
\maketitle              % typeset the header of the contribution

\begin{abstract}

This paper considers the problem of decentralized, personalized federated learning. For centralized personalized federated learning, a penalty that measures the deviation from the local model and its average, is often added to the objective function. However, in a decentralized setting this penalty is expensive in terms of communication costs, so here, a different penalty --- one that is built to respect the structure of the underlying computational network --- is used instead. We present lower bounds on the communication and local computation costs for this problem formulation and we also present provably optimal methods for decentralized personalized federated learning. Numerical experiments are presented to demonstrate the practical performance of our methods.

\normalsize

\keywords{Federated Learning, Decentralized Optimization, Distributed Optimization, Lower and Upper Bounds, Accelerated Algorithms}

\end{abstract}

\section{Introduction}
Today's data revolution is transforming the world, with vast amounts of data collected daily from a wide range of sources. Automation is necessary when processing and extracting information from such large quantities of data, and machine learning has proven to be a useful tool to assist with this task. The essence of machine learning is to build, and then train, models. To speed up the process of training, modern computer architectures can be used, where instead of one single computing device, the problem and associated data is shared among many devices/agents. This leads to the following distributed learning/optimization problem formulation 
\begin{equation}\label{main_problem0}
    \min_{ \substack{  x \in \R^{d} }} \sum^n_{i = 1} f_i(x),
\end{equation}
where each agent/device $i \in [n]:=\{1,2,\dots,n\}$, has an associated  local loss function $f_i:\R^d\to\R$, as well as its own locally stored data. 

Federated Learning (FL) \cite{konevcny2016federated,mcmahan2017communication} is a subset of distributed machine learning, where one assumes that computing agents are simply general user devices, (for example, smartphones, tablets, laptops, personal computers), and where different devices may have different memory capacity and computing power. This leads to many new and important problems and questions that did not arise previously in the classical distributed setting \cite{kulkarni2020survey}. For example, data may be spread unequally between devices, privacy considerations may prohibit the sharing of data between certain devices on the network, poor or unreliable connectivity may inhibit the flow of data, and data on certain devices may be of poorer quality compared with others.

% In this paper, we focus on one of your biggest federated learning issues - personalization. In \eqref{main_problem0}, we find parameters of the model $x$, based on global data (from all agents). Meanwhile, each individual user can be absolutely unimportant how the global model works, because his local data differs from the data of other users and in his particular case the global model can work poorly. On the other hand, each user may have very little local data, and the process of training a model solely on local data may also be of poor quality. The question of  how to combine the two extremes: global and local, as well as how to vary the proportions of their participation, gave a start to the direction of {\it Personalized  Federated Learning (PFL)}   \cite{hanzely2020federated,hanzely2020lower,hanzely2021personalized}. 
Given these issues, this work focuses on \emph{personalization} for federated learning. Notice that in \eqref{main_problem0}, the model parameter $x$ is found using global data (from all agents). However, the inclusion of a particular individual agent might negatively impact the global training process if their local data differs markedly from the global data, or if they have low quality local data; in this case training the global model may lead to a poor solution. On the other hand, each user may have very little local data, and the process of training a model solely on local data may also lead to a poor quality solution. The question of  how to balance the two extremes, \emph{global versus local}, gave rise to {\it Personalized  Federated Learning (PFL)}   \cite{hanzely2020federated,hanzely2020lower,hanzely2021personalized}.

For PFL, each agent $i$ has their own parameter $x_i$, but the discrepancy between the parameters held on different devices is penalized. Correspondingly, PFL can be formulated as the following regularized optimization problem:
%  \begin{equation}\label{eq:gen_probFL}
% \min_{\x = [x_1, \dots, x_n] \in \R^{n d}} \sum^n_{i = 1} f_i(x_i) +   \frac{\lm}{2} r(\x),
% \end{equation}
 \begin{equation}\label{eq:gen_probFL}
\min_{\x \in \R^{n d}} \sum^n_{i = 1} f_i(x_i) +   \frac{\lm}{2} r(\x),
\end{equation}
where the vector $\x = [x_1^T, \dots, x_n^T]^T \in \R^{nd}$ is the concatenation of the $n$ local vectors $\{x_i \in \R^d\}_{i=1}^n$, $r(\x)$ is a convex penalty function, and the weight parameter $\lambda$ balances the degree of personalization. There are many possible choices for the penalty function $r(\x)$. A simple option is to let $r(\x)$ be the deviation between the local models and their average \cite{hanzely2020lower,hanzely2020federated,hanzely2021personalized}:
\begin{equation}
    r(\x) = \sum_{i=1}^n \|x_i - \bar x\|^2, 
\quad \mbox{where} \   \bar x = \frac1n (x_1+\dots+x_n).
    \label{cent_penalt}
\end{equation}
% This setting is good in a centralized distributed case, when devices communicate with the central server, send information to it and receive responses without failures. In this case, calculating $\bar x$ is easy to realize -- one just need to send all $x_i$ to the server and get an averaged response. In this paper we consider more general setup.
This is a reasonable choice in the centralized distributed setting, where devices communicate with a central server, sending and receiving information without failures. In this case, calculating the average $\bar x$ is easy: all agents $i$ simply send their local $x_i$ to the central server, which then calculates the average $\bar x$, and communicates it back to every agent. In this paper we consider a more general setup, where different penalty functions might be more appropriate.

% \paragraph{Decentralized Setting.}  We assume that now there is no main node (server), and all devices are connected into some large network, which is an undirected connected graph. Moreover, each agent can only communicate with neighbors on this network. Such a setup includes a centralized one - it is enough to consider the complete graph. But in fact, a decentralized setup is more general and closer to federated setting, when devices can communicate only with some other agents and are not accessible to others due to the quality of connection or remoteness.
%\paragraph{Decentralized Setting.}  
Throughout this work we consider \emph{decentralized distributed learning}, where there is no main server (node), but instead all devices are connected via some large network. Moreover, each agent in the network can only communicate with its neighbors. Mathematically, the network is represented by a \ab{fixed}, un-directed, connected graph, where each node corresponds to an agent, and connections between agents are represented by edges. Although a decentralized setting is assumed here, our problem formulation is general enough to include a centralized set-up as a special case (simply take a complete graph). However, a decentralized setup perhaps better captures the federated learning setting, where each device only communicates with a limited number of other agents, corresponding to an incomplete graph. As previously mentioned, communication links between certain agents may be inaccessible, for example, due to the (poor) quality of connection between agents, or due to remoteness of location, and this leads to missing edges in the graph.

In a decentralized setting, using the penalty in \eqref{cent_penalt} is not sensible because of the impracticality of calculating the average $\bar x$. (Note that to calculate $\bar x$, all local $x_i$'s must be sent to one device (node) and then the average $\bar x$ broadcast back to every node, which is a long and expensive operation, especially for large networks.) 
%Including because of this, in classical decentralized optimization, it is unacceptable to calculate the average value of the network, and all agents store their local copies, which converge to the same value during the operation of the algorithm.
With this in mind, here we propose the use of a different penalty $r(\x)$, which is more suitable for a decentralized setup. Hence, the problem formulation considered in this work is:
\begin{equation}
\label{rewrite_main_problem}
    \min_{\x \in \R^{nd}}F(\x) =  \underbrace{\sum^n_{i = 1} f_i(x_i)}_{f(\x)} + \underbrace{\frac{\lambda}{2}\langle \x, W\x\rangle}_{g(\x)},
\end{equation}
\ab{where $W$ is a communication matrix that reflects the properties of the network (see Section \ref{sec:matrix} for a formal definition of $W$). The function $r(\x) = \langle \x, W\x\rangle$ penalizes the difference between neighboring local models in the network, and is computationally friendlier than \eqref{cent_penalt} in a decentralized setting. The matrix $W$ determines how much an agent depends on each of the other nodes in the learning process. This is achieved due to the fact that $W$ represents the structure of the communication graph, gives information about the remoteness of the nodes, the speed of transfer between them, and carries weights of how much to rely on one or another neighbour in the network.}
% Moreover, the penalization for the dissimilarity is fulfilled by the fact $W\x = 0$ if and only if $x_1=\dots = x_n$. The idea of using this type of penalty is not new and has been used in the literature in several contexts, in particular for classical decentralized minimization \cite{li2020decentralized,gorbunov2019optimal,beznosikov2021decentralized} with large $\lambda$ and for multitask PFL \cite{smith2017federated,wang2018distributed,beznosikov2021decentralized} with small $\lambda$.
Note that $W\x = 0$ (and consequently $r(\x) = 0$) if and only if $x_1=\dots = x_n$. This penalty function is not new and has been used in the literature in several contexts, for example, for classical decentralized minimization with large $\lambda$ \cite{li2020decentralized,gorbunov2019optimal,beznosikov2021decentralized}, and for multitask PFL with small $\lambda$ \cite{smith2017federated,wang2018distributed,beznosikov2021decentralized}.

The parameter $\lambda$ balances the `global vs local' trade-off. For example, consider the following extremes:
 \begin{itemize}[label=$\ast$]
    \item If $\lambda = 0$, then \eqref{rewrite_main_problem} becomes $\min_{\x \in \R^{nd}}\sum^n_{i = 1} f_i(x_i)$,
    % \begin{equation*}
    %     \min_{\x \in \R^{nd}}\sum^n_{i = 1} f_i(x_i).
    % \end{equation*}
    where the local function $f_i$ held by agent $i$ is minimized by $x_i^*$, and $x_i^*$ is likely to be different than that obtained for agent $j$. 
    This is equivalent to independent local training of the models.
    % If $\lambda = 0$, the problem \eqref{rewrite_main_problem} turns to the  problem with probably different arguments $x_i$'s minimizing the local functions $f_i$'s
    % \begin{equation*}
    %     \min_{\x \in \R^{nd}}\sum^n_{i = 1} f_i(x_i).
    % \end{equation*}
    % This is equivalent to independent local training of the models.
    % In this case, the 
    % This means that each local model $f_k$ have own variables that can be different. Another words, we give preference to each local model instead of global. 
    % \item When  $\lambda \in \left(0, +\infty\right)$, we have mixed model, where we take attention to both the local models and the global models. Easy to spot, that with increase of parameter $\lambda$, we give more preference to global model.
   \item
   As $\lambda \to +\infty$,  \eqref{rewrite_main_problem} tends to the distributed problem where the local arguments are constrained to be equal: i.e., $\min_{x_1 = \dots = x_n \in \R^d}\sum^n_{i = 1} f_i(x_i)$.
    % \begin{equation*}
    %     \min_{x_1 = \dots = x_n \in \R^d}\sum^n_{i = 1} f_i(x_i).
    % \end{equation*}
    This is equivalent to problem \eqref{main_problem0} and the training of one global model.
\end{itemize}

\subsection{Preliminaries}

% \subsection{Notation and Definitions}\label{sec:notation}

% Bold variables $\x$ denote vectors of length $\R^{nd}$, representing a concatenation of the vectors $\{x_i \}_{i=1}^n$, where $x_i \in \R^d$. 
% We use $\la x,y \ra = x^T y$ to denote the inner product of $x,y\in\R^d$, and the corresponding  definition of the $\ell_2$-norm in $\R^d$ is $\|x\|_2 = \sqrt{\la x, x \ra}$. The inner product and $\ell_2$-norm in $\R^{nd}$ are defined analogously.

%\subsection{Assumptions} \label{sec:as}

Throughout this work the following assumption is made regarding the functions in \eqref{rewrite_main_problem}.
% \begin{assumption} \label{as:func}
% It is assumed that each $f_i$ in problem \eqref{rewrite_main_problem} is
%  \begin{itemize}[label=$\ast$]
% \item  $L$-smooth w.r.t the $\ell_2$-norm, i.e. for all $u, v \in \R^d$ we have
% \[
% \left\| \nabla f_i(u)-\nabla f_i(v)
% \right\|_2  \leq L\left\|u - v\right\|_2;
% \]
%     \item $\mu$-strongly-convex w.r.t. the $\ell_2$-norm, i.e. for all $u, v \in \R^d$ we have
% \[
% f_i(u) - f_i(v) \geq \langle\nabla f_i(v) , u - v\rangle + \frac{\mu}{2}\left\|u- v\right\|^2_2. 
% \]
% \end{itemize}
% \end{assumption}
\begin{assumption} \label{as:func}
It is assumed that each $f_i:\R^d\to\R$ in problem \eqref{rewrite_main_problem} is
 \begin{itemize}[label=$\ast$]
\item  $L$-smooth w.r.t the $\ell_2$-norm, i.e. for all $u, v \in \R^d$,
$\| \nabla f_i(u)-\nabla f_i(v)
\|_2  \leq L\|u - v\|_2$; and
    \item $\mu$-strongly-convex w.r.t. the $\ell_2$-norm, i.e. $\forall u, v \in \R^d$,
$f_i(u) - f_i(v) \geq \langle\nabla f_i(v) , u - v\rangle + \tfrac{\mu}{2}\left\|u- v\right\|^2_2$. 
\end{itemize}
\end{assumption}
By Assumption~\ref{as:func}, $f$ in \eqref{rewrite_main_problem} is $L$-smooth and $\mu$-strongly convex, and subsequently $F$ is also $\mu$-strongly convex.

\subsection{Communication} \label{sec:matrix}

The communication network  is modeled as a fixed, connected, undirected graph,   $\mathcal{G} = (\mathcal{V}, \mathcal{E})$, where $\mathcal{V} = \{1,\ldots, n\}$ are vertices (devices) and  $\mathcal{E} = \{(i,j) \, |\, i,j \in \mathcal{V} \}$ are edges (connections between devices). Note that  $(i,j) \in \mathcal{E}$ if and only if there exists a communication link between agents $i$ and $j$. For such a graph, a gossip matrix $\hat W$ is defined as follows.
% For such a graph we can define  the gossip matrix $\hat W$, which  satisfies the following three assumptions:
% \begin{enumerate}
%     \item $\hat W$ is symmetric positive semi-definite;
%     \item The kernel $\hat W$ consists of the vector $\mathbf{1} = (1, \dots, 1)^{\top}$;
%     \item $\hat W$ is defined on the edges of the communication network: $\hat w_{i,j} \neq 0$ if and only if $i = j$ or $(i, j)\in E$.
% \end{enumerate}
\begin{definition}[Gossip matrix]\label{gossipmatrix}
A matrix $\hat W \in \R^{d\times d}$, associated with a graph $\mathcal{G}$, is called a gossip matrix, if it satisfies the following conditions: 
\begin{enumerate}
    \item $\hat W$ is symmetric positive semi-definite;
    \item The kernel of $\hat W$ consists of the vector $\mathbf{1} = (1, \dots, 1)^{\top}$;
    \item $\hat W$ is defined on the edges of the communication network: $\hat w_{i,j} \neq 0$ if and only if $i = j$ or $(i, j)\in \mathcal{E}$.
\end{enumerate}
\end{definition}
% In turn, the communication matrix $W$ in \eqref{rewrite_main_problem} is defined  as the Kronecker product of matrix $\hat W$ and the identity matrix $I_d$ to take into the consideration that all $x_k \in \R^d$ ($k=1,\dots,n$):  $W = \hat{W}\otimes I_d$. As mentioned above, only neighboring agents can communicate.  In the decentralized optimization this  communications is made by gossip protocol \cite{boyd2006randomized,nedic2009distributed}. Is is realized as a  matrix-vector multiplication by a matrix $W$ and vector $\bar x$. 

The communication matrix $W$ in \eqref{rewrite_main_problem} is $W = \hat{W}\otimes I_d$, i.e., $W$ is the Kronecker product of a gossip matrix $\hat W$ and the identity matrix $I_d$. \abc{Because only neighboring agents can communicate in this decentralized optimization setting, it is assumed that communication is made via a gossip protocol \cite{boyd2006randomized,nedic2009distributed}, i.e., $\hat W$ is a gossip matrix (Definition~\ref{gossipmatrix}) and communication is realized via matrix-vector multiplication with $W$.} \ab{During one communication/communication round, for every node, full local vectors of dimension $d$ (e.g. variables $\{ x_i\}$ or gradients $\{\nabla f_i(x_i) \}$) are exchanged with all neighbours. This work supposes that the network remains unchanged, all connections are stable, and no interruptions nor asynchronous/delayed transmissions are considered.}

Here, $\eigWmax$ denotes the maximum eigenvalue of $W$, $\eigWmin$  denotes the minimum positive eigenvalue of $W$ and $\chi \geq \eigWmax / \eigWmin$ is an upper bound on the condition number. Because $W = \hat{W}\otimes I_d$, it holds that $\eigWmax = \eigWhatmax$ and $\eigWmin = \eigWhatmin$. The quantity $\chi$ reflects how quickly information is transmitted through the graph; a small $\chi$ corresponds to fast transmission, while a large $\chi$ corresponds to slow transmission. 

\ab{A simple example of a matrix $\hat W$ satisfying Definition~\ref{gossipmatrix} is the 
Laplacian matrix. For example, the Laplacian of a linear graph (chain) is}
% \begin{eqnarray*}
% \ab{\hat W = \left(
% \begin{array}{cccccccc}
% 1& -1 & & & & & &  \\
% -1 & 2 &-1 & & & & &  \\
% &-1 &2 &-1 & & & & \\
% & &-1 &2 &-1 & & & \\
% & & & &\ldots &\ldots & & \\
% & & & &-1 &2  &-1   & \\
% & & &   & &-1 &2 &-1 \\
% & & &  & & &-1 &1 \\
% \end{array}
% \right).}
% \end{eqnarray*}
\begin{eqnarray*}
\ab{\hat W = \left(
\begin{array}{cccccccc}
1& -1 & & & & &   \\
-1 & 2 &-1 & & & &   \\
&-1 &2 &-1 & & &  \\
 & & \ddots&\ddots &\ddots & & \\
 & & &-1 &2  &-1   & \\
 & &   & &-1 &2 &-1 \\
 & &  & & &-1 &1 \\
\end{array}
\right).}
\end{eqnarray*}
\ab{In terms of personalization, this means that the model on the first node relies directly on the 2nd node. In turn the 2nd node depends on the 1st and 3rd nodes, and so on. In particular, the 1st and last nodes depend on each other weakly and only indirectly through the whole chain. 

However, it is also possible to define $\hat W$ in a more complex way. For example, in the case of a linear graph, one can add weights that represent how much a given node relies upon its neighbours:}  
\begin{eqnarray*}
\ab{\hat W = \left(
\begin{array}{cccccccc}
1& -1 & & & & & &  \\
-1 & 1.5 &-0.5 & & & & &  \\
&-0.5 &1 &-0.5 & & & & \\
& &-0.5 &1.2 &-0.7 & & & \\
& & & \ddots&\ddots &\ddots & & \\
& & & &-0.5 &1.5  &-1   & \\
& & &   & &-1 &1 &-1 \\
& & &  & & &-1 &1 \\
\end{array}
\right).}
\end{eqnarray*}
\ab{In this example, the second node trusts the 1st node more than the 3rd, while the third node trusts the 2nd and 4th nodes equally, etc.}

\section{Contributions}

In this paper, we study the personalized federated learning formulation \eqref{rewrite_main_problem}. Lower complexity bounds for communication and local computation are proposed, and we develop several algorithms capable of achieving the lower bounds. Our results extend the work in \cite{hanzely2020lower}, which used the penalty \eqref{cent_penalt}, to problem \eqref{rewrite_main_problem}, which involves a penalty more amenable to the decentralized setting. Our contributions are summarized now.

\begin{itemize}
    \item \emph{Lower bounds.} We present lower bounds for the decentralized personalized federated learning problem \eqref{rewrite_main_problem} in the deterministic case (i.e., when we have access to full gradients for each $f_i$); see Section \ref{sec:lower}. The lower bounds are valid for all values of the parameter $\lambda$. In particular, in the smooth strongly convex case with small $\lambda$, the lower bounds are of the order $\sqrt{\lambda\lambda_{\max}\left(W\right)/\mu}$, which can be a significant improvement on the bound $\sqrt{\chi L/\mu}$ in the general, non-personalized case, \cite{scaman2017optimal}. 
This reflects a key advantage of the formulation \eqref{rewrite_main_problem}, because it is then  possible to both solve the problem of personalizing the models, and also to significantly reduce the total number of communications. This is an important factor not only in federated learning, but also in general distributed learning. Note that the lower bounds obtained in the work \cite{hanzely2020lower} are a special case of our lower bounds, when the communication network is represented by a fully connected graph. A summary of these lower bounds is presented in Table~\ref{tab:summary}.
\renewcommand{\arraystretch}{1.7}
\renewcommand{\tabcolsep}{6pt} 
\begin{table}[h!]
\vspace{-0.3cm}
\begin{center}
\begin{tabular}{c|c|c}
\hline
& \textbf{Lower bounds} & \textbf{Upper bounds}  \\ \hline
{\tt comm}& $\widetilde{\Omega}\left(\min\left\{ \sqrt{\frac{\lambda\eigWmax}{\mu}}, \sqrt{\frac{L}{\mu}\chi}\right\}\right)$ & $\widetilde{\mathcal{O}}\left(\min\left\{\sqrt{\frac{\lambda\eigWmax}{\mu}}, \sqrt{\frac{L}{\mu} \chi}\right\}\right)$  \\ \hline
{\tt local}& $\widetilde{\Omega}\left(\sqrt{\frac{L}{\mu}}\right)$ & $\widetilde{\mathcal{O}}\left(\sqrt{\frac{L}{\mu}}\right)$ \\ \hline
\end{tabular}
\vspace{0.3cm}
\caption{\small Summary of complexity results (upper and lower bounds) on communications ({\tt comm}) and local computations ({\tt local}) for finding an $\varepsilon$-optimal solution of \eqref{rewrite_main_problem} in the deterministic (gradient) case. 
}
\label{tab:summary}
\end{center}
\vskip-1cm
\end{table}

\item \emph{Near-optimal algorithm.} Another contribution is the development of optimal algorithms that match the theoretical lower bounds. The Accelerated Meta-Algorithm of \cite{gasnikov2020AM} (for general composite problems), is used as the base algorithm. The application of this algorithm to our problem formulation \eqref{rewrite_main_problem} is discussed, and specific implementation modes are suggested depending on small and large values of the regularization parameter $\lambda$. The analysis of the convergence in these modes shows that using this approach we achieve the lower optimal bounds up to logarithmic factors (Section \ref{sec:analys}). Hence, our algorithm is `near-optimal' in the deterministic case; see Section~\ref{sec:algorithms}.

\item \emph{Stochastic case.} We extend the previously reported results from the deterministic case (when the full gradient for all $f_i$ is available), to the stochastic setting. In particular, we consider the case when each local function $f_i$ is a finite sum (for example, the sum of batches), i.e. $f_i = \tfrac{1}{M} \sum_{m=1}^{M} f_{i,m}$. In this case, for one call of the oracle we can get only the gradient of one term $f_{i,m}$. We provide lower bounds, as well as a stochastic modification of our near-optimal deterministic algorithm; see Section \ref{sec:stoch}.

\item \emph{Experiments.} We present numerical experiments to demonstrate the benefits of our approach. In particular, we used several datasets from the benchmark LIBSVM library, and we considered several different graph structures. We also run the experiments for several values of the penalty parameter $\lambda$, to better  understand the impact of personalization; see Section~\ref{sec:experiments}.
\end{itemize}

\section{Lower Bounds} \label{sec:lower}

In this section, optimal algorithms for problems of the form \eqref{rewrite_main_problem} are described, and lower bounds on the local computation and communication costs for such optimal algorithms, are presented.
% In this section we present lower bounds on problem
% \eqref{rewrite_main_problem}. But before that, we restrict the class of algorithms for which the following lower bounds are valid.
We begin with the following assumption, which describes the properties of algorithms relevant for this work, (i.e., the properties of the algorithms for which the lower bounds, developed later in this section, are valid).  Such an assumption is common in the literature; see, for example, \cite{hanzely2020lower,Hendrikx2020,scaman2017optimal}.
% \begin{assumption}[Class of Optimal Algorithms]\label{as:lower}
% Let $\{\x^k\}^{N}_{k = 1}$ be iterates generated by algorithm $\cA$. For each node of graph $\mathcal{G}$ we define sequence of local memory $\{\cM_{i, k}\}^{N}_{k = 1}$ for $1 \leq i \leq n$:
% \begin{eqnarray*}
% \cM_{i, 0} &=&\text{span}\left\{x^0_i\right\},\\
% \cM_{i, k+1} &=&\begin{cases}
%  \text{span}\left\{\cM_{i, k}, \nabla f_i(y_i)\right\}, \forall y_i \in \cM_{i, k}& \text{if local comp.  at iteration }k\\
%  \text{span}\left\{\bigcup\limits_{j: (i,j) \in\mathcal{E} }\cM_{j, k}\right\}& \text{if communication at iteration } k.
% \end{cases}
% \end{eqnarray*}
% \end{assumption}
% \begin{assumption}\label{as:lower}
% Let $\{\x^k\}^{N}_{k = 1}$ denote the iterates generated by an Algorithm~$\cA$, belonging to the class of optimal algorithms. For each node of graph $\mathcal{G}$ we define sequence of local memory $\{\cM_{i, k}\}^{N}_{k = 1}$ for $1 \leq i \leq n$:
% \begin{eqnarray*}
% \cM_{i, 0} &=& {\rm span}\left\{x^0_i\right\},\\
% \cM_{i, k+1} &=&\begin{cases}
%  {\rm span}\left\{\cM_{i, k}, \nabla f_i(y_i)\right\}, \forall y_i \in \cM_{i, k}& \text{if local comp.  at iteration }k\\
%  {\rm span}\left\{\bigcup\limits_{j: (i,j) \in\mathcal{E} }\cM_{j, k}\right\}& \text{if communication at iteration } k.
% \end{cases}
% \end{eqnarray*}
% \end{assumption}
\begin{assumption}\label{as:lower}
Consider an Algorithm~$\cA$, for problem \eqref{rewrite_main_problem}. Then, the iterates $\{\x^k\}^{K}_{k = 1}$ of Algorithm~$\cA$ are generated using only components available in local memory, where, for each node of graph $\mathcal{G}$ the sequence of local memory $\{\cM_{i, k}\}^{K}_{k = 1}$ for $1 \leq i \leq n$ is:
\begin{eqnarray*}
\cM_{i, 0} &=&  \left\{x^0_i\right\},\\
\cM_{i, k+1} &=&\begin{cases}
 {\rm span}\left\{\cM_{i, k}, \nabla f_i(y_i)\right\}, \forall y_i \in \cM_{i, k}& \text{if local comp.  at iteration }k\\
 {\rm span}\left\{\bigcup\limits_{j: (i,j) \in\mathcal{E} }\cM_{j, k}\right\}& \text{if communication at iteration } k.
\end{cases}
\end{eqnarray*}
\end{assumption}
% \todo{Clarify whether $\cM_{i, 0} = {\rm span}\left\{x^0_i\right\}$ or $\cM_{i, 0} = {\rm span}\left\{x^0_1,\dots,x^0_n\right\}$}
Assumption~\ref{as:lower} can be interpreted as follows. Initially, each agent $i$ (corresponding to a node on graph $\mathcal{G}$) has local memory $\cM_{i, 0}$, which comprises of the initial point $x_i^0$. At any iteration $k\geq 1$, the algorithm can either perform a computation using the locally available memory, or it can carry out a communication step. If the algorithm performs a local computation, then each device can calculate the gradient at any point from its current memory $\cM_{i, k}$ and take a linear combination of this gradient with the previously generated points stored in $\cM_{i, k}$. If the algorithm performs a communication step, then information is exchanged with neighbors and the current local memory $\cM_{i, k}$ is combined with that held by its neighbors. Such an algorithm is first order, because it generates its iterates using linear combinations of local points and gradients.
%  \begin{itemize}[label=$\ast$]
% \item Initially, $\cM_{i, 0}$ are simply all possible linear combinations of the initial points $x^0_i$ of all agents. 

% \item Further, at each iteration $k$, the algorithm can make either local calculations or communication. 

% \item If the algorithm choose local computations, then each device can calculate the gradient at any point from the current $\cM_{i, k}$ and add a linear combination of this gradient to the already reached points $\cM_{i, k}$. 

% \item If the algorithm does communication, then it exchanges information with a neighbor and combines sets $\cM_{i, k}$ with neighbors. 
% \end{itemize}

We are now ready to present our first theorem, which gives a lower bound on the number of communications needed by an algorithm whose iterates are generated according to Assumption~\ref{as:lower}. (The proof can be found in Appendix \ref{ap:T1}.)
% \begin{theorem}
% \label{theorem_lower_bounds}
% Let $\sqrt{\chi} \geq 6$, $L \geq \mu$, and $\eigWminplus \geq \mu$. Then, there exists a graph $\mathcal{G}$, with a corresponding matrix $\hat W$ satisfying Definition~\ref{gossipmatrix}, a set of $L$-smooth $\mu$-strongly convex functions $f_1, f_2, \dots, f_n: \R^d \rightarrow \R$ and an initial point $\x^0 \in \R^{nd}$, such that after $N$ iterations the output $\x^N \in \bigcup_i \cM_{i, N}$ of any algorithm $\cA$ satisfying Assumption \ref{as:lower} is
% \begin{equation*}
% \|\x^{N} - \x^*\|^2 \geq
% \left(1 - 10\max\left\{ \sqrt{\frac{\mu}{(L -\mu) \chi}},\sqrt{\frac{\mu}{\lambda \lambda_{\max}(\hat W)}}\right\}\right)^{q}\frac{\|\x^{0} - \x^*\|^2}{4},
% \end{equation*}
% where $q$ is number of communication rounds among $N$ iterations.
% \end{theorem}
\begin{theorem}
\label{theorem_lower_bounds}
Let $\chi \geq 3$, $L \geq 2\mu$, and $\lambda\eigWmin \geq \mu$. Then there exist functions $f_1, f_2, \dots, f_n: \R^d \rightarrow \R$ satisfying Assumption~\ref{as:func}, a graph $\mathcal{G}$ with associated matrix $\hat W$ satisfying Definition~\ref{gossipmatrix}, and an initial point $\x^0 = [(x^0_1)^T, \ldots (x^0_n)^T]^T \in \R^{nd}$, such that 
any algorithm $\cA$ (satisfying Assumption \ref{as:lower}) among $K$ iterations need to make at least
\begin{align*}
     \Omega\left( \min\left\{\sqrt{\frac{\lambda \lambda_{\max}(W)}{\mu }}, \sqrt{\frac{(L - \mu)\chi}{\mu}}\right\}  \log\frac{1}{\varepsilon}\right) \quad \text{communications}
\end{align*}
to achieve $\varepsilon$-optimal solution in the outputs ($\| x_j - x^*_j\|^2 \varepsilon$ for all $j$).
\end{theorem}
The proof of this Theorem is placed in Appendix \ref{ap:T1}. 

It remains to develop lower bounds for the local computation costs for any algorithm $\cA$ satisfying Assumption~\ref{as:lower}. Hence, consider a special instance of problem \eqref{rewrite_main_problem}, where $\x^0\in \R^{nd}$, $f_1 = f_2 = \cdots = f_n$, and $\hat W$ is the Laplace matrix for a fully connected graph. Then \eqref{rewrite_main_problem} reduces to the minimization of the single local function $f_1$ (communication is unnecessary, irrespective of $\lambda$, because the functions are all identical). Now, if $f_1$ is chosen to be the worst-case quadratic from \cite{nesterov2003introductory}, then the lower bound of at least
\begin{eqnarray}\label{temp1}
    N^{\text{loc}} =
\Omega\left(\sqrt{\tfrac{L}{\mu}}\log \tfrac1{\varepsilon}\right)
\end{eqnarray} gradient calls are needed to find an $\varepsilon$-optimal solution. 

% Theorem \ref{theorem_lower_bounds} and the local computation cost can be combined into the following result, which gives the complexity bounds for the communication and local computation costs incurred by any algorithm $\cA$ satisfying Assumption~\ref{as:lower} to obtain an $\varepsilon$-optimal solution to problem \eqref{rewrite_main_problem}.
% \begin{corollary}
% \label{cor_lower_bounds}
% Let $\sqrt{\chi} \geq 6$, $L \geq 2\mu$, $\lambda\lambda^+_{\min}(W) \geq \mu$. Then there exist functions $f_1, f_2, \dots, f_n: \R^d \rightarrow \R$ satisfying Assumption~\ref{as:func}, a graph $\mathcal{G}$ with associated matrix $\hat W$ satisfying Definition~\ref{gossipmatrix}, and an initial point $\x^0 \in \R^{nd}$, such that any algorithm $\cA$ satisfying Assumption~\ref{as:lower} needs to make a least 
% % \begin{equation*}
% % \Omega\left(\min\left\{ \sqrt{\frac{\lambda\lambda_{\max}\left(\hat W\right)}{\mu}}, \sqrt{\frac{L}{\mu}\chi}\right\} \log \frac{1}{\varepsilon}\right)~~\text{comm. and} ~~~ \Omega\left(\sqrt{\frac{L}{\mu}}\log \frac{1}{\varepsilon}\right) ~~\text{local comp.}
% % \end{equation*}
% \begin{equation*}
% \Omega\left(\min\left\{ \sqrt{\tfrac{\lambda\lambda_{\max}( W)}{\mu}}, \sqrt{\tfrac{L}{\mu}\chi}\right\} \log \tfrac{1}{\varepsilon}\right)~~\text{comm. and} ~~~ \Omega\left(\sqrt{\tfrac{L}{\mu}}\log \tfrac{1}{\varepsilon}\right) ~~\text{local comp.}
% \end{equation*}
% to obtain an output $\x^K$ such that $\|\x^{K} - \x^*\|^2 \leq \varepsilon$, where $\x^*$ denotes the optimal solution to problem \eqref{rewrite_main_problem}.
% \end{corollary}

\section{Algorithms}\label{sec:algorithms}

The goal of this section is to develop an optimal algorithm for problem \eqref{rewrite_main_problem}, i.e., to develop an algorithm whose iterates satisfy the lower bounds in Section~\ref{sec:lower}. In Section \ref{sec:ama}, we discuss an algorithm that can be applied to general composite optimization problems. In Section \ref{sec:analys}, this algorithm is specialized to the application considered in this work, that of decentralized personalized federated learning \eqref{rewrite_main_problem}. The algorithms in Sections~\ref{sec:ama} and  \ref{sec:analys} can be applied to deterministic problems, and the extension to a stochastic setting is considered in Section \ref{sec:stoch}. In particular, the case when the function at each node has finite sum structure is considered, and two approaches, both equipped with convergence results, are described and compared.

\subsection{Accelerated Meta-Algorithm} \label{sec:ama}

In this section, consider the general composite optimization problem
\begin{equation}
\label{composite_problm}
    \min_{\x \in \mathbb{R}^{nd}} H(\x) = h_1(\x)+h_2(\x).
\end{equation}
% \todo{Is $H\equiv F$?}
(Later we will consider how $h_1$ and $h_2$ in  \eqref{composite_problm} are related to $f$ and $g $ in \eqref{rewrite_main_problem}.) The following assumption is made about problem \eqref{composite_problm}.
\begin{assumption}\label{assumptionh1h2}
For the problem \eqref{composite_problm}, it is assumed that $h_1$ is convex and $L(h_1)$-smooth, that $h_2$ is convex and $L(h_2)$-smooth, and that $H$ is $\mu$-strongly convex.
\end{assumption}

There are many efficient algorithms that can be applied to problem \eqref{composite_problm}, including the Accelerated-Meta-Algorithm (see Algorithm \ref{meta_algorithm}) proposed in \cite{gasnikov2020AM}, as well as its restarted version (see Algorithm \ref{restart_AM_algorithm}).

\begin{algorithm} [h!]
\caption{Accelerated Meta-Algorithm (\texttt{MA}) \cite{gasnikov2020AM}.}
\label{meta_algorithm}
\hspace*{0.2cm} {\bf Input:} starting point $\x^0 \in \R^{nd}$, no. of iterations $K$,  parameter $\gamma>0$, accuracy $\delta > 0$\\
\hspace*{0.3cm}{\bf Initialization:} $A^0 = 0,\quad \y^0 = \x^0, \quad \tau = \frac{1}{2\gamma}$
	\begin{algorithmic}
		\For{ $k = 0, \ldots, K-1$}
		\State $ a^{k+1} = \displaystyle \frac{\tau+\sqrt{\tau^2+4\tau A^k}}{2}$
		\State $A^{k+1} = A^k+a^{k+1}$
		\State $\w^k = \displaystyle \frac{A^k}{A^{k + 1}}\y^k + \frac{a^{k+1}}{A^{k+1}} \x^k
		$
		\State Find $\y^{k+1}\in \R^{nd}$, such that $\|\hat \y^{k+1} -  \y^{k+1}\|_2^2 \leq \delta$, where
		\begin{equation}\label{auxiliary_prob}
		    \hat \y^{k+1} = \argmin_{\y\in \R^{nd}} \left\{ \langle \nabla h_1 (\w^k), \y-\w^k\rangle  + h_2(\y)   +\frac{\gamma}{2}\|\y-\w^k\|_2^{2} \right\}
		\end{equation}
		\State $\x^{k+1} = \x^k-a^{k+1} \nabla H(\y^{k+1})$
		\EndFor
	\end{algorithmic}
	\hspace*{0.2cm} {\bf Output:} $\y^{K}$
\end{algorithm}
% \begin{algorithm} [h!]
% \caption{Restarted Accelerated Meta-Algorithm (\texttt{Restarted-MA}) \cite{gasnikov2020AM}}
% \label{restart_AM_algorithm}
% \hspace*{0.2cm} {\bf Input:} initial point $\x_0 \in \R^{nd}$, no. of iterations $S$,  parameter $\gamma>0$, accuracy $\delta > 0$\\
% \hspace*{0.3cm}{\bf Initialization:} $\displaystyle N_s= \max \left\{\left\lceil4\cdot \sqrt{\tfrac{2\gamma}{\mu}} \right\rceil, 1 \right\}$
% 	\begin{algorithmic}
% 		\For{ $s = 0, \ldots, S-1$}
% 		\State  $\x_{s+1} =$ \texttt{AM}$(\x_s,N_s, \gamma, \delta)$
% 		\EndFor
% 	\end{algorithmic}
% 	\hspace*{0.2cm} {\bf Output:} $\x_S$
% \end{algorithm}
\begin{algorithm} [h!]
\caption{Restarted Accelerated Meta-Algorithm (\texttt{Restarted-MA}) \cite{gasnikov2020AM}}
\label{restart_AM_algorithm}
\hspace*{0.2cm} {\bf Input:} initial point $\x^0 \in \R^{nd}$, no. of iterations $S$,  parameter $\gamma>0$, accuracy $\delta > 0$\\
\hspace*{0.3cm}{\bf Initialization:} $\displaystyle N^s= \max \left\{\left\lceil4\cdot \sqrt{\tfrac{2\gamma}{\mu}} \right\rceil, 1 \right\}$
	\begin{algorithmic}
		\For{ $s = 0, \ldots, S-1$}
		\State  $\x^{s+1} =$ \texttt{MA}$(\x^s,N^s, \gamma, \delta)$
		\EndFor
	\end{algorithmic}
	\hspace*{0.2cm} {\bf Output:} $\x^S$
\end{algorithm}

\begin{remark}
Note that \eqref{auxiliary_prob} can be solved, for example, by Accelerated Gradient Descent \cite{nesterov2003introductory}.
\end{remark}

The following convergence results hold for \texttt{Restarted-MA} (Algorithm~\ref{restart_AM_algorithm}) applied to  problem \eqref{composite_problm}.
% \begin{theorem}[\cite{gasnikov2020AM}]
% \label{theorem_meta_algorithm}
%     Let $\x_S$ denote the output of Algorithm \ref{restart_AM_algorithm} after $S$ iterations. If  we choose $\gamma \geq 2L(h_1)$, then to achieve $F(\x_S) - F(x_*) \leq \varepsilon$, we need to make
%     \begin{equation}
%     \label{eq:MA_conv}
%     S = {\mathcal{O}} \left(\sqrt{\tfrac{L(h_1)}{\mu}}\log\tfrac{1}{\varepsilon} \right) \quad \text{iterations,}
%     \end{equation}
%     and also to choose the accuracy $\delta$ of \ref{auxiliary_prob}, satisfies to the following inequality
%     \begin{equation*}
%         \delta \leq \frac{\e\mu}{864^2(L(h_1)+ L(h_2) +\gamma)^2}.
%     \end{equation*}
% \end{theorem}
\begin{theorem}[Theorem 3 in \cite{tominin2021accelerated}]
\label{theorem_meta_algorithm}
    Let Assumption~\ref{assumptionh1h2} hold, let $\gamma \geq 2L(h_1)$, let $\varepsilon >0$, and let 
    \begin{equation*}
        \delta \leq \frac{\e\mu}{864^2(L(h_1)+ L(h_2) +\gamma)^2}.
    \end{equation*}
    If Algorithm \ref{restart_AM_algorithm} runs for
    \begin{equation}
    \label{eq:MA_conv}
    S = {\mathcal{O}} \left(\sqrt{\tfrac{L(h_1)}{\mu}}\log\tfrac{1}{\varepsilon} \right)
    \end{equation}
    iterations, generating output $\x^S$, then $H(\x^S) - H(\x^*) \leq \varepsilon$, where $\x^*$ denotes the optimal solution to \eqref{composite_problm}.
\end{theorem}

%%%%%%%%%%%%%%%%%%%%%%%%%%%%%%%%%%%%%%%%%%%%%%%%%%%%%%%%%%%%%
%\section{Main Results. Near-Optimal Algorithm}

\subsection{Convergence Analysis -- Near-Optimal Algorithm} \label{sec:analys}

Section~\ref{sec:ama} introduced an accelerated algorithm for the general problem \eqref{composite_problm}, with associated convergence results. The purpose of this section is to make a connection between the results in Section~\ref{sec:ama}, and how they are applicable in the context of personalized federated learning (i.e., problem \eqref{rewrite_main_problem}). Moreover, the lower bounds established in Section~\ref{sec:lower} related to local computation and communication costs for an optimal algorithm for problem \eqref{rewrite_main_problem}. Thus, another goal is to show that Algorithms~\ref{meta_algorithm}+\ref{restart_AM_algorithm} is an optimal algorithm for \eqref{rewrite_main_problem}, by showing that it achieves the lower bounds on communication and local computation costs presented in Section~\ref{sec:lower}.

\abc{By comparing problems \eqref{rewrite_main_problem} and \eqref{composite_problm}, it can be seen that they are both convex and composite. The key here is that we do not make a one-to-one correspondence between $(f,g)$ and $(h_1,h_2)$. That is, depending on the parameter $\lambda$, two different cases --- one in which $f\equiv h_1$ and $g\equiv h_2$, while the other in which $f\equiv h_2$ and $g\equiv h_1$ --- are considered.} \ab{Practical versions of Algorithm \ref{meta_algorithm} for the problem \eqref{rewrite_main_problem} are presented in Appendix~\ref{app:alg_prac} (Algorithms \ref{meta_algorithm_big_lambda} and \ref{meta_algorithm_small_lambda}).}

Regardless, to apply Algorithm~\ref{meta_algorithm}+\ref{restart_AM_algorithm}, it is necessary to compute the gradients for both $h_1$ and $h_2$ (recall subproblem \eqref{auxiliary_prob}), and therefore for both $f$ and $g$ when extending to the original problem \eqref{rewrite_main_problem}. So, let us study how to compute the gradients $\nabla f$ and $\nabla g$ for \eqref{rewrite_main_problem}, and try to understand where the communications arise. Note that the computation of $\nabla f$ does not require communication. Indeed, each block $i$ has a corresponding gradient $\nabla f_i$ (taken with respect to the local variables $x_i$), and the `long' gradient $\nabla f$ is simply the concatenation of the block gradients. On the other hand, $\nabla g(\x) = \lambda W \x$, and to compute the matrix-vector product $W\x$ requires communication with neighbors (recall the gossip protocol described in Section \ref{sec:matrix}, and see also \cite{boyd2006randomized,nedic2009distributed}). It can be shown that computing $\nabla g(\x)$ is equivalent to one communication. Therefore, if we know how many times $\nabla f(\x)$ and $\nabla g(\x)$ are called by Algorithm~\ref{meta_algorithm}, then  complexities for the number of local computations and communications, respectively, can be obtained.

We are now ready to present the main convergence theorem of this paper, which provides complexity results for the local computation and communication costs for an optimal algorithm for problem \eqref{rewrite_main_problem}.
\begin{theorem}\label{th:alg_pfl_minimization}
Let Assumption~\ref{as:func} hold and let the graph $\mathcal{G}$ have an associated matrix $\hat W$ that satisfies Definition~\ref{gossipmatrix}. Then, to obtain an $\varepsilon$-optimal solution to problem \eqref{rewrite_main_problem}, solving by Algorithm \ref{restart_AM_algorithm} with
\begin{equation}
    \label{eq:delta}
    \delta = \frac{\e\mu}{3000^2(L + \lambda\lambda_{\max}(W))^2},
\end{equation}
requires the number of communications and local computations, respectively, to be of the order
\begin{equation}\label{comms}
    N^{\text{comm}} = \mathcal{O}\left(\min\left\{\sqrt{\tfrac{\lambda\eigWmax}{\mu}}, \sqrt{\tfrac{L}{\mu} \chi}\right\}\log\tfrac{1}{\varepsilon}\log\tfrac{1}{\delta}\right),
\end{equation}
and 
\begin{equation}\label{loc}
    N^{\text{loc}} = \mathcal{O}\left(\sqrt{\tfrac{L}{\mu}}\log\tfrac{1}{\e}\log\tfrac{1}{\delta}\right).
\end{equation}
\end{theorem}

\begin{proof} First, note that $\mathcal{G}$ is a quadratic function with a positive semi-definite Hessian, so it is $\eigWmax$-smooth and convex. Moreover, it is $\lambda\lambda^+_{\min}(W)$-strongly convex on the subspace $\left(\textbf{Ker } W\right)^{\perp}$. By Assumption~\ref{as:func}, $f$ is $L$-smooth and $\mu$-strongly convex. Hence, $F$ is strongly convex. Thus, the conditions of Theorem~\ref{theorem_meta_algorithm} hold, and the application of its analysis is valid. The remainder of the analysis is split into two cases.

\textbf{Case 1:} $\lambda\eigWmax \geq L$. 
Here, let $h_1(\x) = f(\x)$ and $h_2(\x) = g(\x)$. Theorem \ref{theorem_meta_algorithm} gives the complexity for the function $h_1 = f$, i.e. the number of local computations $N^{\text{loc}}$ is given in \eqref{loc}.
% \begin{equation}
% \label{temp1}
%     N^{\text{loc}} = \mathcal{O}\left(\sqrt{\tfrac{L}{\mu}}\log\tfrac{1}{\e}\right).
% \end{equation}
Also, it can be shown that $\delta$ in \eqref{eq:delta} satisfies the condition in Theorem~\ref{theorem_meta_algorithm}.

\abc{Next, consider the auxiliary problem \eqref{auxiliary_prob}. By Definition~\ref{gossipmatrix}, $\textbf{Ker } W$ is not empty, and the function $g(\x)$ takes a zero on this subspace. Then we can divide our problem into two subproblems: minimization of a quadratic form with matrix $\gamma \cdot I$ on $\textbf{Ker } W$ and minimization of a quadratic form with matrix $\lambda W + \gamma \cdot I$ on $\left(\textbf{Ker } W\right)^{\perp}$. The complexity of the first problem is  $\mathcal{O}\left(1\right)$.} The second problem is $\lambda\lambda^+_{\min}(W)$-strongly convex, and if the Accelerated Gradient Method \cite{nesterov2003introductory} is used to solve this subproblem, then the complexity is 
\begin{equation}\label{complexityAGM}
    \mathcal{O}\left(\sqrt{\tfrac{\gamma + \lambda\eigWmax}{\max\{\gamma, \lambda \lambda^+_{\min}(W)\}}}\log\tfrac{1}{\delta}\right).
\end{equation}
This is the complexity for a single subproblem \eqref{auxiliary_prob} solve, but \eqref{auxiliary_prob} is solved \eqref{temp1} times. Overall, this means that the total number of calls of $\nabla g$ is:
\begin{equation*}
    N^{\text{comm}} = \mathcal{O}\left(\sqrt{\tfrac{L}{\mu}} \sqrt{\tfrac{\gamma + \lambda\eigWmax}{\max\{\gamma, \lambda \lambda^+_{\min}(W)\}}} \log\tfrac{1}{\e}\log\tfrac{1}{\delta}\right).
\end{equation*}
Noting that
\begin{equation*}
    \sqrt{\tfrac{\gamma + \lambda\eigWmax}{\max\{\gamma, \lambda \lambda^+_{\min}(W)\}}}  = \min\left\{ \sqrt{\tfrac{\gamma + \lambda \lambda_{\max}(W)}{\gamma}}, \sqrt{\tfrac{\gamma + \lambda \lambda_{\max}(W)}{\lambda \lambda^+_{\min}(W)}} \right\},
\end{equation*}
and taking $\gamma=2L$, gives \eqref{comms}.
% \begin{equation*}
%     N^{\text{comm}} = \mathcal{O}\left(\min\left\{\sqrt{\frac{\lambda \lambda_{\max}(W)}{\mu}}, \sqrt{\frac{L}{\mu} \frac{\lambda \lambda_{\max}(W)}{\lambda \lambda^+_{\min}(W)}}\right\}\log\frac{1}{\e}\log\frac{1}{\delta}\right).
% \end{equation*}

\textbf{Case 2:} $\lambda\eigWmax < L$.
Here, let $h_1 = g$ and $h_2 = f$. Theorem~\ref{theorem_meta_algorithm} gives the complexity for the function $h_1 = g$, i.e. the number of communications is
\begin{equation}
    \label{temp2}
    N^{\text{comm}} = \mathcal{O}\left(\sqrt{\tfrac{\lambda\eigWmax}{\mu}}\log\tfrac{1}{\e}\right) = \mathcal{O}\left(\min\left\{\sqrt{\tfrac{\lambda\eigWmax}{\mu}}, \sqrt{\tfrac{L}{\mu} \chi}\right\}\log\tfrac{1}{\varepsilon}\log\tfrac{1}{\delta}\right).
\end{equation}
In last step we additionally use that $\chi \geq 1$. Also, it can be shown that $\delta$ in \eqref{eq:delta} satisfies the condition in Theorem~\ref{theorem_meta_algorithm}.
If the Accelerated Gradient Method \cite{nesterov2003introductory} is used to solve subproblem \eqref{auxiliary_prob}, the complexity for a single subproblem solve is again given by \eqref{complexityAGM},
% \begin{equation*}
%     \mathcal{O}\left(\sqrt{\frac{\gamma + \lambda \lambda_{\max}(W)}{\max\{\gamma, \lambda \lambda^+_{\min}(W)\}}}\log\frac{1}{\delta}\right),
% \end{equation*} 
and this subproblem is solved \eqref{temp2} times. Then we can find the number of calls for $\nabla f$:
\begin{equation*}
    N^{\text{loc}} = \mathcal{O}\left(\sqrt{\tfrac{\lambda \lambda_{\max}(W)}{\mu}}\sqrt{\tfrac{L + \gamma}{\mu + \gamma}}\log\tfrac{1}{\e}\log\tfrac{1}{\delta}\right).
\end{equation*}
Taking $\gamma = 2\lambda\eigWmax$ gives \eqref{loc}.
% \begin{equation*}
%     N^{\text{loc}} = \mathcal{O}\left(\sqrt{\tfrac{L }{\mu }}\log\tfrac{1}{\e}\log\tfrac{1}{\delta}\right).
% \end{equation*}

Finally, combining the two cases establishes the theorem statement.\EndProof
% \begin{equation*}
%     N^{\text{comm}} = \mathcal{O}\left(\min\left\{\sqrt{\tfrac{\eigWmax}{\mu}}, \sqrt{\tfrac{L}{\mu} \tfrac{\eigWmax}{\lambda \lambda^+_{\min}(W)}}\right\} \log\tfrac{1}{\e}\log\tfrac{1}{\delta}\right).
% \end{equation*}
% and 
% \begin{equation*}
%     N^{\text{loc}} = \mathcal{O}\left(\sqrt{\frac{L}{\mu}}\log\frac{1}{\e}\log\frac{1}{\delta}\right).
% \end{equation*}

\end{proof}

\begin{remark}
Note that in the centralized case (with a completely connected communication network) we have that $\chi = 1,~ \lambda_{\max}(W) = 1$ and our method converges with the following rates: 
\begin{equation*}
    N^{comm} = \widetilde{\mathcal{O}}\left(\min\left\{\sqrt{\tfrac{\lambda}{\mu}}, \sqrt{\tfrac{L}{\mu}}\right\}\right),\quad N^{loc} = \widetilde{\mathcal{O}}\left(\sqrt{\tfrac{L}{\mu}}\right).
\end{equation*}
These bounds coincide with lower bounds for centralized PFL \cite{hanzely2020lower}.
\end{remark}

\subsection{Stochastic case} \label{sec:stoch}

% Here we extend the work previously presented and consider the stochastic case of the problem \eqref{rewrite_main_problem}. In particular, it is assumed that each local function has a sum structure (sum of batches):
% \begin{equation}
% \label{rewrite_main_problem_sum}
% f_i(x_i) = \frac{1}{M}\sum_{m=1}^Mf_{i,m}(x_i).
% \end{equation}
Here we extend the work previously presented and consider the stochastic case of problem \eqref{rewrite_main_problem}. In particular, it is assumed that each local function has a sum structure, so that \eqref{rewrite_main_problem} becomes
\begin{equation}
\label{rewrite_main_problem_sum}
    \min_{\x \in \R^{nd}} \sum^n_{i = 1} \underbrace{\frac{1}{M}\sum_{m=1}^M f_{i,m}(x_i)}_{f_i(x_i)} + \frac{\lambda}{2}\langle \x, W\x\rangle.
\end{equation}
\ab{This setup often arises when we consider $f_i (x_i) = \EE_{\xi_i \sim \mathcal{D}_i} [f_{\xi_i} (x_i)]$, where $\mathcal{D}_i$ is an unknown distribution, $f_{\xi_i} (x_i)$ represents the loss of model $x_i$ on sample $\xi_i$, and $f_i (x_i)$ is the generalization error. Since we do not know the distribution $\mathcal{D}_i$, we cannot work with $f_i (x_i) $ directly, and typically replace it with an approximation via Monte Carlo integration $f_i(x_i) = \frac{1}{M}\sum_{m=1}^M f_{i,m}(x_i)$. In this context, the problem is known as empirical risk minimization. This formulation is currently the main setting for solving supervised learning problems \cite{shalev2014understanding}.
Usually it is expensive to compute the full gradients $\nabla f_i(x_i)$ at each iteration, so instead, each node independently and uniformly chooses an index (batch number) $m_i$ and calculates the gradient $\nabla f_{i,m_i}(x_i)$ for that batch only. It turns out that we obtain the stochastic gradient typical of learning processes. Moreover $\nabla f_{i,m_i}(x_i)$ is an unbiased estimator of $\nabla f_{i}(x_i)$.}

The following assumption (a modification of Assumption \ref{as:func}) is used here.
% \begin{assumption}\label{ass_sum} We assume that each $f_i$ in problem \eqref{rewrite_main_problem_sum} is:
% \begin{itemize}[label=$\ast$]
% \item  $L$-average smooth w.r.t. $l_2$-norm, i.e. for all $u, v \in \mathbb{R}^{d}$, it holds that
% $$\frac{1}{M} \sum_{m=1}^M \|\nabla f_{i, m} (u) - \nabla f_{i, m} (v) \|^2 \leq L^2  \|u - v \|^2 ;$$ 
% \item  $\mu$ - strongly-convex w.r.t. $l_2$-norm i.e. for all $u, v \in \mathbb{R}^d$, it holds that
% $$\langle \nabla f_i(u) - \nabla f_i (v) , u - v\rangle \geq 2\mu \|u - v \|^2 .$$
% \end{itemize}
% \end{assumption}
\begin{assumption}\label{ass_sum} It is assumed that each $f_i:\R^d\to \R$ in problem \eqref{rewrite_main_problem_sum} is:
\begin{itemize}[label=$\ast$]
\item  $L$-average smooth w.r.t. $l_2$-norm, i.e. $\forall u, v \in \mathbb{R}^{d}$, 
$\frac{1}{M} \sum_{m=1}^M \|\nabla f_{i, m} (u) - \nabla f_{i, m} (v) \|^2 \leq L^2  \|u - v \|^2 ;$
\item  $\mu$-strongly-convex w.r.t. $l_2$-norm i.e. $\forall u, v \in \mathbb{R}^d$,
$\langle \nabla f_i(u) - \nabla f_i (v) , u - v\rangle \geq 2\mu \|u - v \|^2 .$
\end{itemize}
\end{assumption}
We present two approaches for solving problem \eqref{rewrite_main_problem_sum}. These approaches are efficient in the case of small $\lambda$. The key idea of the first approach (which uses the Accelerated Meta-Algorithm combined with L-Katyusha as the subproblem solver) is that problem \eqref{rewrite_main_problem_sum} is considered as composite problem \eqref{composite_problm}. In the second approach (Accelerated Randomized Algorithm for Decentralized Minimization) the ideas of variance reduction and importance sampling are used.  

\subsubsection{Accelerated Meta-Algorithm + L-Katyusha} 

As previously mentioned, the main idea behind this approach is to view problem \eqref{rewrite_main_problem_sum} as the  composite problem \eqref{composite_problm}. In particular, Section \ref{sec:analys} showed that  \eqref{rewrite_main_problem} can be solved by the Accelerated Meta-Algorithm with $h_1 = g$, $h_2 = f$. With this choice of $h_1$ and $h_2$ communications occur only in the outer loop, when we compute $\nabla g(\x) = \lambda W\x$. The local computations of $\nabla f_i (x_i)$ take place in the inner loop. But now, the inner problem \eqref{auxiliary_prob} has a finite-sum structure \ab{(since $f_i$ has a finite-sum structure and hence $h_2$ does as well). As previously mentioned, it is computationally expensive to use the full gradient for $h_2$, so typically for the subproblem \eqref{auxiliary_prob}, stochastic methods, such as the classical SGD method, are employed. Note that SGD converges only to a neighborhood of the solution, but for the finite-sum type problem it is known that one can use a variance reduction technique \cite{NIPS2013_ac1dd209,allen2017katyusha,L_Katyusha} to achieve convergence to an exact solution. 
%The essence of this approach is that we rarely compute the full gradient at some reference point. 
For this reason, we chose to use an accelerated and practical method that incorporates a variance reduction approach -- L-Katyusha \cite{L_Katyusha}.} 

\begin{theorem} \label{t4}
Let Assumption \ref{ass_sum} hold and let the graph $\mathcal{G}$ have an associated matrix $\hat W$ that satisfies Definition~\ref{gossipmatrix}. Then, to obtain an $\varepsilon$-optimal solution to problem \eqref{rewrite_main_problem}, solving by Algorithm \ref{meta_algorithm} combined with L-Katyusha, with
\begin{equation*}
    \delta = \frac{\e\mu}{3000^2(L + \lambda\eigWmax)^2},
\end{equation*}
requires the number of communications and local computations, respectively, to be of the order
\begin{equation*}
    N^{comm} = \mathcal{O}\left(\sqrt{\tfrac{\lambda\eigWmax}{\mu}}\log\tfrac{1}{\e}\log\tfrac{1}{\delta}\right),
\end{equation*}
and
\begin{equation*}
    N^{loc} = \mathcal{O}\left(\left( M\sqrt{\tfrac{\lambda\eigWmax}{\mu}}+ \sqrt{\tfrac{ML}{\mu }}\right)\log\tfrac{1}{\e}\log\tfrac{1}{\delta}\right).
\end{equation*}
\end{theorem}
The proof of this theorem is similar to the proof of Theorem \ref{th:alg_pfl_minimization} and can be found in Appendix \ref{app:th_ma_LK}. 

\begin{remark}
Note that the Accelerated Meta-Algorithm + L-Katyusha is suboptimal when $ M\lambda \lambda_{\max}(W) \leq L$.
\end{remark}

\subsubsection{Accelerated Randomized Algorithm for Decentralized Minimization.} In contrast with the previous approach, Algorithm \ref{alg_pfl_sum} uses variance reduction and importance sampling techniques and is based on L-Katyusha \cite{L_Katyusha}. We now view problem \eqref{rewrite_main_problem_sum} as being the sum of $M+1$ functions: there are $M$ functions $f_i$, as well as the composite term $g$. In Line 3 of Algorithm \ref{alg_pfl_sum} the value of a random variable $\xi^k$ determines what to choose: $f$ (make a local computation with probability $1-p$) or $g$ (make a communication with probability $p$). If the outcome is a local computation, then we choose index $i$ of the function $f$. \ab{We give a practical version of Algorithm \ref{alg_pfl_sum} in Appendix \ref{app:alg_prac} (Algorithm \ref{alg_pfl_sum_prac}).}

At each iteration of the algorithm, between 0 and 2 communications are made. As noted above, the first communication can take place if $\xi^k = 0$. And then the value of a random variable $\xi^{k  +\frac{1}{2}}$ determines whether to update $\hat \g^k$ or not. If $\hat \g^k$ is updated, then Algorithm  \ref{alg_pfl_sum} makes a  communication and a local computation. 
The following theorem states the convergence rate of Algorithm \ref{alg_pfl_sum}.
\begin{algorithm}
	\caption{Accelerated Randomized for Decentralized Minimization (\texttt{ARDM})}
	\label{alg_pfl_sum}
	\hspace*{0.2cm} {\bf Input:} starting point $\x_0 \in \R^{nd}$, number of iterations $K$,  parameters $0<\theta_1, \theta_2 <1$, $\eta, \beta, \gamma > 0$, probabilities $p$, $\rho$\\
    \hspace*{0.3cm}{\bf Initialization:} $ \y_0 = \z_0 = \u_0 = \x_0$ and $\hat \g^{0}=\lambda W\y^{0} + \nabla f(\y^{0})$
	\begin{algorithmic}[1]
\For {$k=0,1, 2, \ldots K-1$ }
\State $\x^k = \theta_1 \z^k + \theta_2 \u^k + (1 - \theta_1 -\theta_2)\y^k$
\State Generate $\xi^k =  \begin{cases}
 1,&  \text{with probability} ~~ 1 - p \\
0 ,& \text{with probability} ~~ p
\end{cases}$
\label{alg3:step5}
\If{$\xi^k = 0$} \label{alg3:step6}
\State $\g^{k} = \frac{\lambda}{p}\left(W\x^k - W\u^k\right) + \hat \g^{k}$
\Else  \label{alg3:step9}
\State \ab{Sample indices $m^k_1, \ldots m^k_n$ for each node independently and uniformly from $[M]$}
\State $\g^{k} = \frac{1}{1-p}\left(\nabla f_{m^k}(\x^k) - \nabla f_{m^k}(\u^k)\right) + \hat \g^{k}$ \ab{with $\nabla f_{m^k}(\x) = (\nabla f^T_{1, m^k_1}(x_1), \ldots, \nabla f^T_{n, m^k_n} (x_n))^T$} 
\EndIf
\State $\y^{k+1} = \x^k - \eta \g^k$
\State $\z^{k+1} = \beta \z^k + (1-\beta)\x^k + \frac{\gamma}{\eta}(\y^{k+1} - \x^k)$
\State Generate $\xi^{k+\frac{1}{2}}=  \begin{cases}
1,&  \text{with prob.} ~~ 1 - \rho \\
0 ,& \text{with prob.} ~~ \rho
\end{cases}$
\If{$\xi^{k+1/2} = 0$}
\State $\u^{k+1} = \y^{k+1}$
\State $\hat \g^{k+1} = \lambda W\y^{k+1} + \nabla f(\y^{k+1})$
\Else 
\State $\u^{k+1} = \u^{k}$
\State $\hat \g^{k+1} = \hat \g^{k}$
\EndIf
\EndFor
\end{algorithmic}
\end{algorithm}

\begin{theorem}\label{theorem_randomized}
Let Assumption \ref{ass_sum} hold and let the graph $\mathcal{G}$ have an associated matrix $\hat W$ that satisfies Definition~\ref{gossipmatrix}. Then, to obtain an $\varepsilon$-optimal solution to problem \eqref{rewrite_main_problem} using Algorithm \ref{alg_pfl_sum}, we can choose parameters $\gamma, \eta, \beta$, $p = \frac{\lambda \lambda_{\max}(W)}{\bar{L} + \lambda \lambda_{\max}(W)}$, and $\rho=\tfrac1M$ such that we need the following number of communications (on average)
\begin{equation*}
    N^{comm} = \mathcal{O}\left(\sqrt{\tfrac{\lambda\eigWmax}{\mu}}\log\tfrac{1}{\e}\right).
\end{equation*}
For $\rho = p$ we can achieve the following number of local computations (on average)
\begin{equation*}
    N^{loc} = \mathcal{O}\left(\left( M+ \sqrt{\tfrac{ML}{\mu }}\right)\log\tfrac{1}{\e}\right).
\end{equation*}
\end{theorem}
The proof of this Theorem can be found in Appendix \ref{app:rand}.

\begin{remark}
Accelerated Meta-Algorithm + L-Katyusha has optimal local computational complexity when $M\lambda \lambda_{\max} (W)\leq L$.  In contrast, the second algorithm has better local computation complexity (on average) if $\lambda \lambda_{\max}(W) < L$ and $M\lambda \lambda_{\max}(W) \geq L$. 
\end{remark}  

\section{Numerical Experiments}\label{sec:experiments}

In this section, we present several numerical experiments to demonstrate the practical advantages of the proposed approach for problem \eqref{rewrite_main_problem}. We study logistic loss functions,
$$f_i(x_i)=
\frac1n 
\left(\frac1{m_i}
  \sum_{j=1}^{m_i} 
   \log(1+e^{- y_j^i a_j^i x_i  })
   \right),
$$
where
$\{ (a_j^i, y_j^i)  \}_{j=1}^{m_i}$
is local dataset stored on each machine $i \in \{1,2,\dots,n\}$,
$a_j^i \in \R^d$ represents the feature vector
and
$y_j^i \in \{-1,1\}$ is the label. 
In the experiments, the power method was used to estimated the smoothness parameter of the objective function, as well as $\eigWmax$. 

\paragraph{Datasets.}
The experiments were performed on datasets from the LIBSVM \cite{chang2011libsvm} database.\footnote{The datasets are available at \url{https://www.csie.ntu.edu.tw/~cjlin/libsvmtools/datasets/}.} 
Table~\ref{tbl:datasets} shows the basic characteristics of the datasets that were used.
\begin{table}
    \centering
    \begin{tabular}{crr}
    \toprule
        dataset &  \# features ($d$) & \# samples \\  
        \midrule
        mushrooms & 112 &
        8,124		
        \\ \hdashline
        a9a  & 123 & 32,561	
        \\ \hdashline
        covtype.scale & 54
        & 581,012		 
        \\ \hdashline
        rcv1.binary & 47,236
        &20,242	 	
        \\
        \bottomrule
    \end{tabular}
    \caption{The number of features and number of samples for each dataset used in the numerical experiments.}
    \label{tbl:datasets}
\end{table}

\paragraph{The communication networks.}
In the experiments, three different network topologies were considered:
\begin{enumerate}
    \item \emph{Cyclic:} In this topology, devices are connected in a cycle, where each device is connected to it's two closest neighbours only. In this communication network, it takes $\sim \frac n2$ iterations to transmit information between two devices on opposite sides of the cycle.
    \item \emph{Grid:} Here devices are organized in a $\sqrt{n} \times \sqrt{n}$ grid, and are connected to their nearest neighbours.\footnote{
\url{https://networkx.org/documentation/networkx-1.10/reference/generated/networkx.generators.classic.grid_2d_graph.html}}
    \item \emph{Erdos:} A
random communication graph, also known as an Erdős-Rényi graph.\footnote{\url{https://networkx.org/documentation/stable/reference/generated/networkx.generators.random_graphs.erdos_renyi_graph.html}}
\end{enumerate}
% \begin{description}
% \item[{\it Cyclic}:] in this topology, devices are connected in a cycle, where each device is connected exactly to two closest devices. In this communication network, it takes $\sim \frac n2$ iterations to transmit information between two devices in the opposite sides of cycle.
% \item[{\it Grid}\footnote{
% \url{https://networkx.org/documentation/networkx-1.10/reference/generated/networkx.generators.classic.grid_2d_graph.html}
% } - ] we assume that devices are organized in $\sqrt{n} \times \sqrt{n}$ grid, where devices are connected to the the nearest neighbours.

% \item[{\it Erdos:}\footnote{\url{https://networkx.org/documentation/stable/reference/generated/networkx.generators.random_graphs.erdos_renyi_graph.html}} - ]
% is a 
% random communication graph, also known as an Erdős-Rényi graph.

% \end{description}
We used the {\bf networkx} python package\footnote{
The {\bf networkx} package \url{https://networkx.org/} is hosted at
\url{https://github.com/networkx/networkx}.
}
to generate random bi-directional graphs with the structures described above.
As highlighted in the theory, the algorithm depends on the parameters $\lambda$, $\lambda_{\max}(W)$ and $L$.
We ran several experiments with varying values of $\lambda$, where
\begin{equation}
\lambda = r \tfrac{L}{\lambda_{\max}(W)}, \qquad r > 0.     \label{eq:lambdaInexperiments}
\end{equation}

\subsection{Solving the sub-probleem}

% In Algorithm~\ref{meta_algorithm}, solving auxiliary problem~\eqref{auxiliary_prob}
% is impossible in closed form solution as it requires communicating steps. 
% In order to solve obtains an approximate solution, we performed $T$ iterations of Nesterov accelerated algorithm.
Algorithm~\ref{meta_algorithm} requires the solution to the auxiliary problem~\eqref{auxiliary_prob}. To avoid communication costs, an approximate solution to \eqref{auxiliary_prob} was obtained by performing $T$ iterations of Nesterov's accelerated gradient method.
\begin{figure}
    \centering

\includegraphics[width=0.245\textwidth]{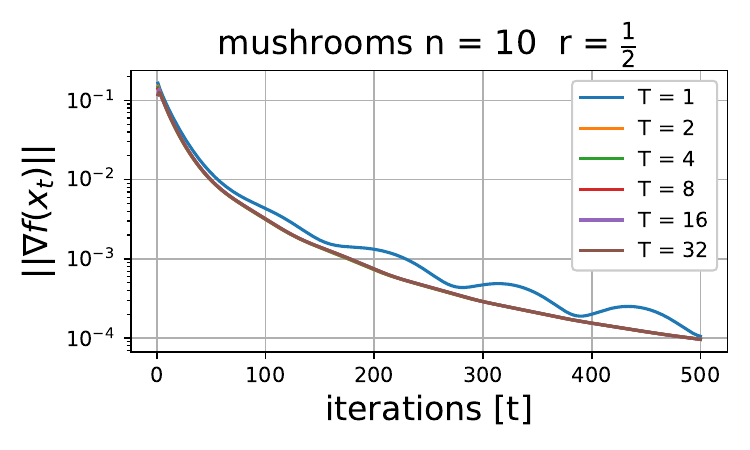}
\includegraphics[width=0.245\textwidth]{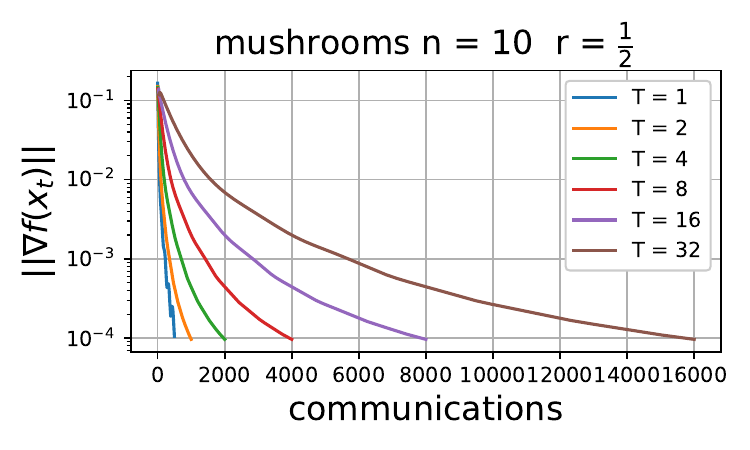}
\includegraphics[width=0.245\textwidth]{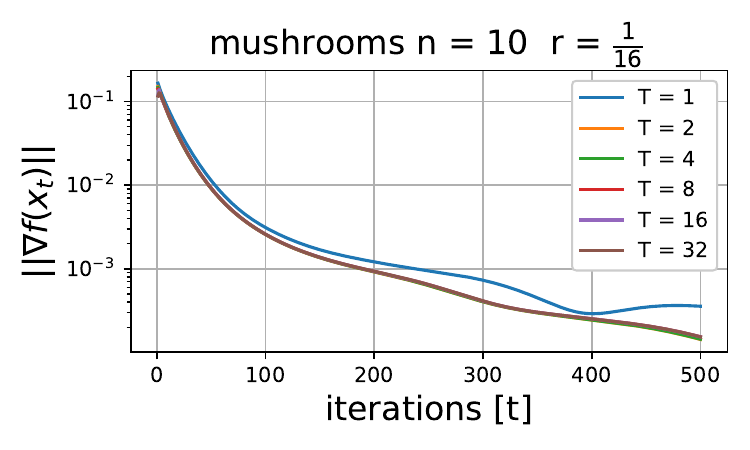}
\includegraphics[width=0.245\textwidth]{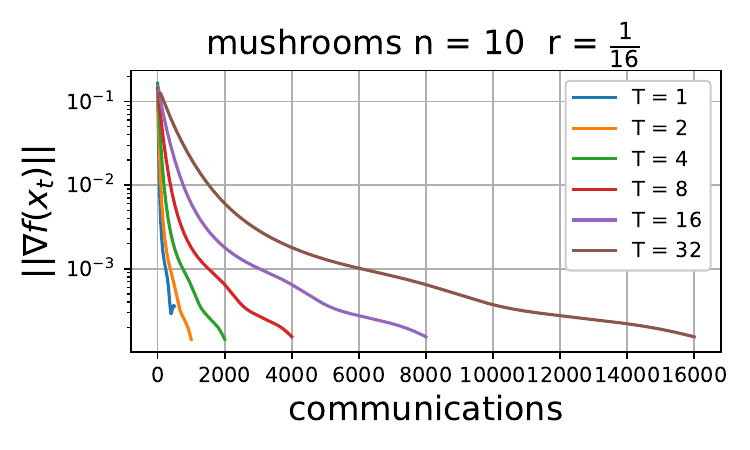}

\includegraphics[width=0.245\textwidth]{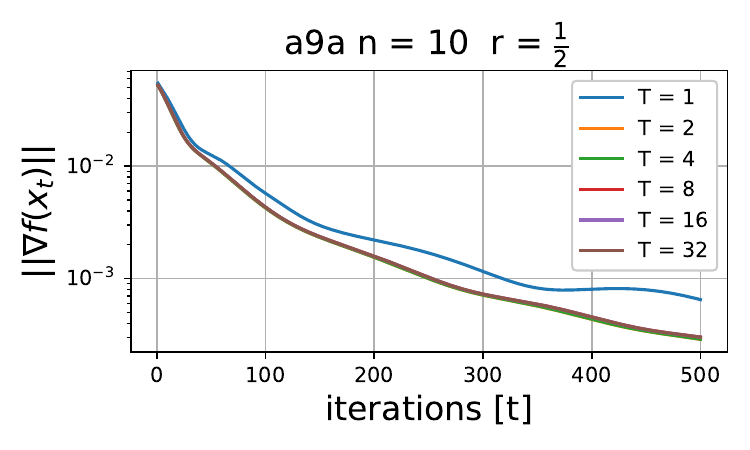}
\includegraphics[width=0.245\textwidth]{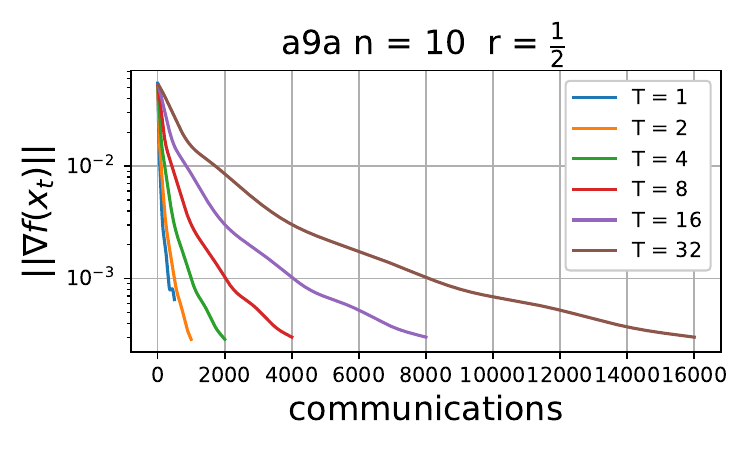}
\includegraphics[width=0.245\textwidth]{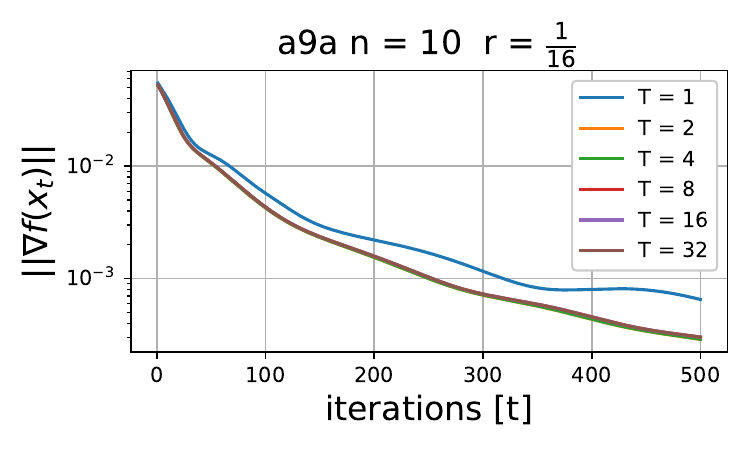}
\includegraphics[width=0.245\textwidth]{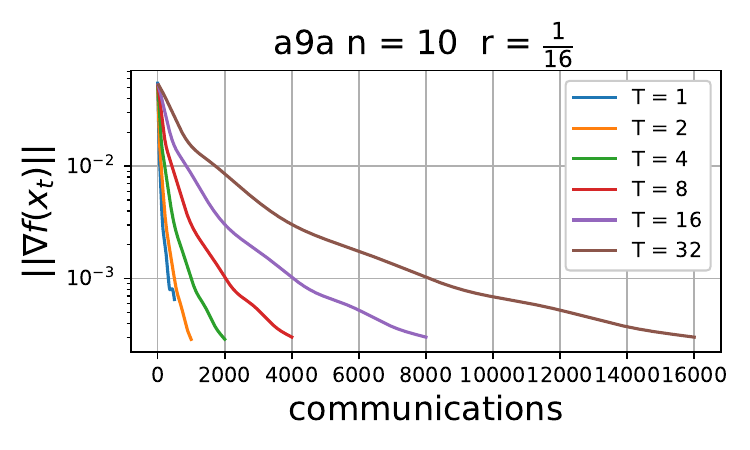}

\includegraphics[width=0.245\textwidth]{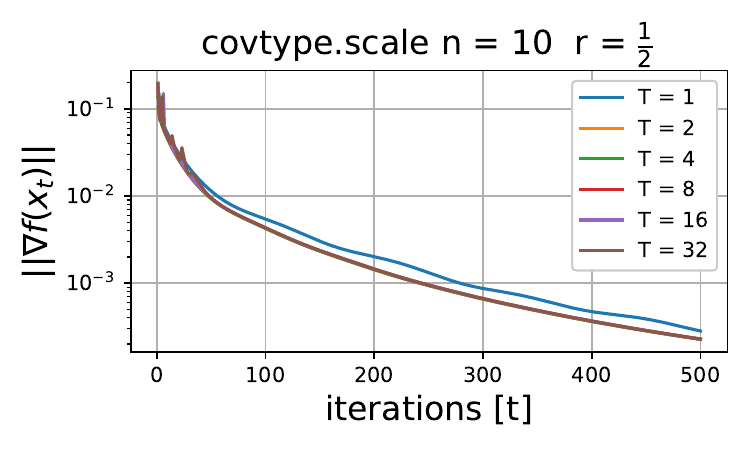}
\includegraphics[width=0.245\textwidth]{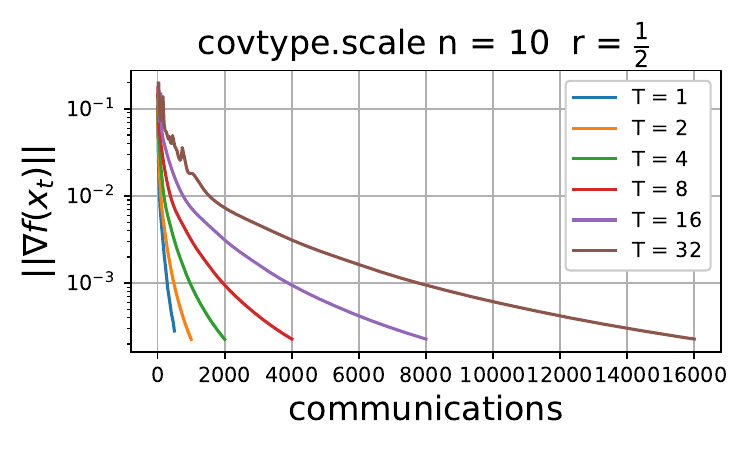}
\includegraphics[width=0.245\textwidth]{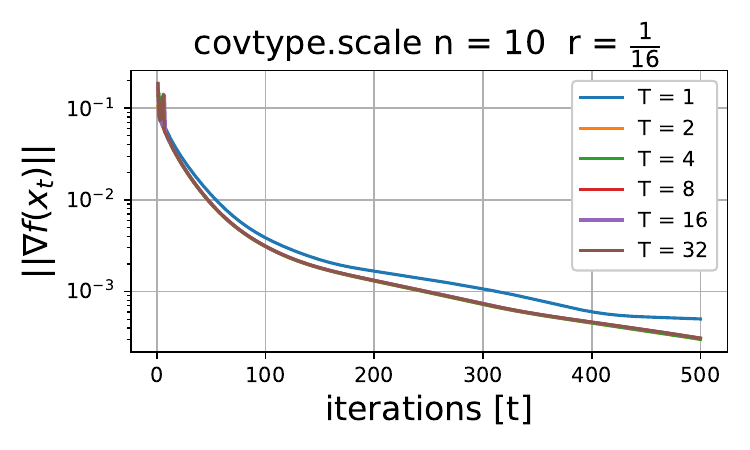}
\includegraphics[width=0.245\textwidth]{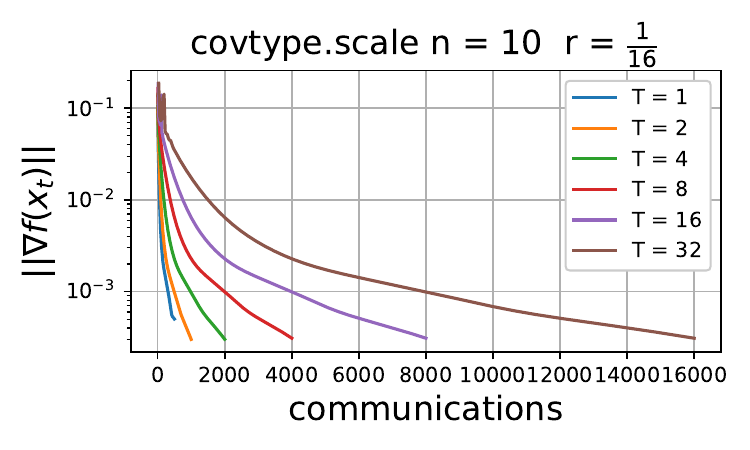}

\includegraphics[width=0.245\textwidth]{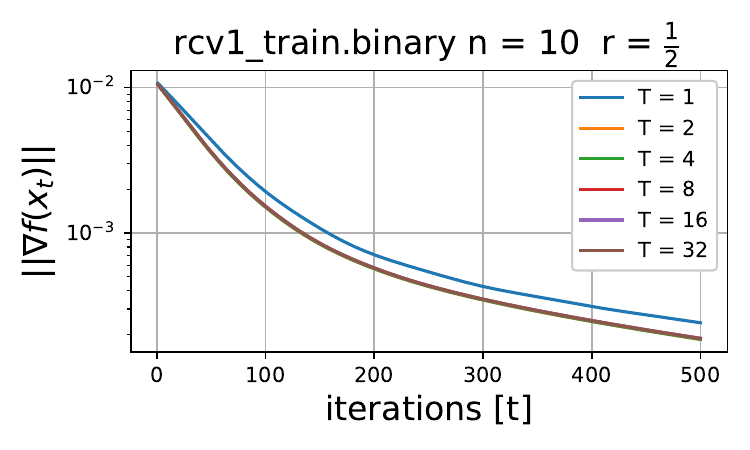}
\includegraphics[width=0.245\textwidth]{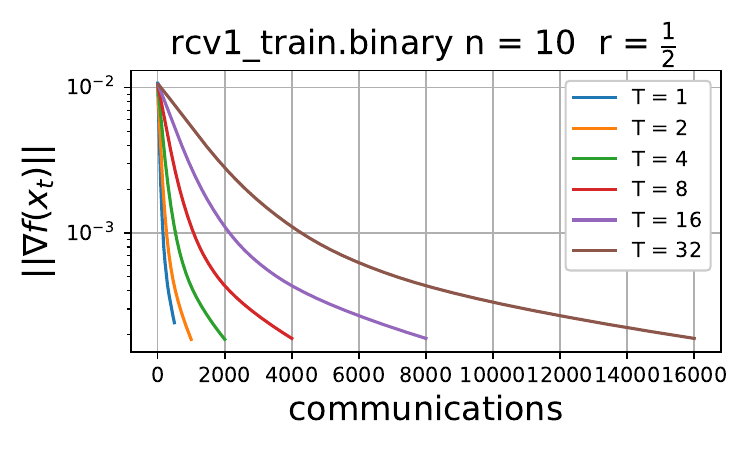}
\includegraphics[width=0.245\textwidth]{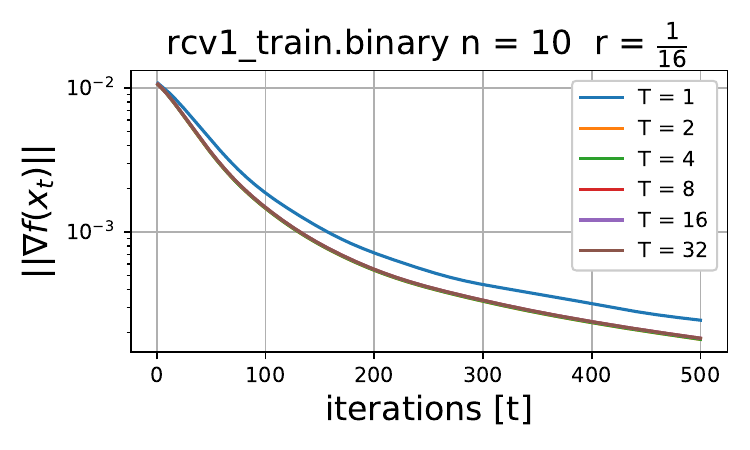}
\includegraphics[width=0.245\textwidth]{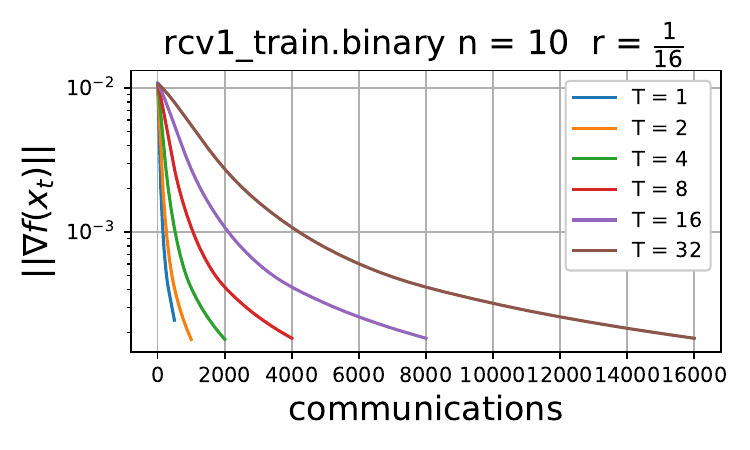}

    \caption{Evolution of $\|\nabla f(x_t)\|$ for different datasets,  regularization parameter and different level of solving the subproblem \eqref{auxiliary_prob} (larger $T$ means we optimize the subproblem better).}
    \label{fig:DependenceOnT}
\end{figure}
In Figure~\ref{fig:DependenceOnT}
we show the evolution of 
$\|\nabla f(x_t)\|$
for various selections of  parameter $T$.
Observe that the behaviour for $T\in\{2,4,\dots,32\}$
is almost identical (in terms of the iterations of the algorithm), however, 
larger $T$ requires additional rounds of communications.
Therefore, in the following experiments we selected $T = 2$.

\subsection{Effect of the regularization parameter}
\label{sec:variousSmallr}
\begin{figure}
    \centering

\includegraphics[width=0.245\textwidth]{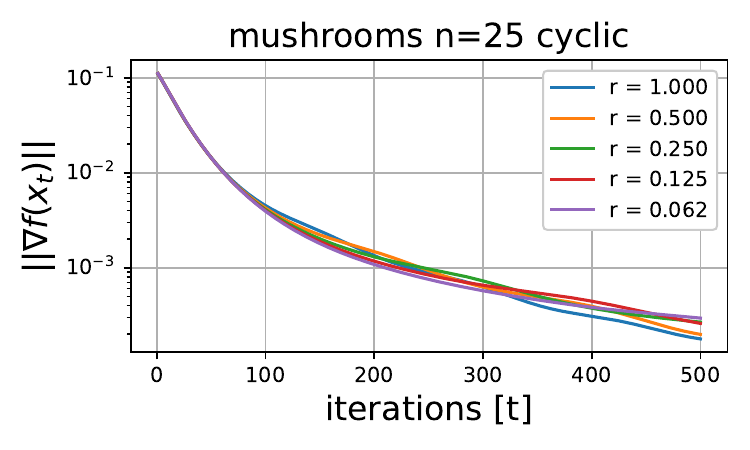}
\includegraphics[width=0.245\textwidth]{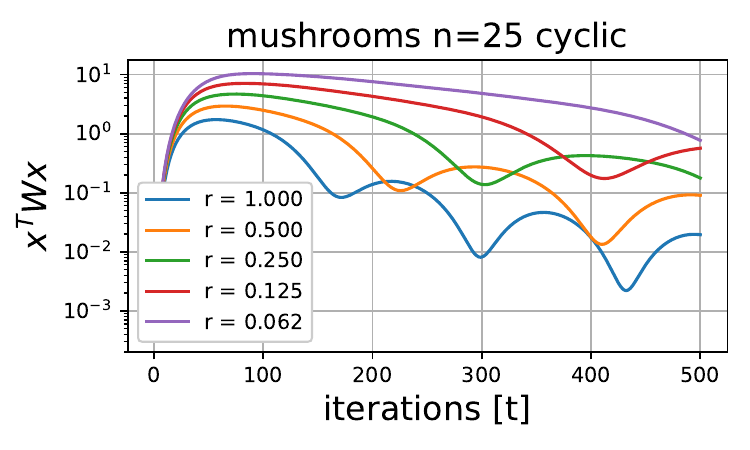}
\includegraphics[width=0.245\textwidth]{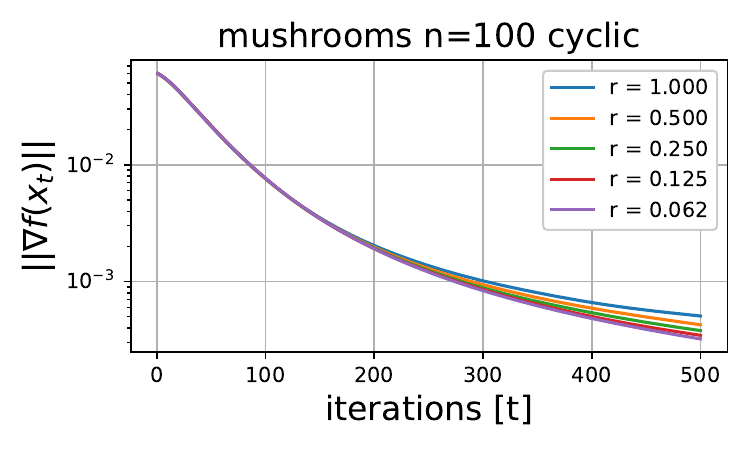}
\includegraphics[width=0.245\textwidth]{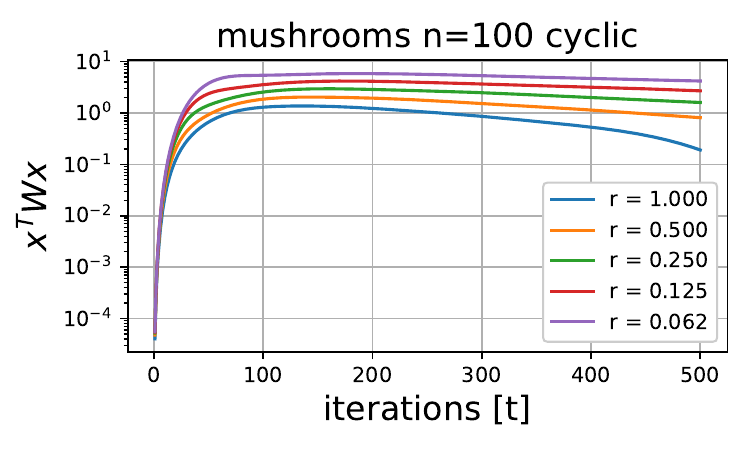}

\includegraphics[width=0.245\textwidth]{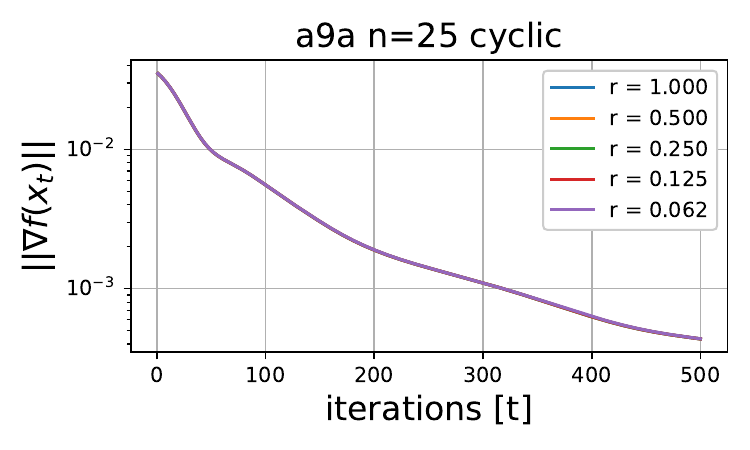}
\includegraphics[width=0.245\textwidth]{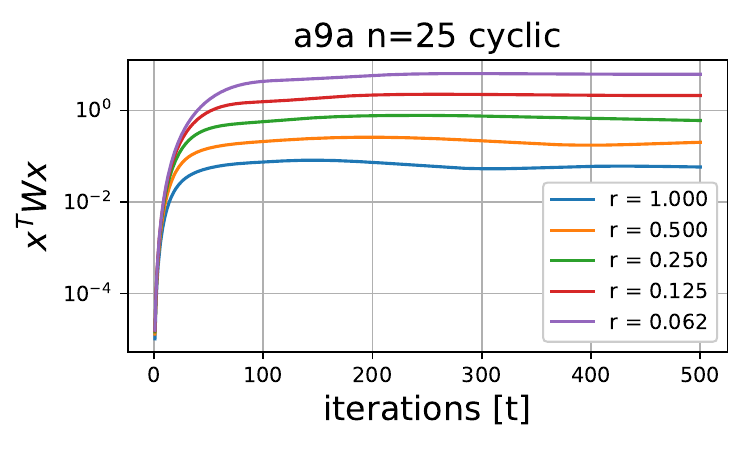}
\includegraphics[width=0.245\textwidth]{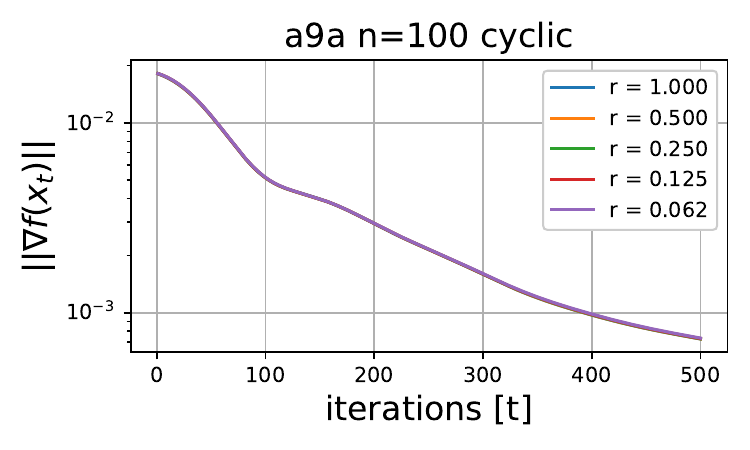}
\includegraphics[width=0.245\textwidth]{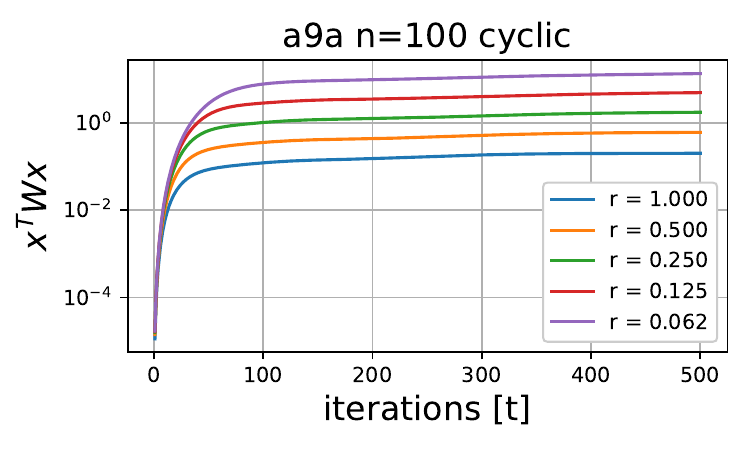}

\includegraphics[width=0.245\textwidth]{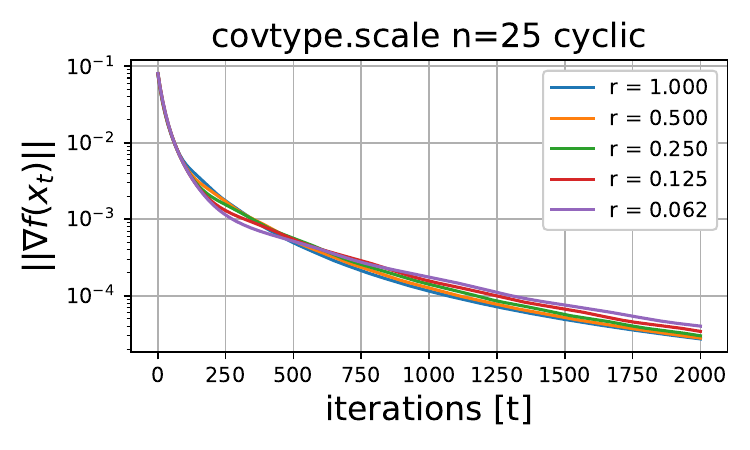}
\includegraphics[width=0.245\textwidth]{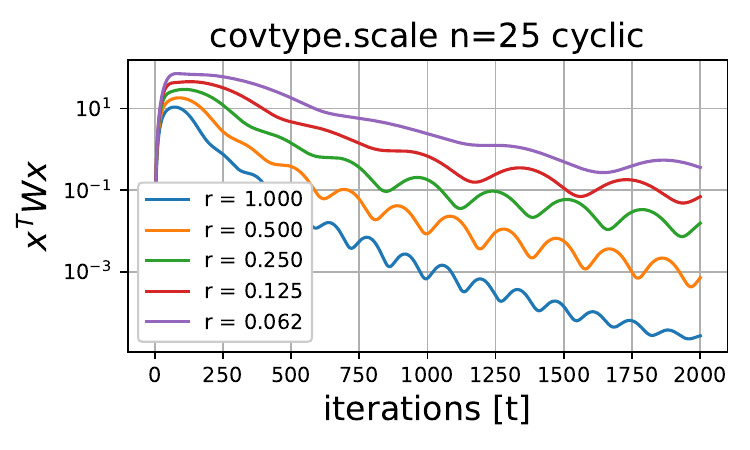}
\includegraphics[width=0.245\textwidth]{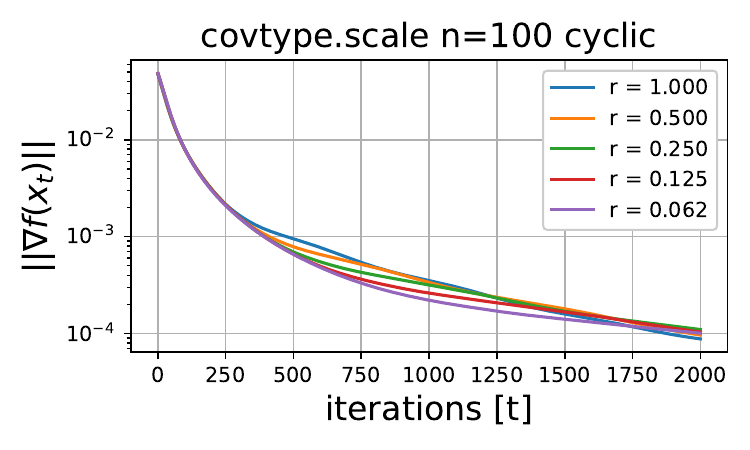}
\includegraphics[width=0.245\textwidth]{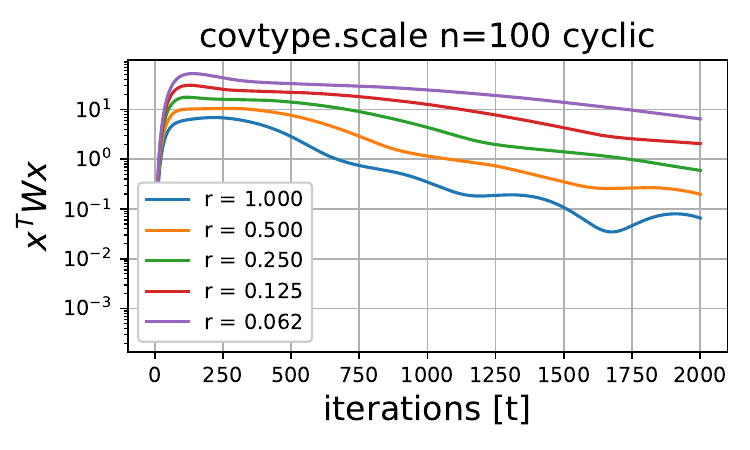}

\includegraphics[width=0.245\textwidth]{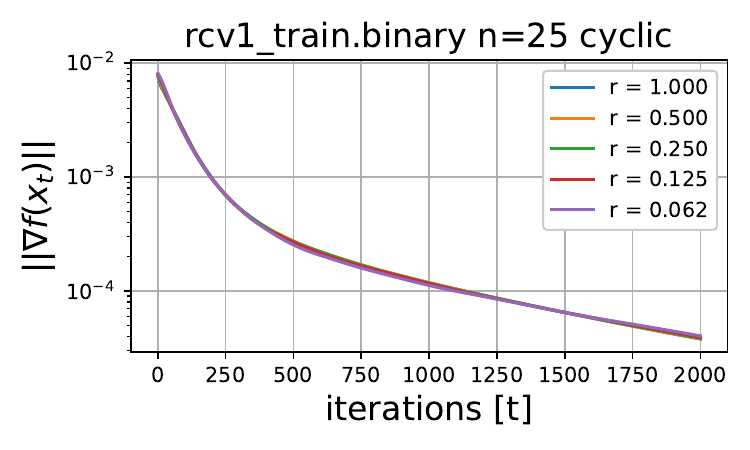}
\includegraphics[width=0.245\textwidth]{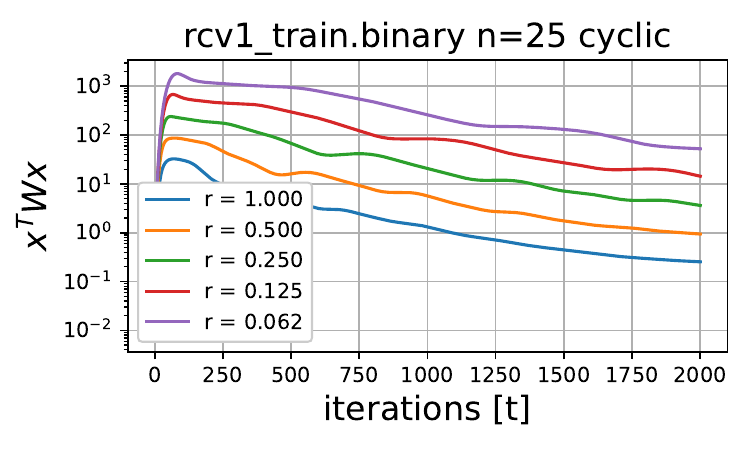}
\includegraphics[width=0.245\textwidth]{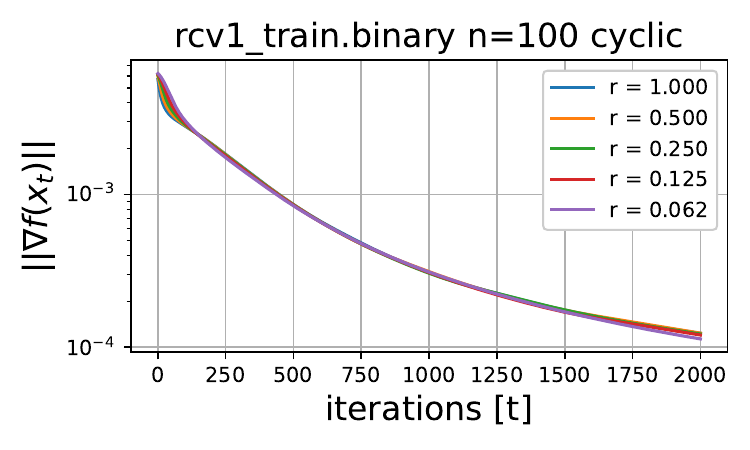}
\includegraphics[width=0.245\textwidth]{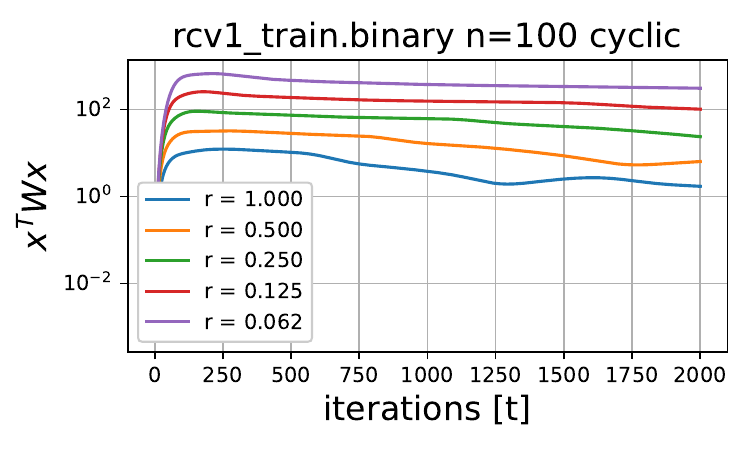}
    \caption{Evolution of $\|\nabla f(x_t)\|$ and $x_t^T W x_t$ for various datasets, levels of $\lambda$ and number of devices $n \in \{25, 100\}$ with cyclic network.}
    \label{fig:effectOfLambda}
\end{figure}
The main benefit of personalized federated learning is the ability to have slightly different local models, $x_i$, for each device $i$. The regularization term
$\lambda x^T W x$ penalizes local models (i.e., the $x_i$'s) that are different from their mean, where the parameter $\lambda$ controls the emphasis placed on this penalty term. When $\lambda$ is large, problem \eqref{rewrite_main_problem} tends to a consensus/classical federated learning problem, because there is a large penalty for models that are different at distinct devices. The current work focuses on \emph{personalized} federated learning, so here we consider the small $\lambda$ regime. 

\begin{figure}
    \centering

\includegraphics[width=0.245\textwidth]{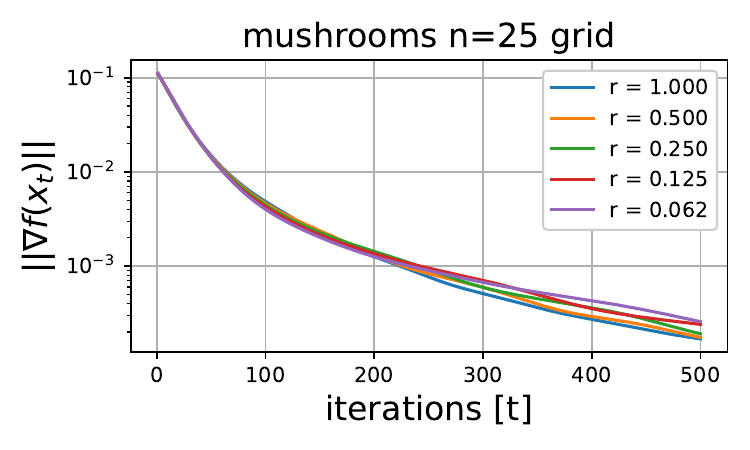}
\includegraphics[width=0.245\textwidth]{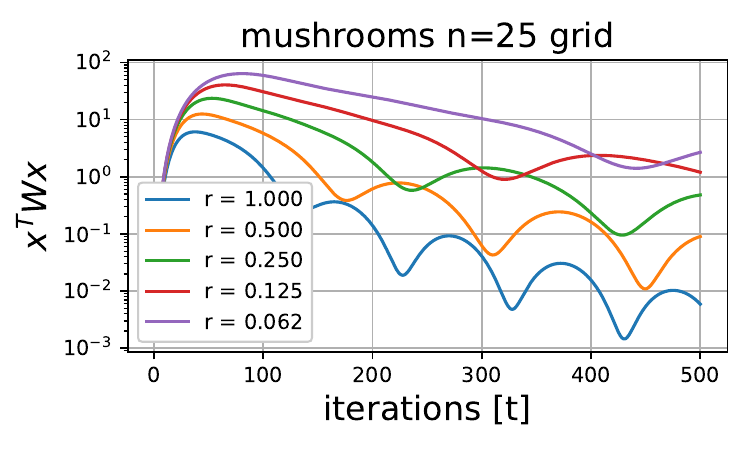}
\includegraphics[width=0.245\textwidth]{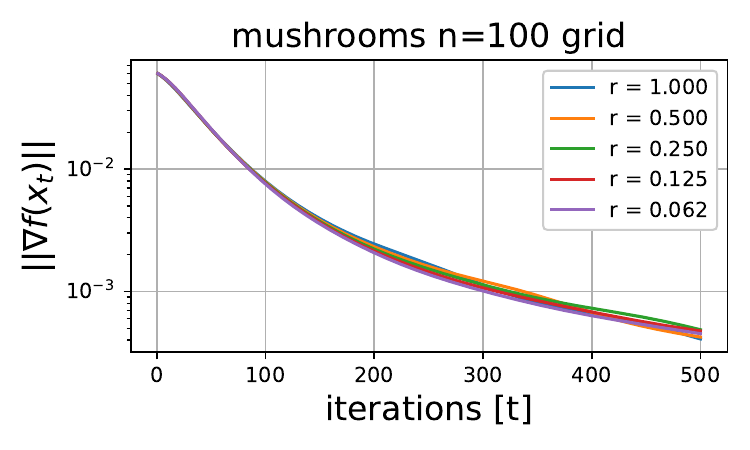}
\includegraphics[width=0.245\textwidth]{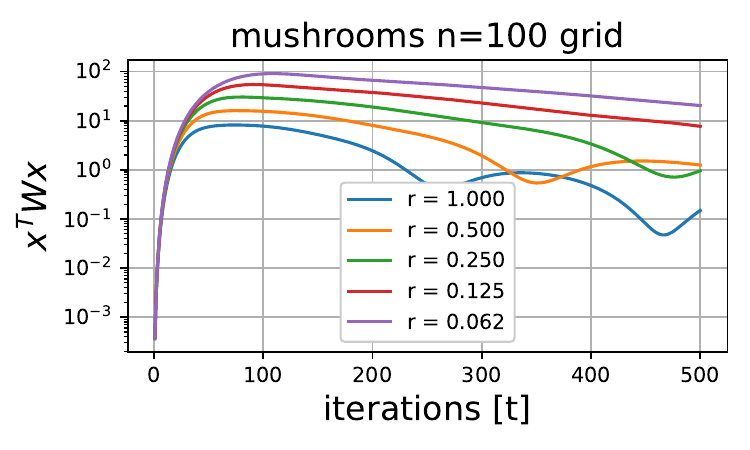}
\includegraphics[width=0.245\textwidth]{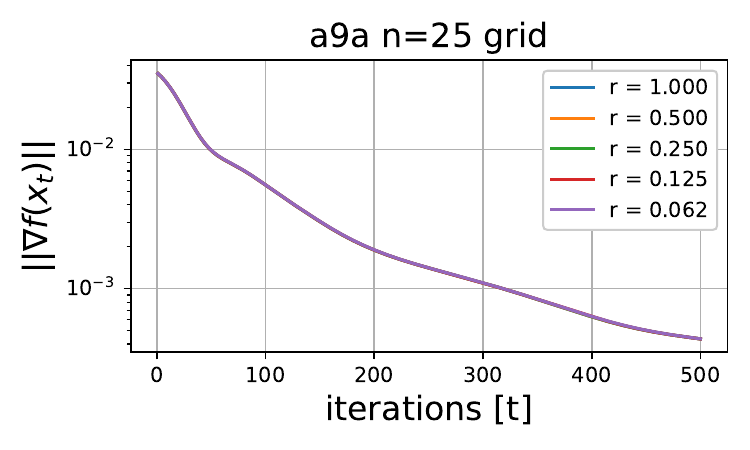}
\includegraphics[width=0.245\textwidth]{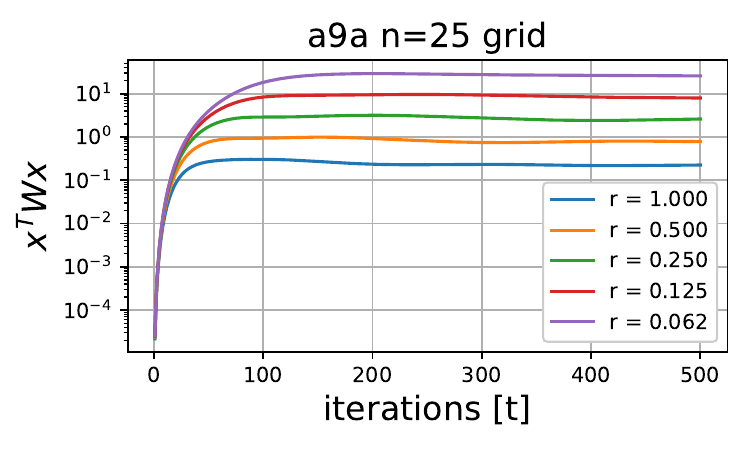}
\includegraphics[width=0.245\textwidth]{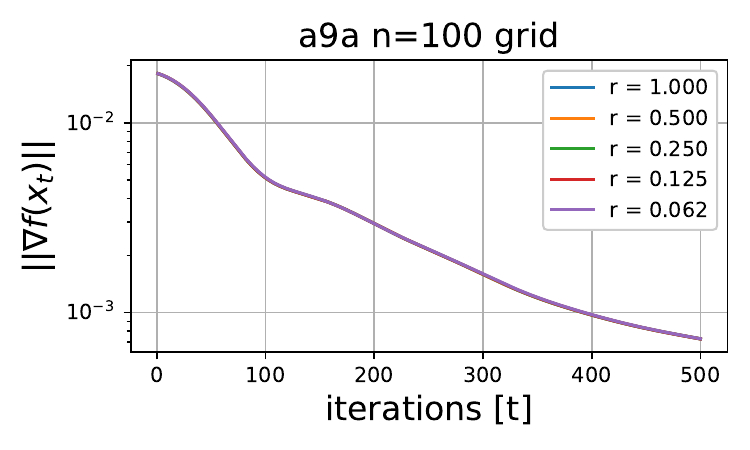}
\includegraphics[width=0.245\textwidth]{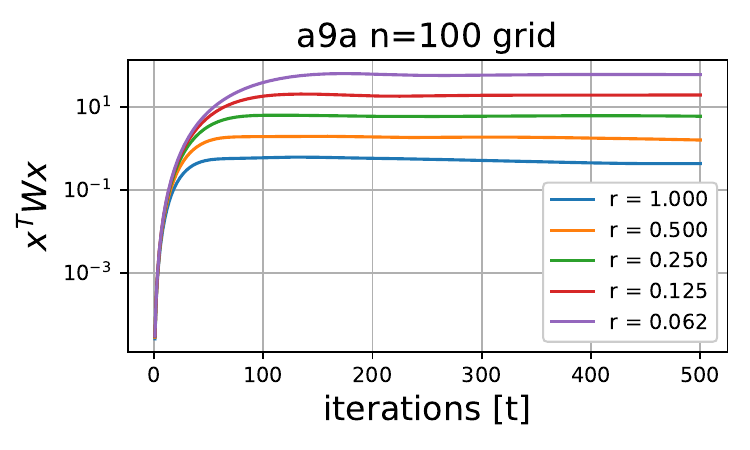}
\includegraphics[width=0.245\textwidth]{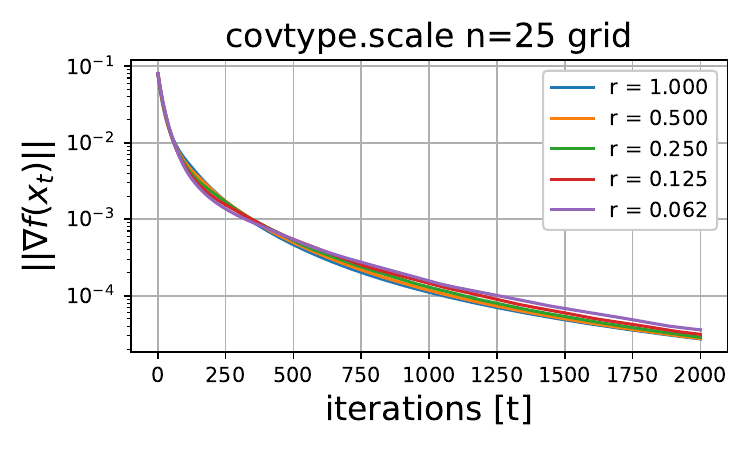}
\includegraphics[width=0.245\textwidth]{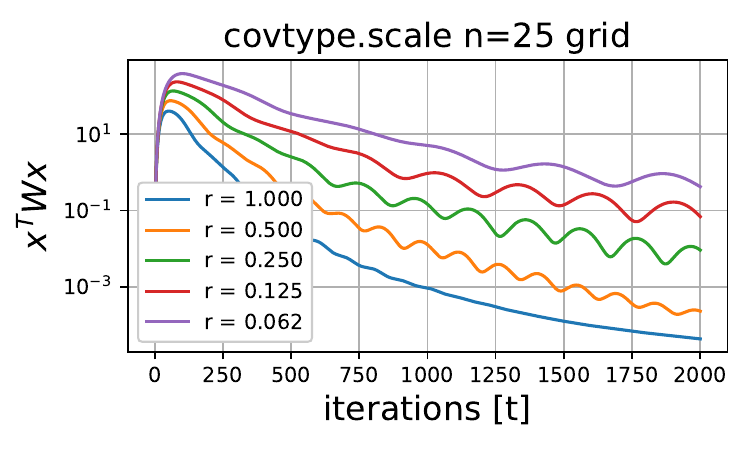}
\includegraphics[width=0.245\textwidth]{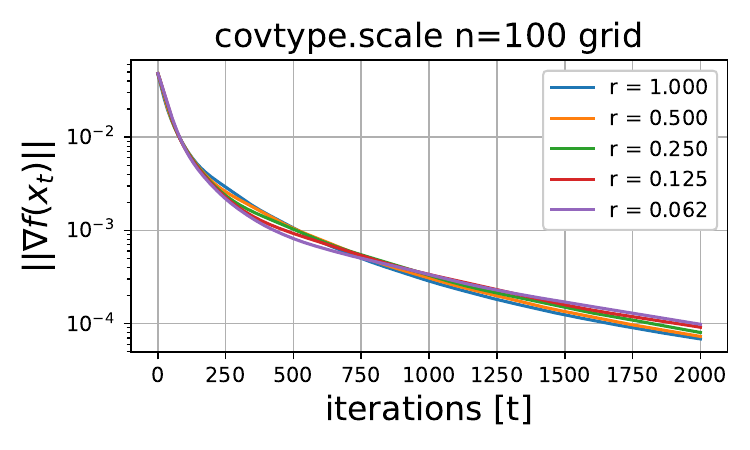}
\includegraphics[width=0.245\textwidth]{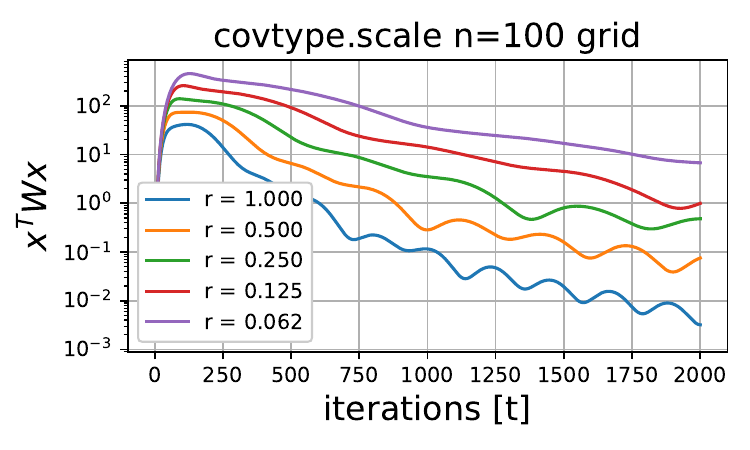}
\includegraphics[width=0.245\textwidth]{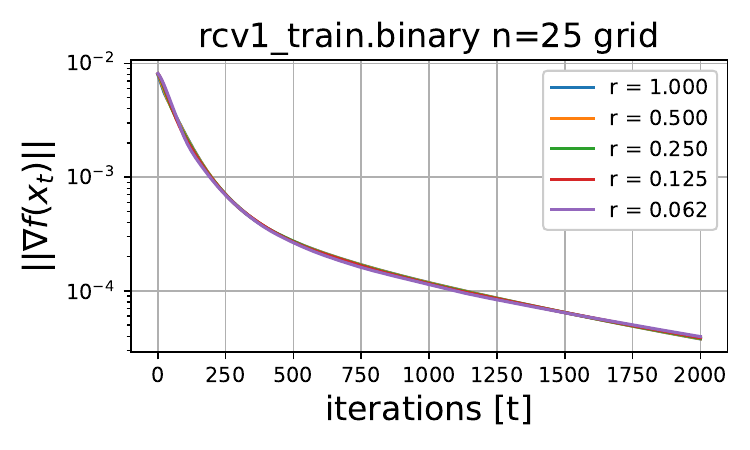}
\includegraphics[width=0.245\textwidth]{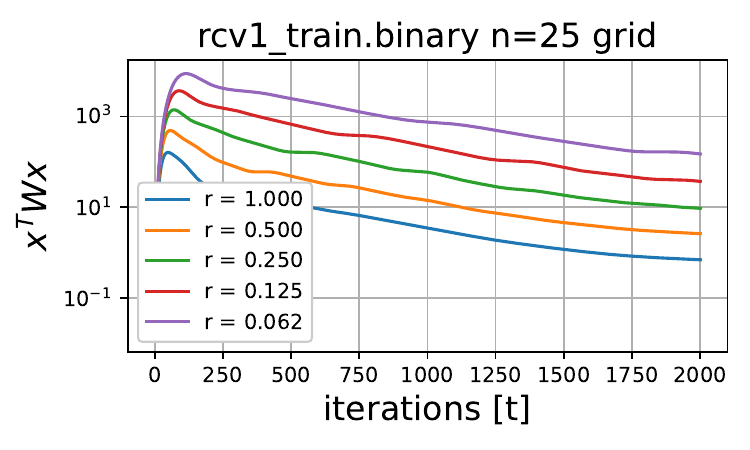}
\includegraphics[width=0.245\textwidth]{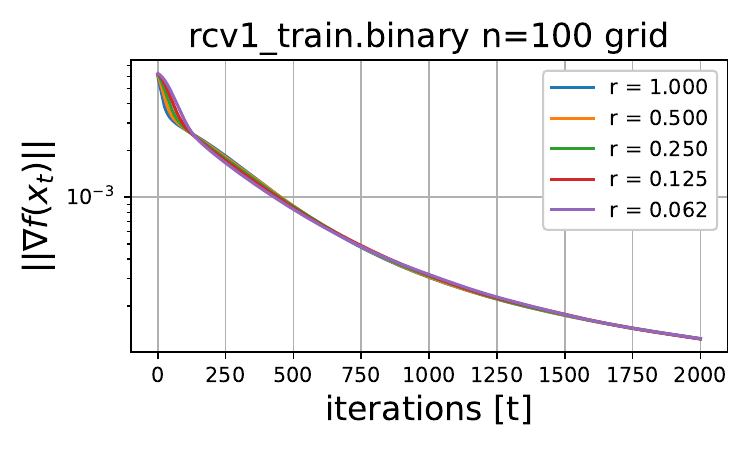}
\includegraphics[width=0.245\textwidth]{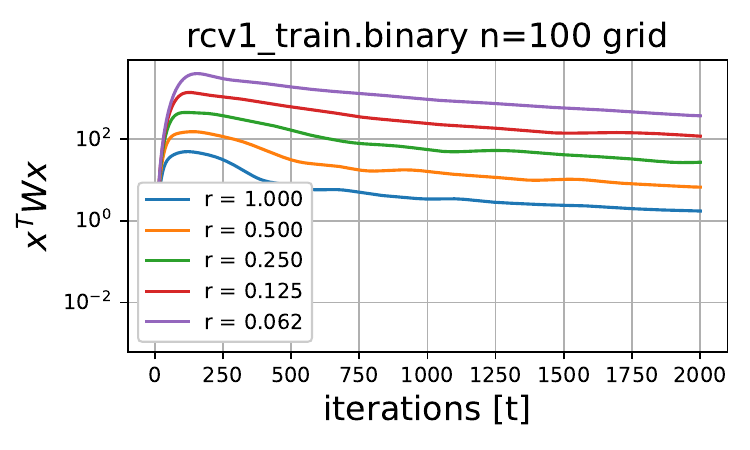}

    \caption{Evolution of $\|\nabla f(x_t)\|$ and $x_t^T W x_t$ for various datasets, levels of $\lambda$ and number of devices $n \in \{25, 100\}$ with grid network.}
    \label{fig:constrainGrid}
\end{figure}

\begin{figure}
    \centering

\includegraphics[width=0.245\textwidth]{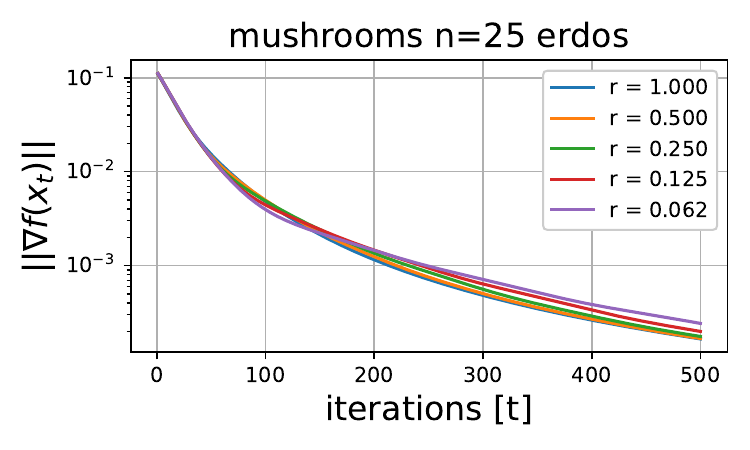}
\includegraphics[width=0.245\textwidth]{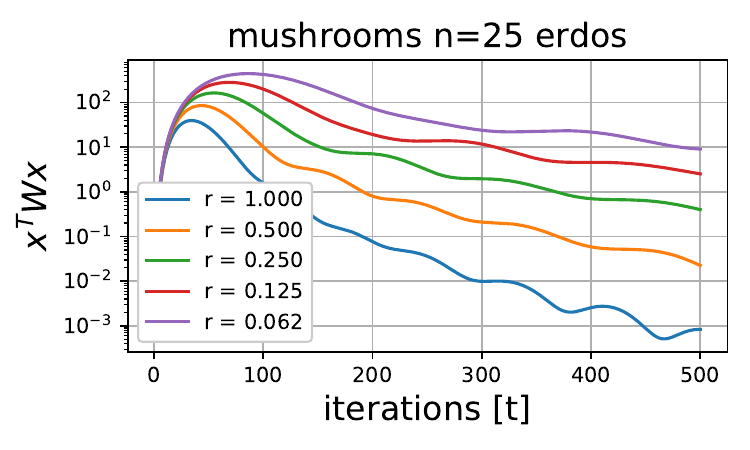}
\includegraphics[width=0.245\textwidth]{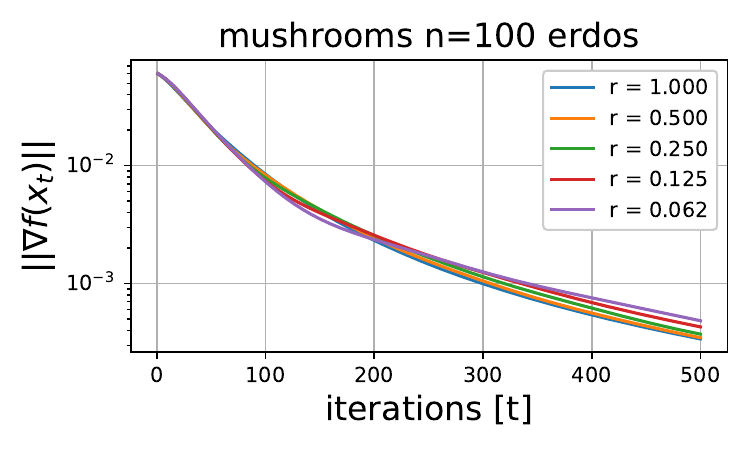}
\includegraphics[width=0.245\textwidth]{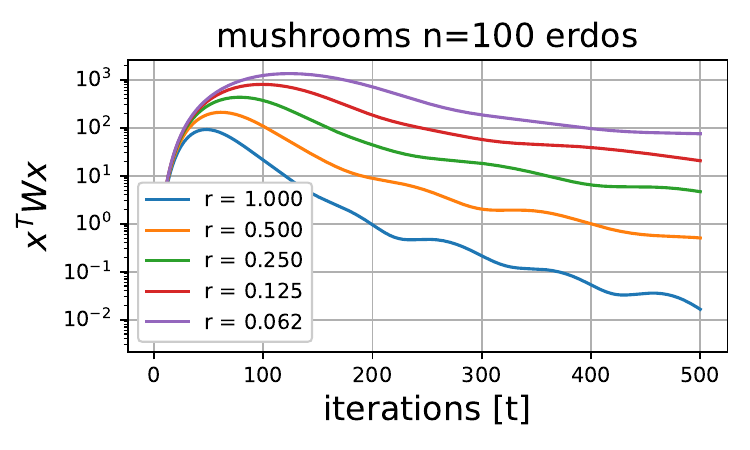}
\includegraphics[width=0.245\textwidth]{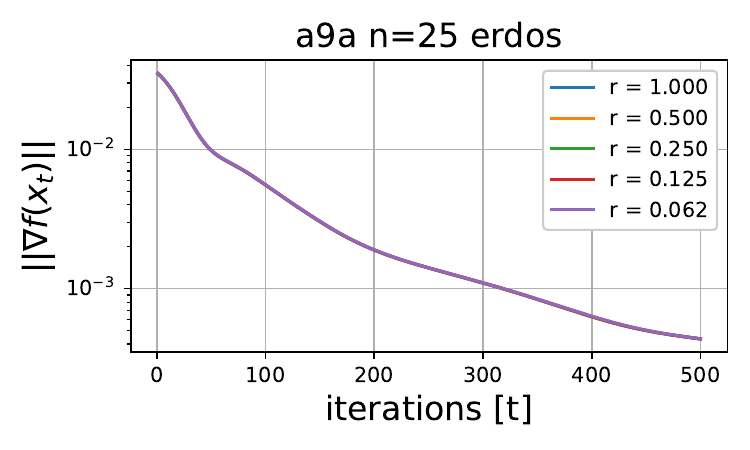}
\includegraphics[width=0.245\textwidth]{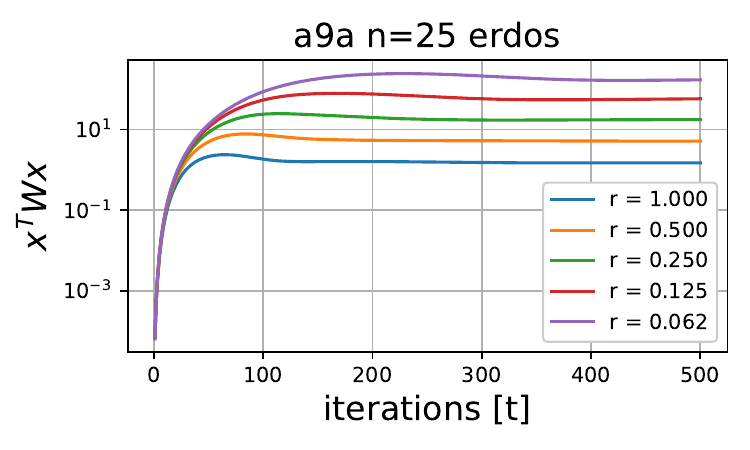}
\includegraphics[width=0.245\textwidth]{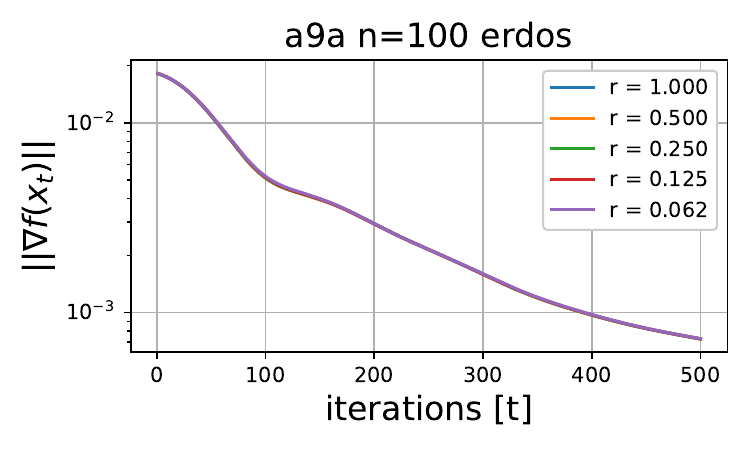}
\includegraphics[width=0.245\textwidth]{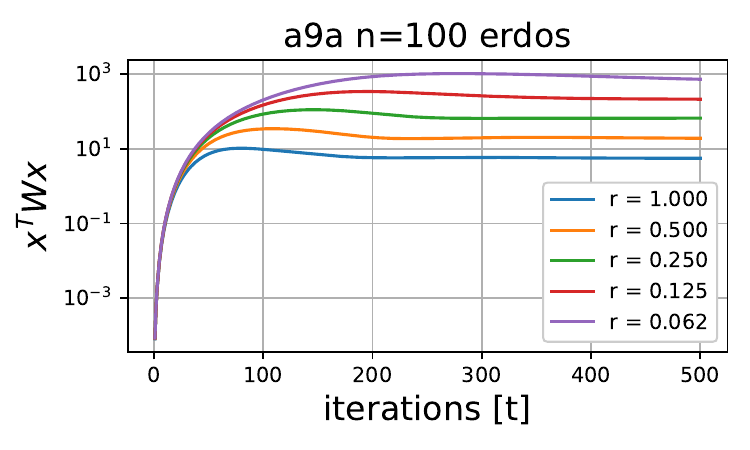}
\includegraphics[width=0.245\textwidth]{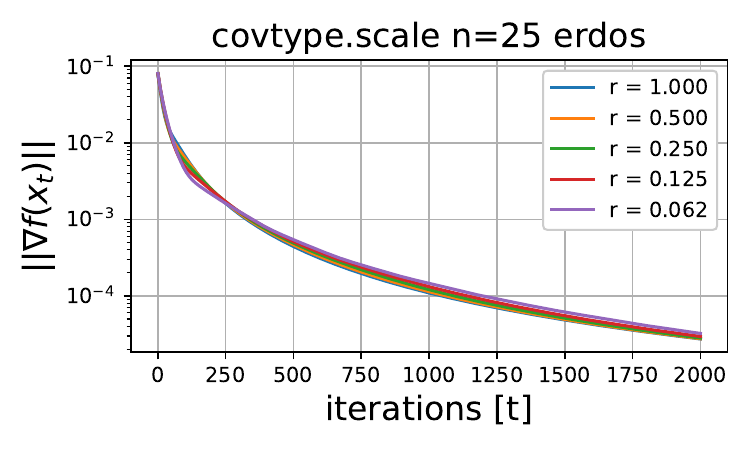}
\includegraphics[width=0.245\textwidth]{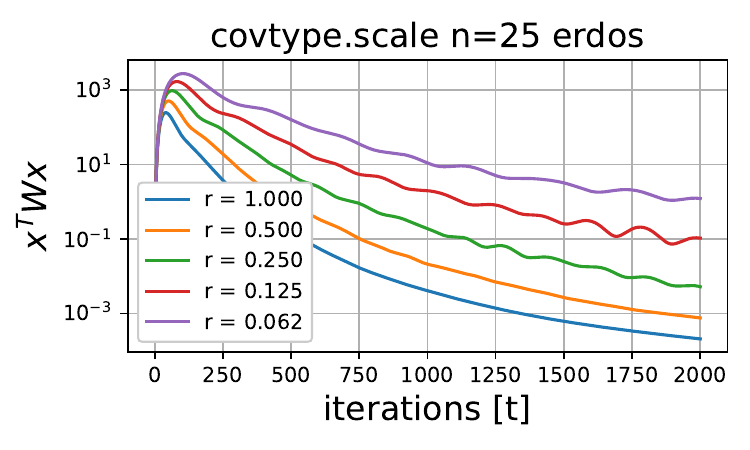}
\includegraphics[width=0.245\textwidth]{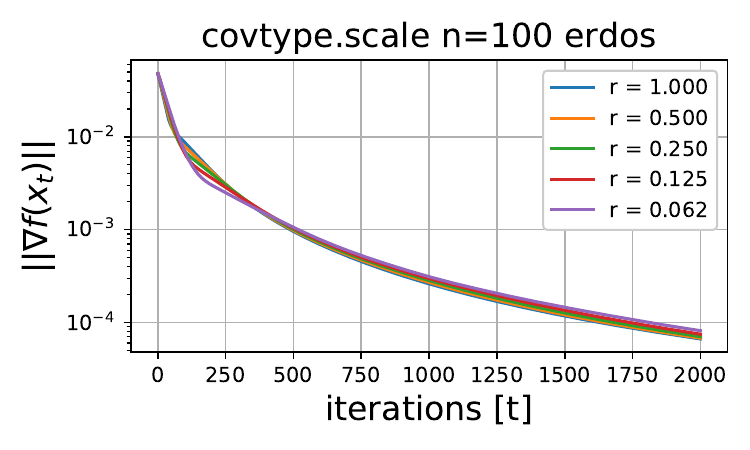}
\includegraphics[width=0.245\textwidth]{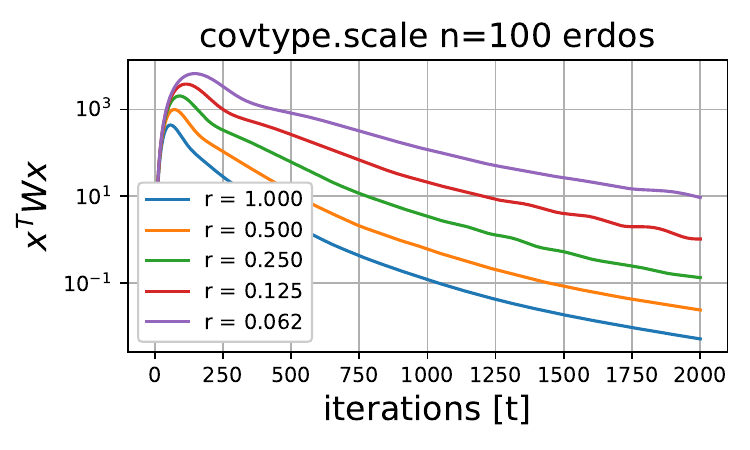}
\includegraphics[width=0.245\textwidth]{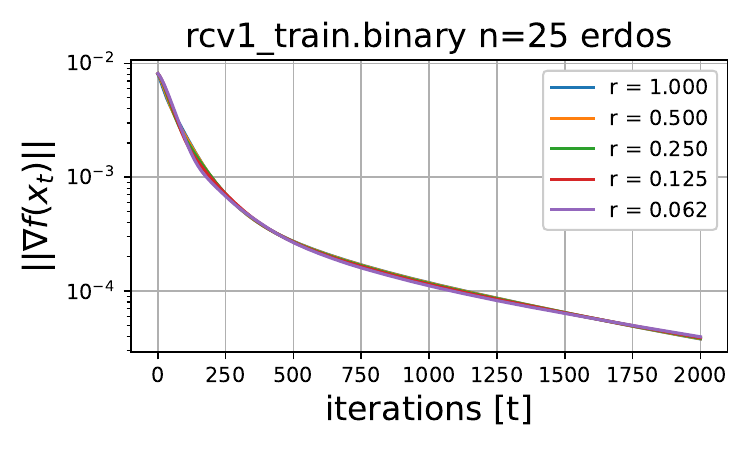}
\includegraphics[width=0.245\textwidth]{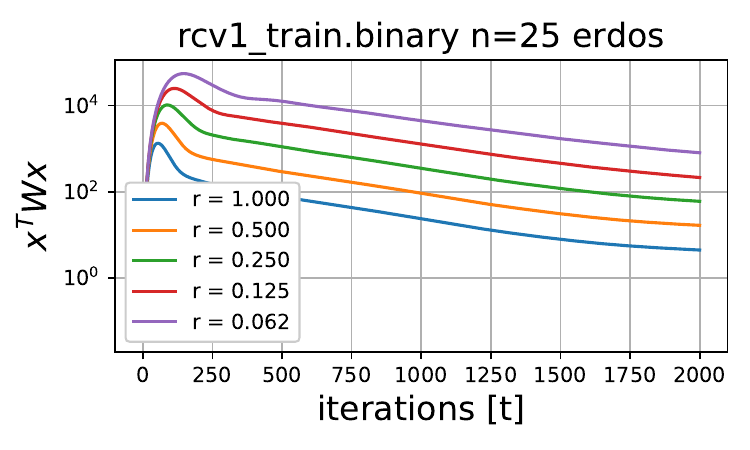}
\includegraphics[width=0.245\textwidth]{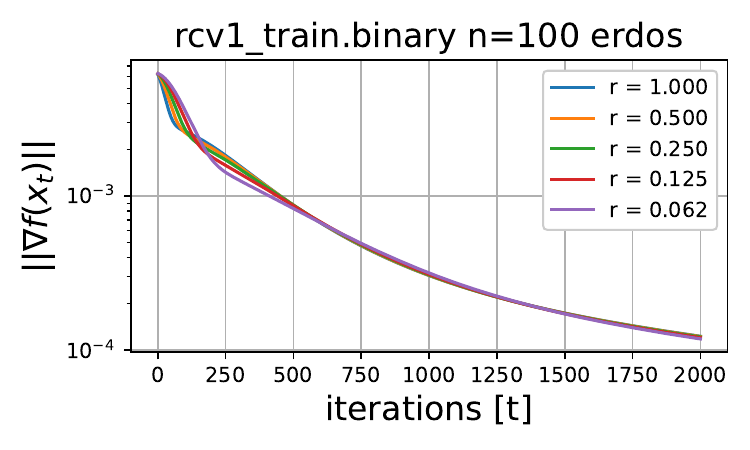}
\includegraphics[width=0.245\textwidth]{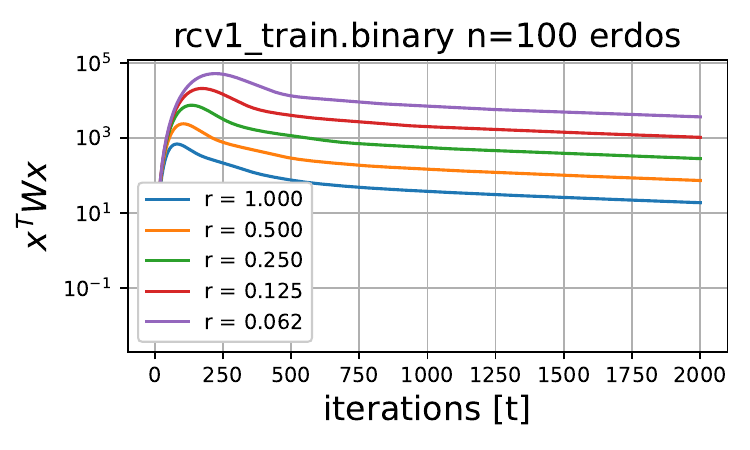}

    \caption{Evolution of $\|\nabla f(x_t)\|$ and $x_t^T W x_t$ for various datasets, levels of $\lambda$ and number of devices $n \in \{25, 100\}$ with erdos network.}
        \label{fig:constrainErdos} 
\end{figure}

Recall that for each problem instance considered, two iterations of accelerated gradient descent ($T=2$) were used to give an approximate solution to subproblem \eqref{auxiliary_prob}.
The parameter $\lambda$ is defined in \eqref{eq:lambdaInexperiments}, and several values of $r \in \{1,\tfrac12,\tfrac14,\tfrac18,\tfrac1{16}\}$ were used.
Let us stress that, as the number of local functions $n$ increases, the matrix $W$ changes, and hence, so too does $\lambda_{\max}(W)$. One can observe that, as expected, larger values of $\lambda$ (that corresponds to larger values of $r$) lead to solutions $x\in \R^{n d}$ that have a smaller value of the penalty term $x^T W x$. Figures~\ref{fig:effectOfLambda} (cyclic network),
        \ref{fig:constrainGrid} (grid network)
          and \ref{fig:constrainErdos}
            (Erdős-Rényi network) show the results of several numerical experiments. 

\subsection{Local Training Accuracy}

\begin{figure}
    \centering
\includegraphics[width=0.32\textwidth]{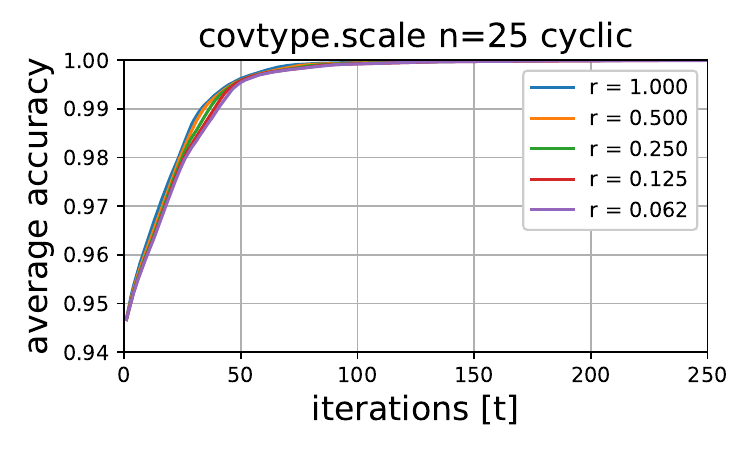}
\includegraphics[width=0.32\textwidth]{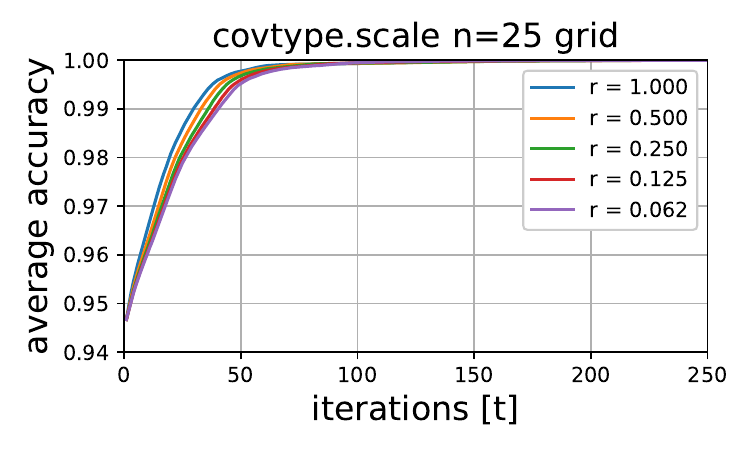}
\includegraphics[width=0.32\textwidth]{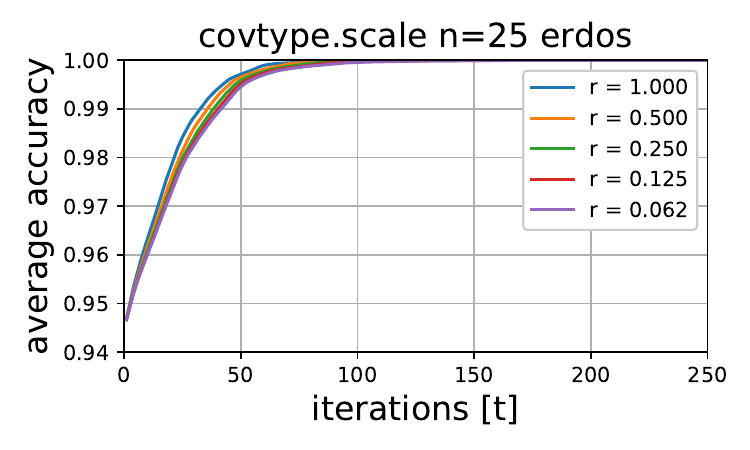}
\includegraphics[width=0.32\textwidth]{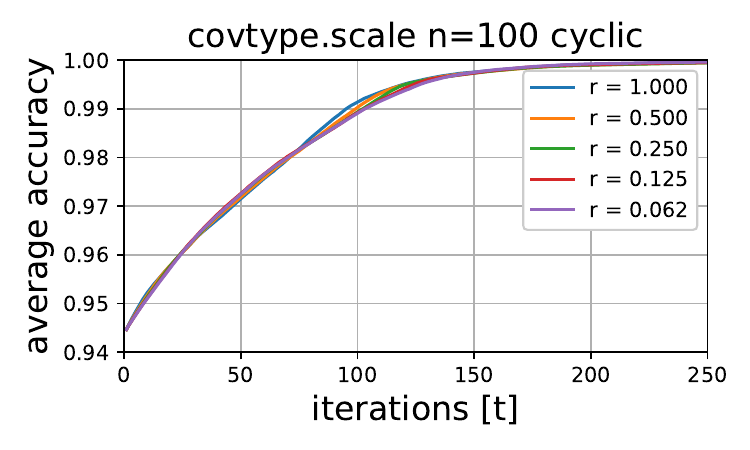}
\includegraphics[width=0.32\textwidth]{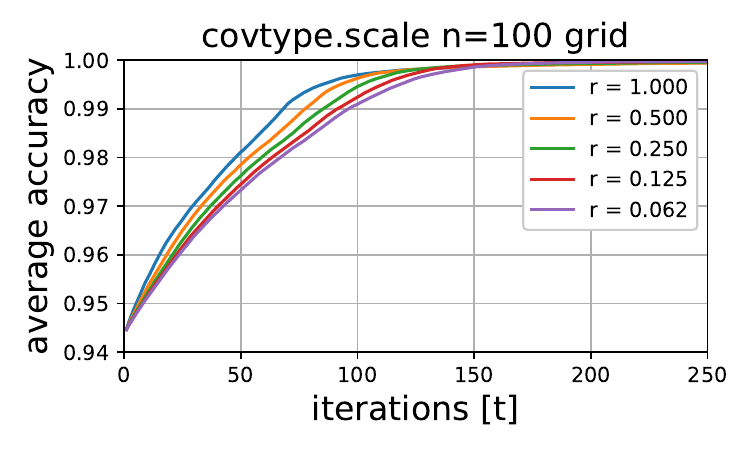}
\includegraphics[width=0.32\textwidth]{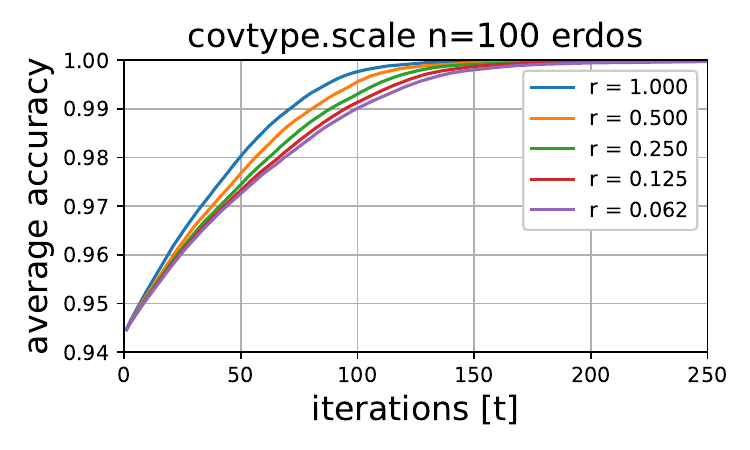}
\includegraphics[width=0.32\textwidth]{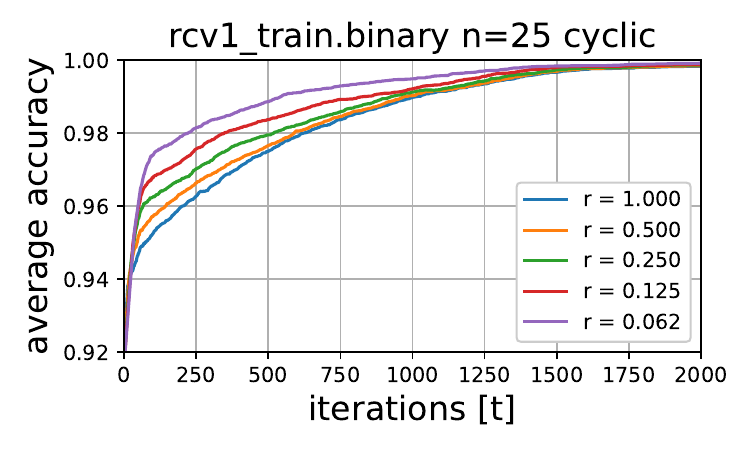}
\includegraphics[width=0.32\textwidth]{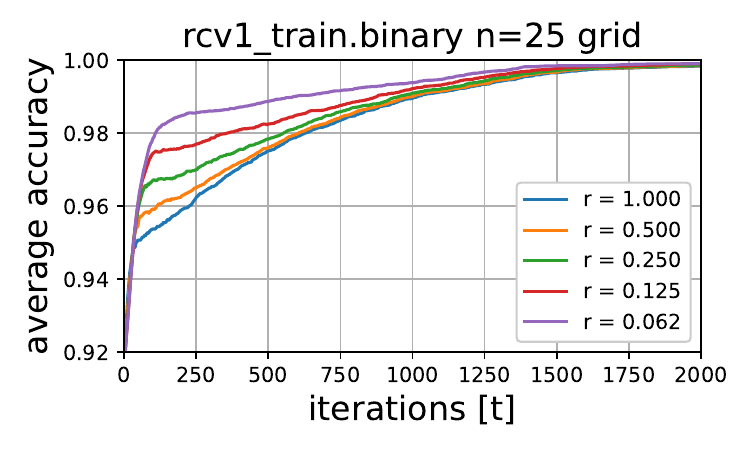}
\includegraphics[width=0.32\textwidth]{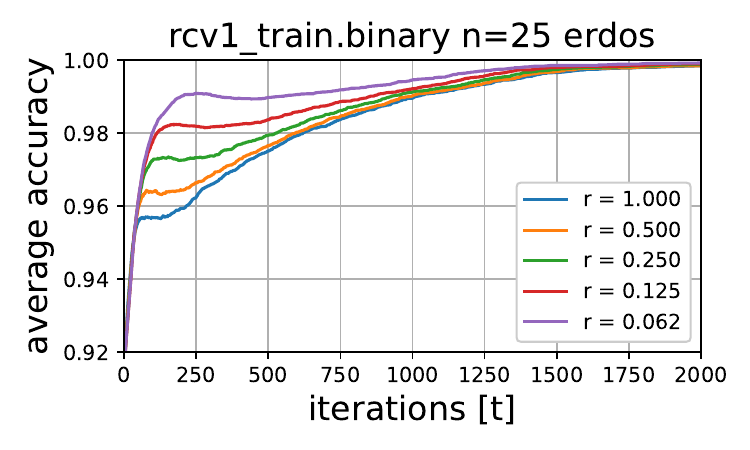}
\includegraphics[width=0.32\textwidth]{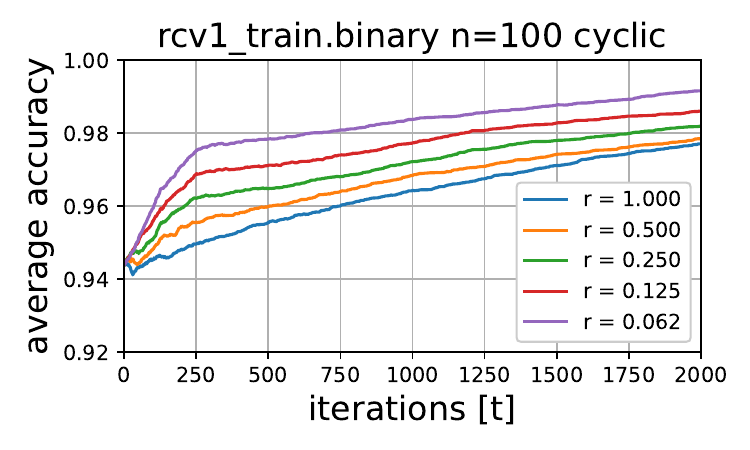}
\includegraphics[width=0.32\textwidth]{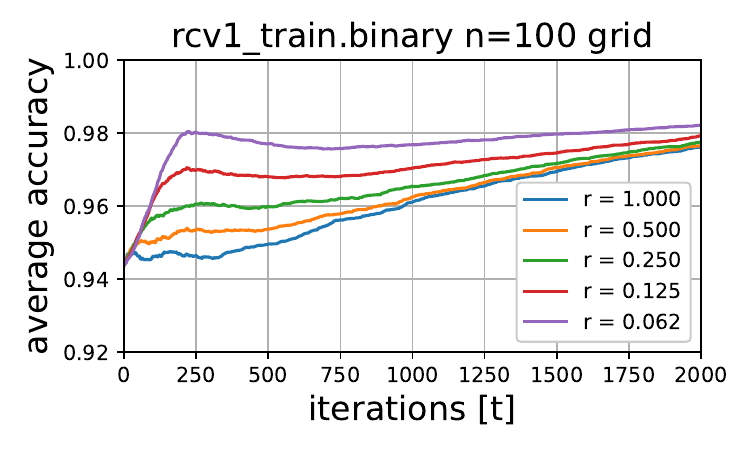}
\includegraphics[width=0.32\textwidth]{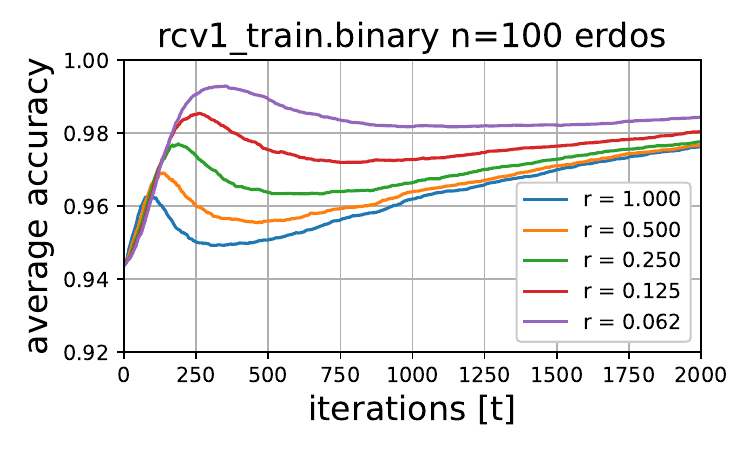}
    \caption{Comparison of average accuracy of local models on local data for various communication networks and levels of $\lambda$ ($r$).}
    \label{fig:networkEffect}
\end{figure} 
In Figure~\ref{fig:networkEffect} we demonstrate the main benefit of using PFL - namely, the ability for each device to have a slightly different local model, thereby capturing small differences in the local data. This is done by selecting various values of $\lambda$ and observing the affect that has on the 
training accuracy over various local functions $f_i$. We plot the average accuracy over local accuracies (each using their own set of parameters).
For the {\it mushrooms} and {\it a9a} datasets, the algorithm quickly achieved very good local accuracy for all local models. However, the {\it covtype} (more samples) and {\it rcv1} (more features) datasets were more challenging. 
For the {\it rcv1} dataset, we can see that initially (mainly due to the over-parametrization of the data) the local models achieve better accuracy (for smaller value of $\lambda$ ($r$)), demonstrating the advantages of PFL.

\begin{figure}
    \centering

\includegraphics[width=0.245\textwidth]{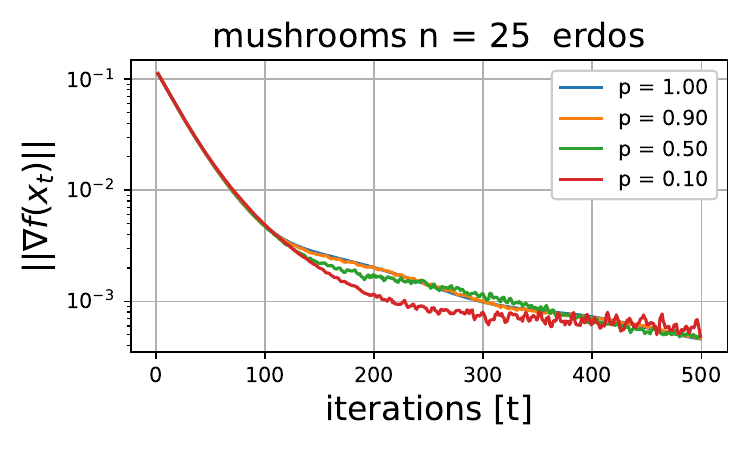}
\includegraphics[width=0.245\textwidth]{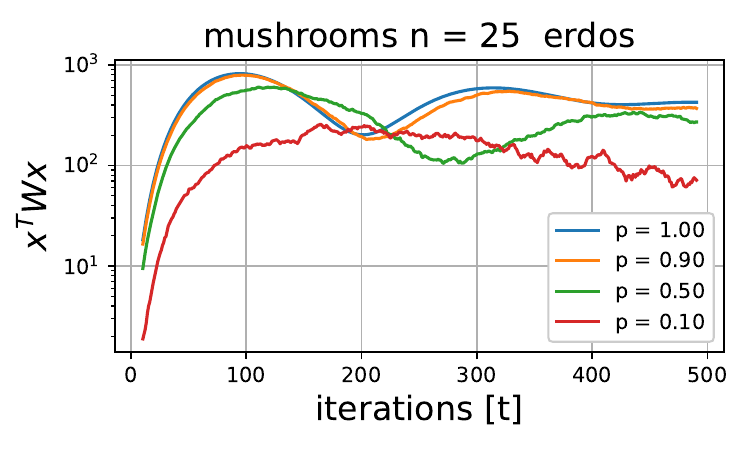}
\includegraphics[width=0.245\textwidth]{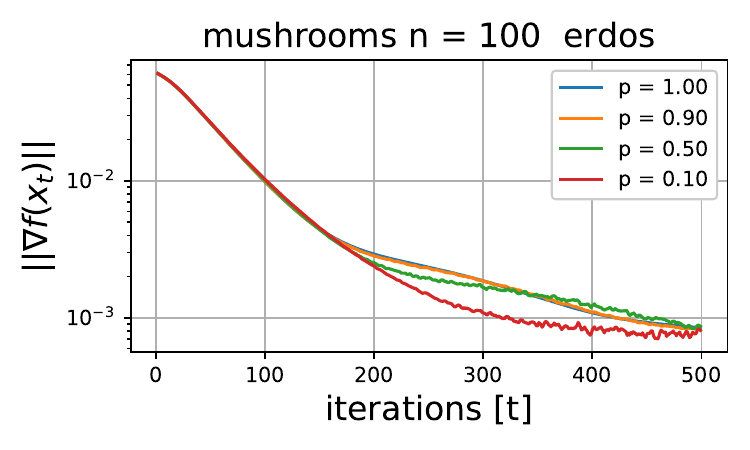}
\includegraphics[width=0.245\textwidth]{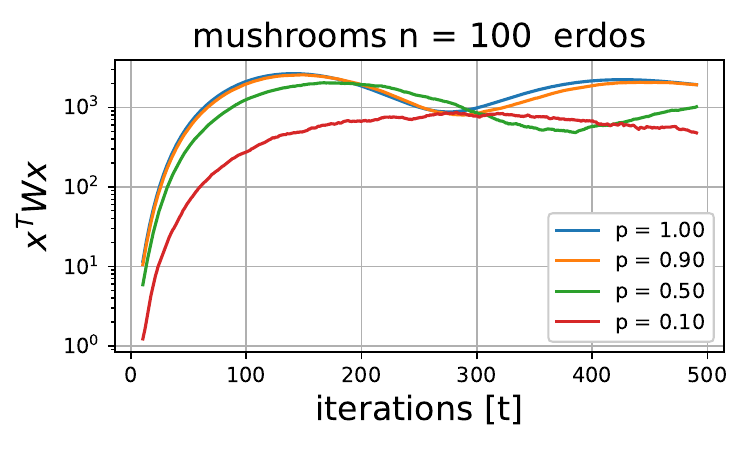}

\includegraphics[width=0.245\textwidth]{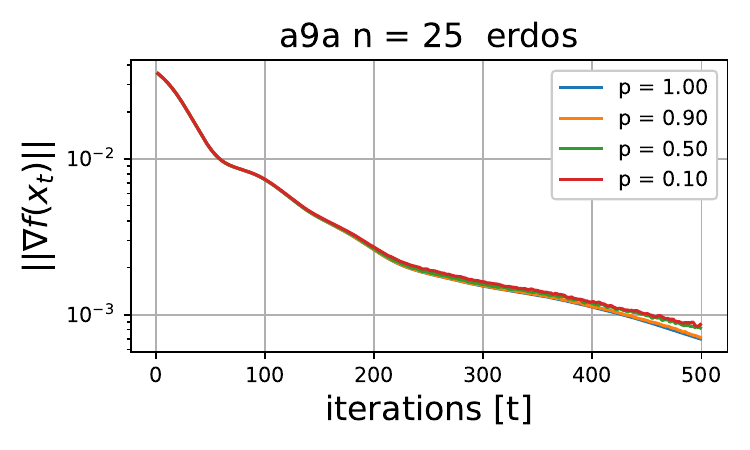}
\includegraphics[width=0.245\textwidth]{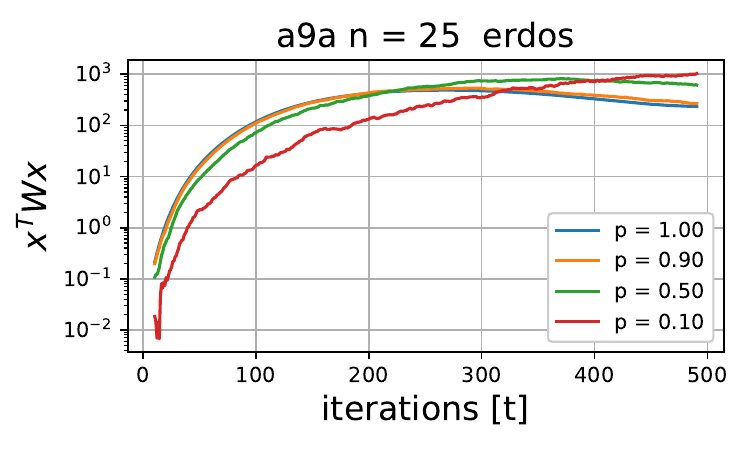}
\includegraphics[width=0.245\textwidth]{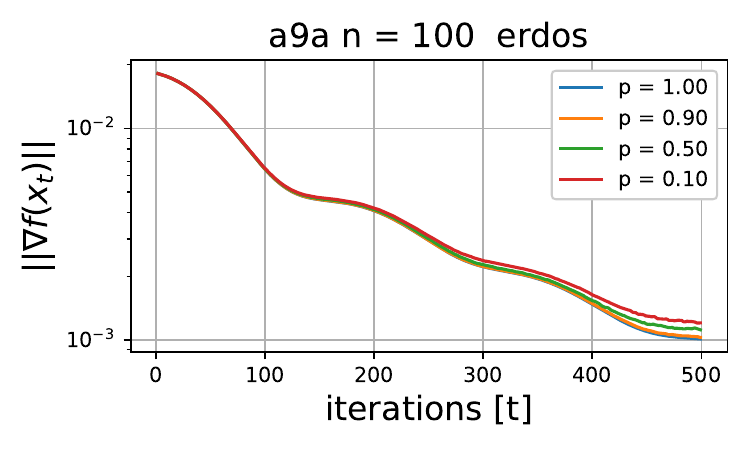}
\includegraphics[width=0.245\textwidth]{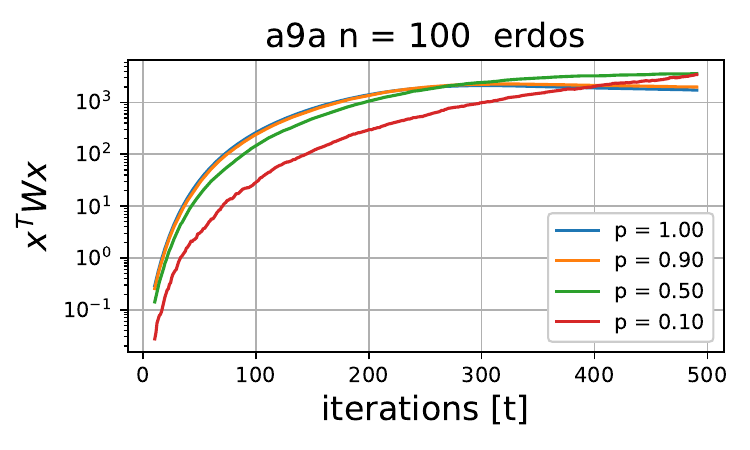}

\includegraphics[width=0.245\textwidth]{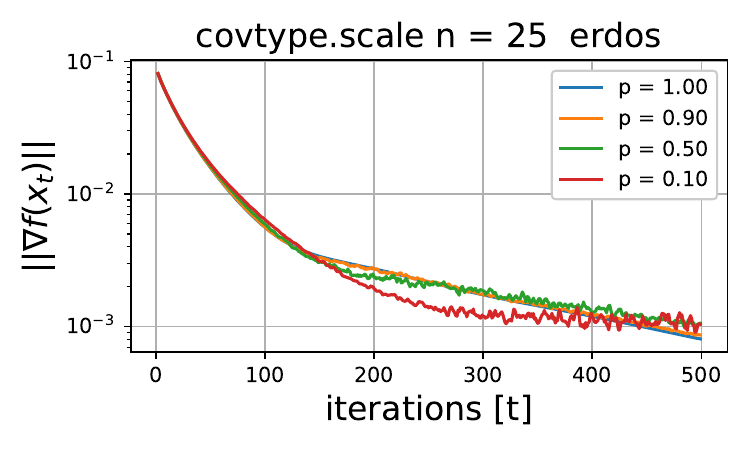}
\includegraphics[width=0.245\textwidth]{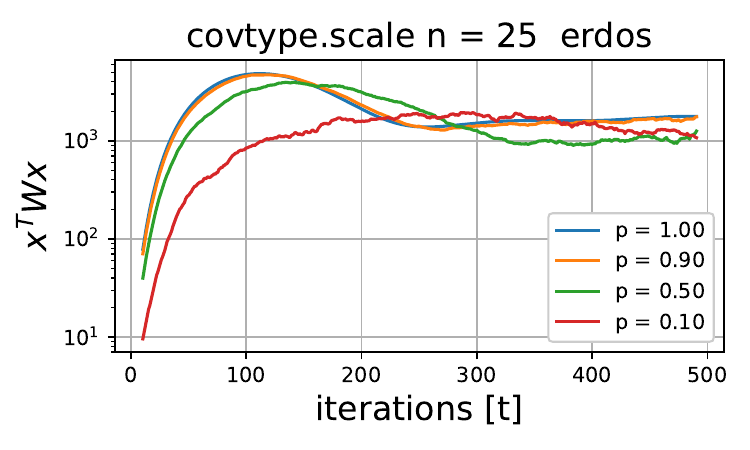}
\includegraphics[width=0.245\textwidth]{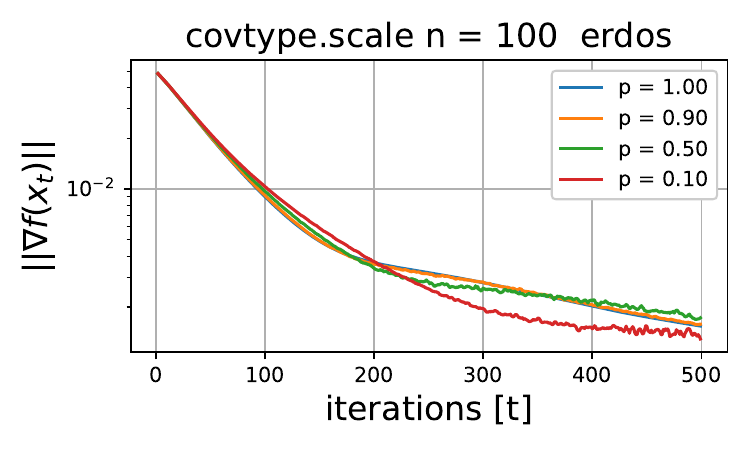}
\includegraphics[width=0.245\textwidth]{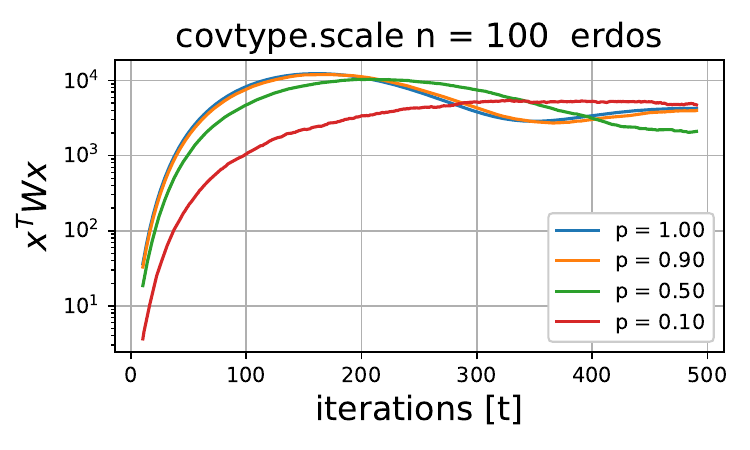}

\includegraphics[width=0.245\textwidth]{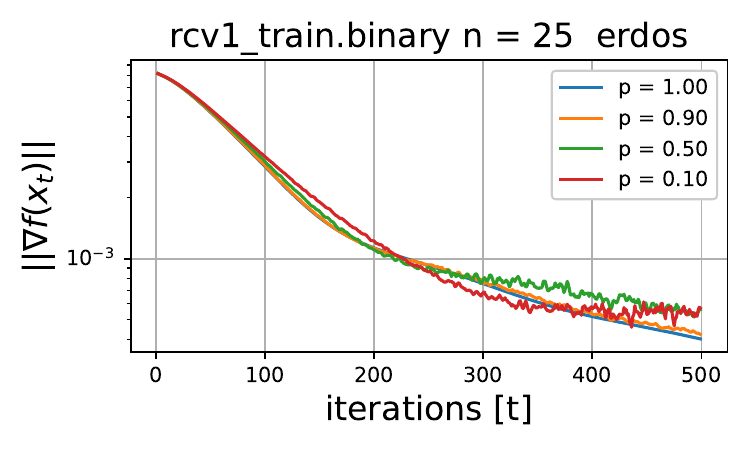}
\includegraphics[width=0.245\textwidth]{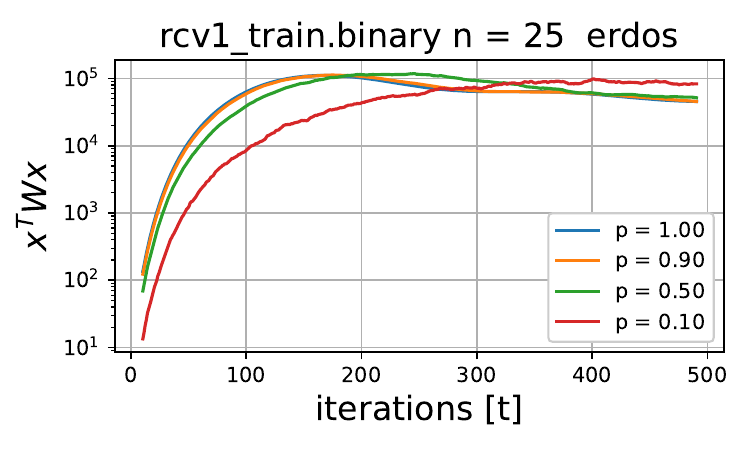}
\includegraphics[width=0.245\textwidth]{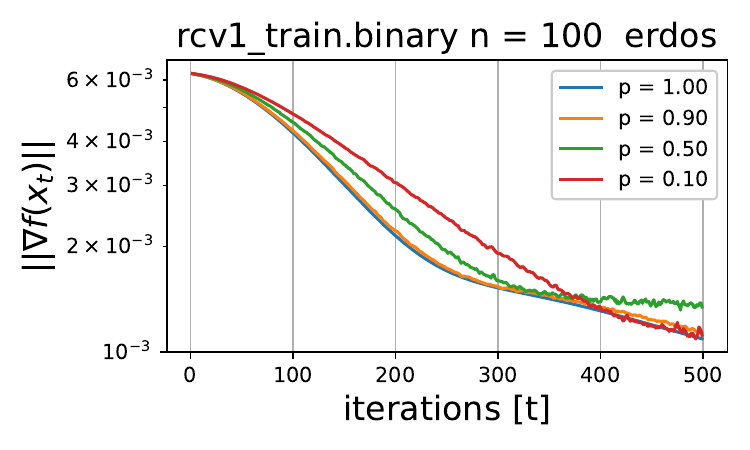}
\includegraphics[width=0.245\textwidth]{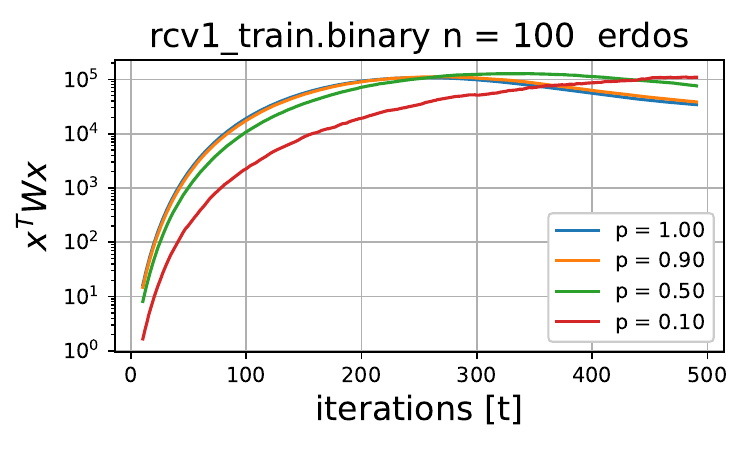}

    \caption{\ab{Evolution of $\|\nabla f(x_t)\|$ and $x_t^T W x_t$ for various datasets, number of devices $n \in \{25, 100\}$
    and different probability of keeping the communication edge $p$
    with the erdos network.}}
    \label{fig:droppingCommunication1}
\end{figure}

\begin{figure}
    \centering

\includegraphics[width=0.245\textwidth]{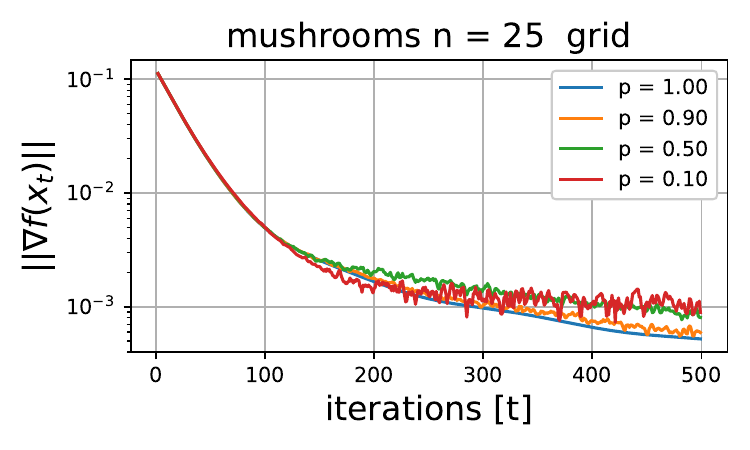}
\includegraphics[width=0.245\textwidth]{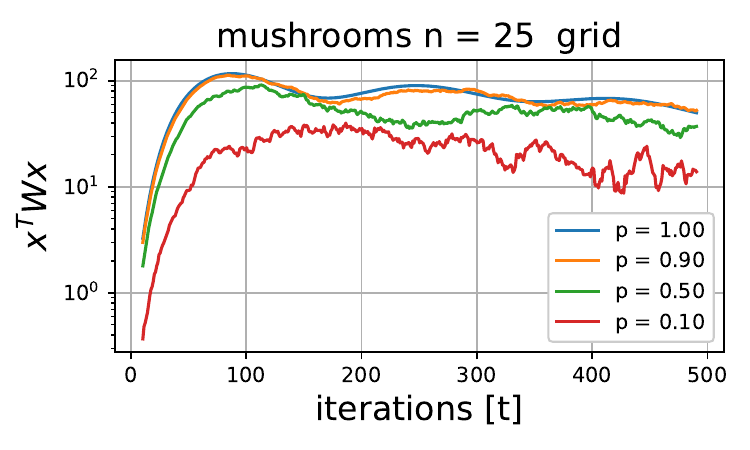}
\includegraphics[width=0.245\textwidth]{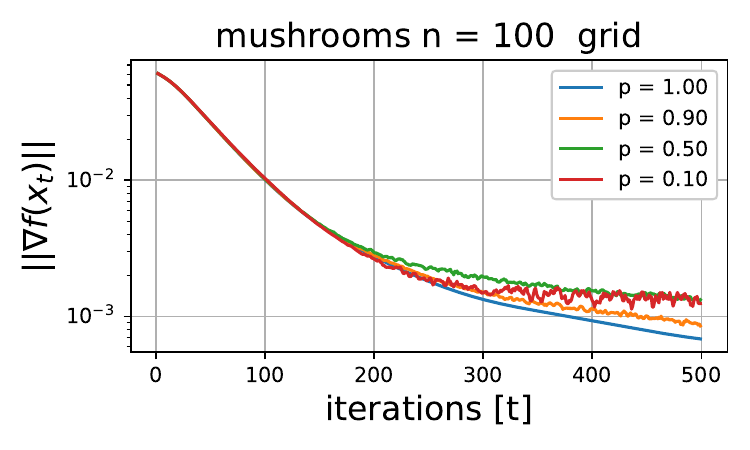}
\includegraphics[width=0.245\textwidth]{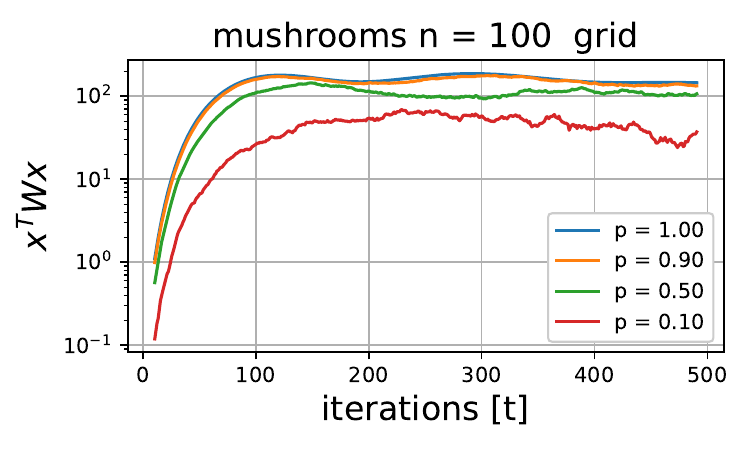}

\includegraphics[width=0.245\textwidth]{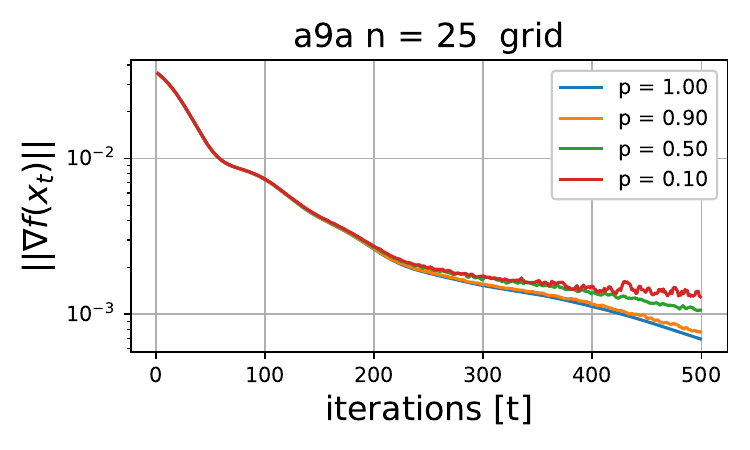}
\includegraphics[width=0.245\textwidth]{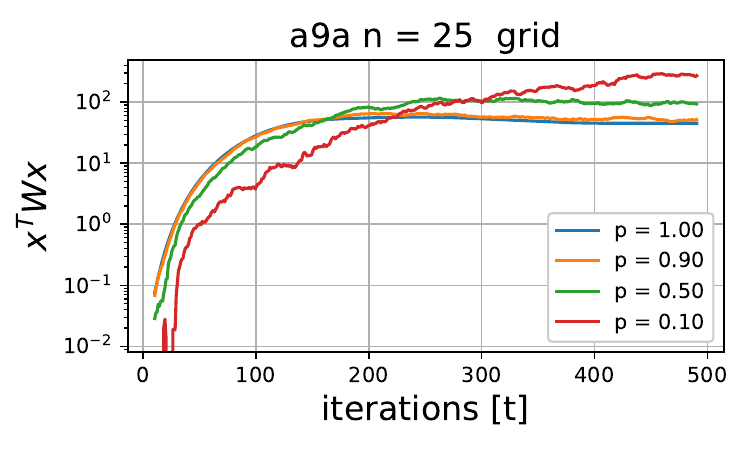}
\includegraphics[width=0.245\textwidth]{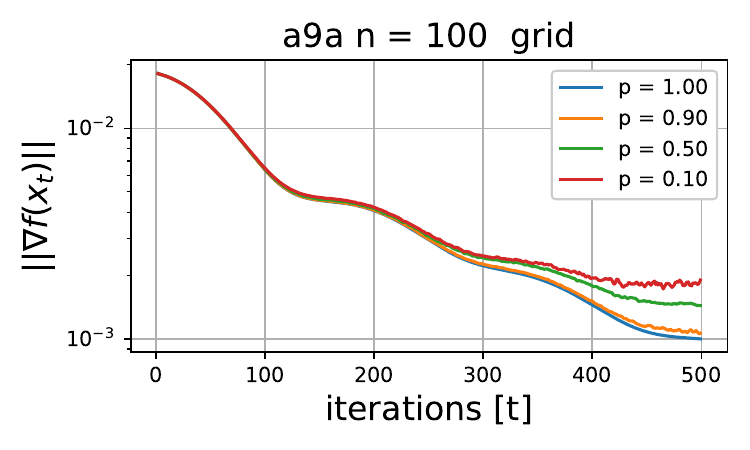}
\includegraphics[width=0.245\textwidth]{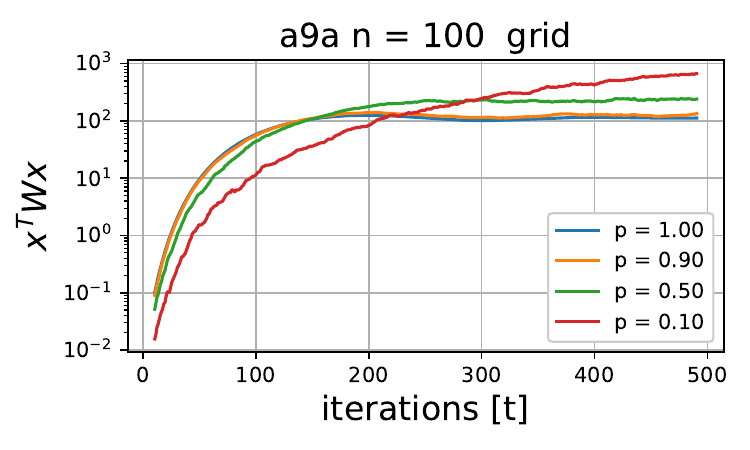}

\includegraphics[width=0.245\textwidth]{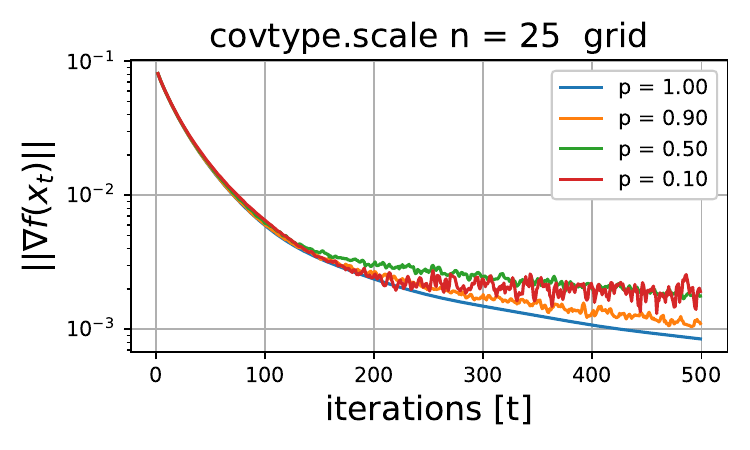}
\includegraphics[width=0.245\textwidth]{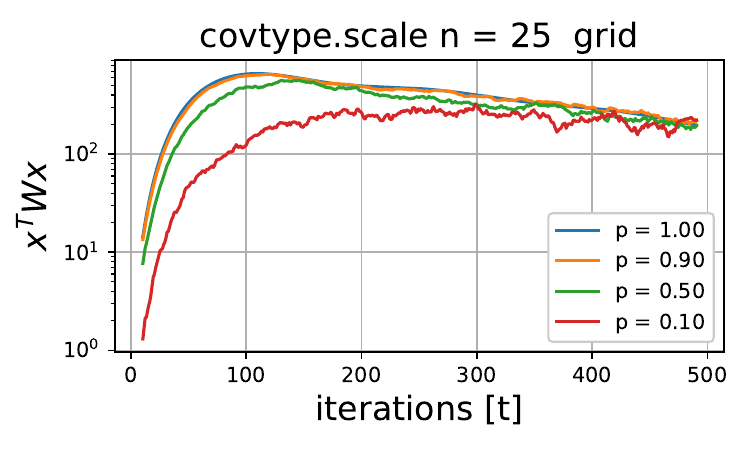}
\includegraphics[width=0.245\textwidth]{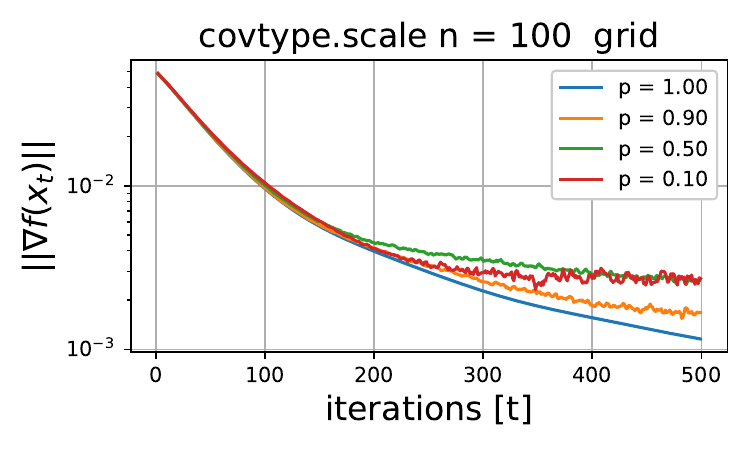}
\includegraphics[width=0.245\textwidth]{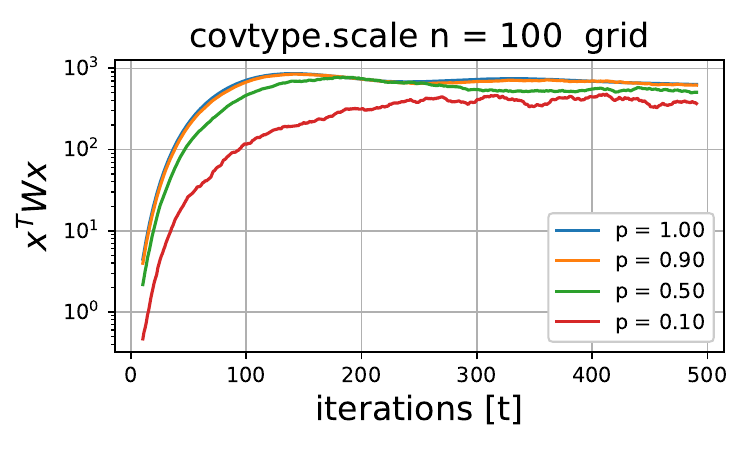}

\includegraphics[width=0.245\textwidth]{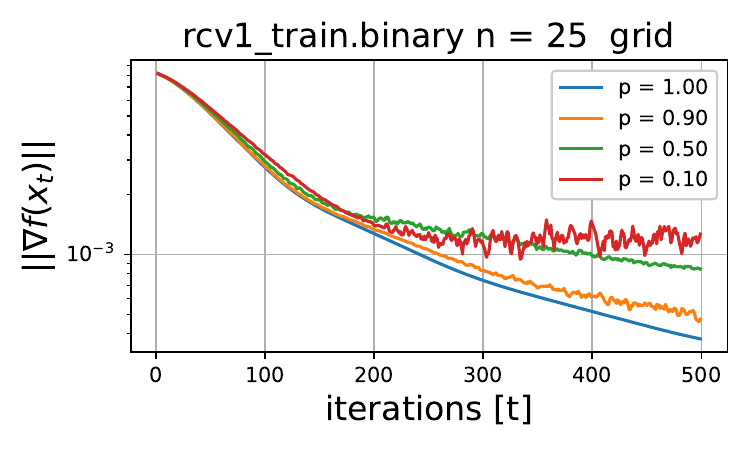}
\includegraphics[width=0.245\textwidth]{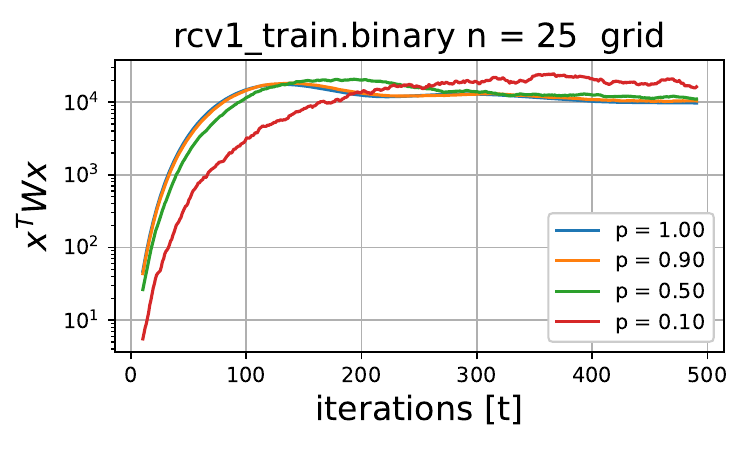}
\includegraphics[width=0.245\textwidth]{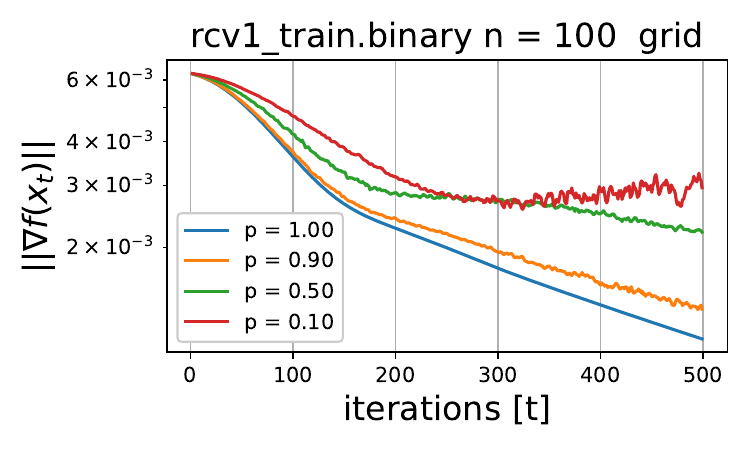}
\includegraphics[width=0.245\textwidth]{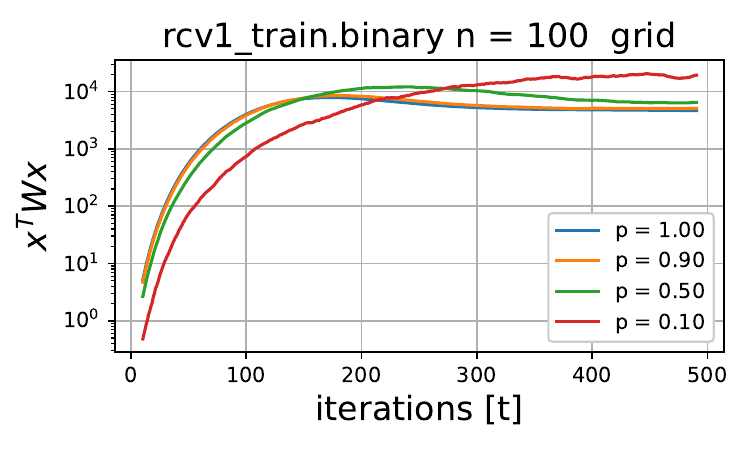}

    \caption{\ab{Evolution of $\|\nabla f(x_t)\|$ and $x_t^T W x_t$ for various datasets, number of devices $n \in \{25, 100\}$
    and different probability of keeping the communication edge $p$
    with grid network.}}
    \label{fig:droppingCommunication2}
\end{figure}

\begin{figure}
    \centering

\includegraphics[width=0.245\textwidth]{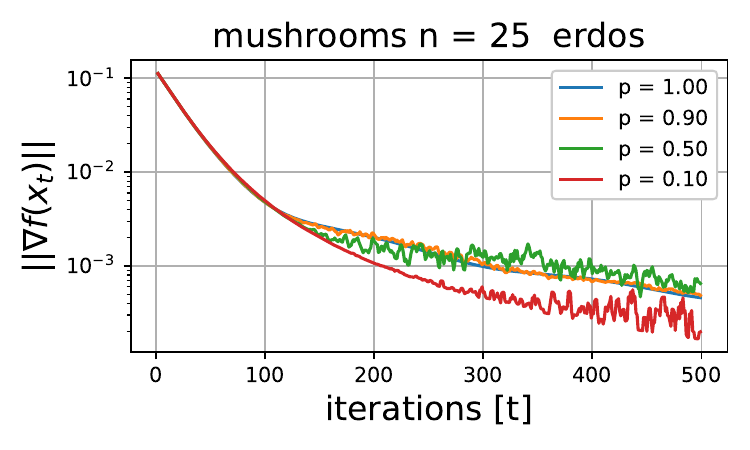}
\includegraphics[width=0.245\textwidth]{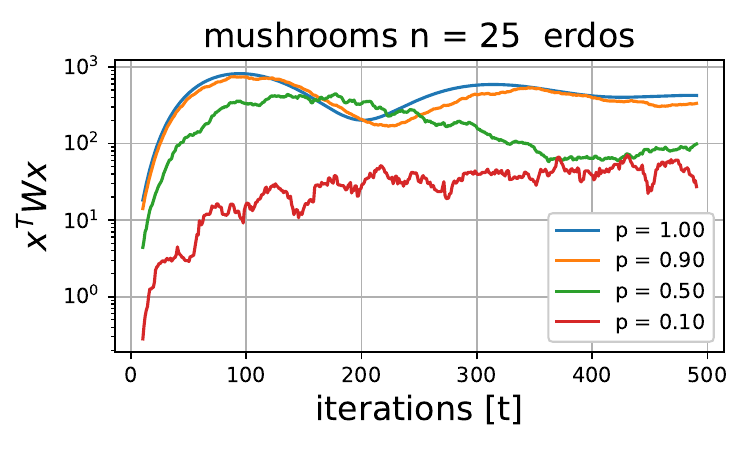}
\includegraphics[width=0.245\textwidth]{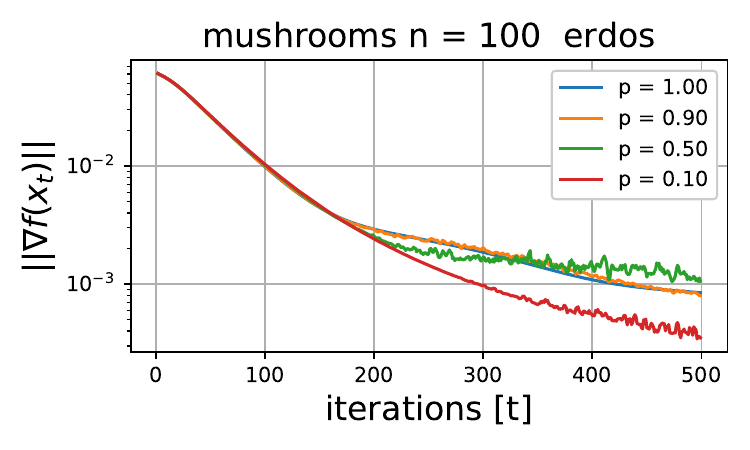}
\includegraphics[width=0.245\textwidth]{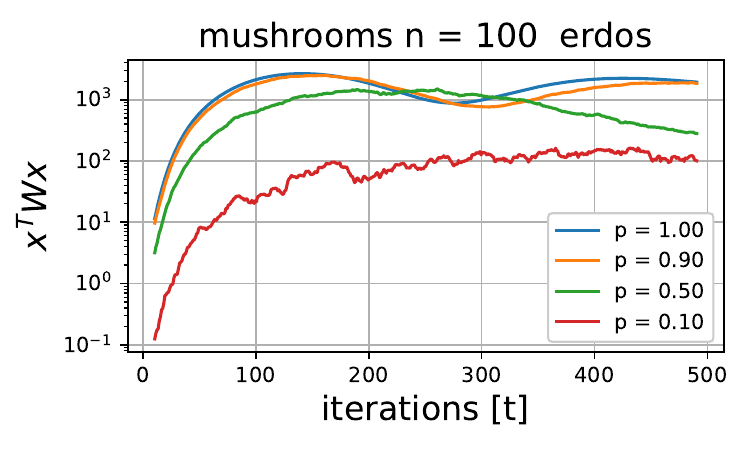}

\includegraphics[width=0.245\textwidth]{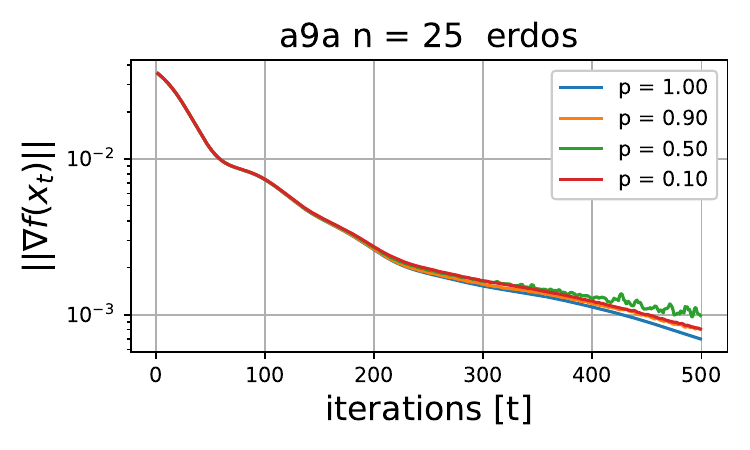}
\includegraphics[width=0.245\textwidth]{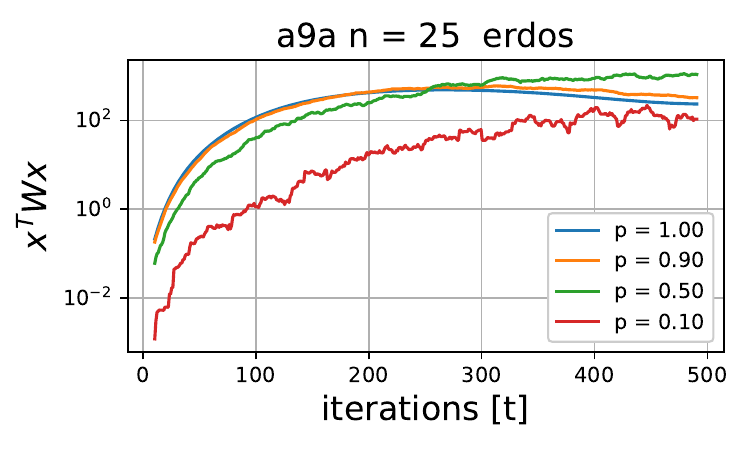}
\includegraphics[width=0.245\textwidth]{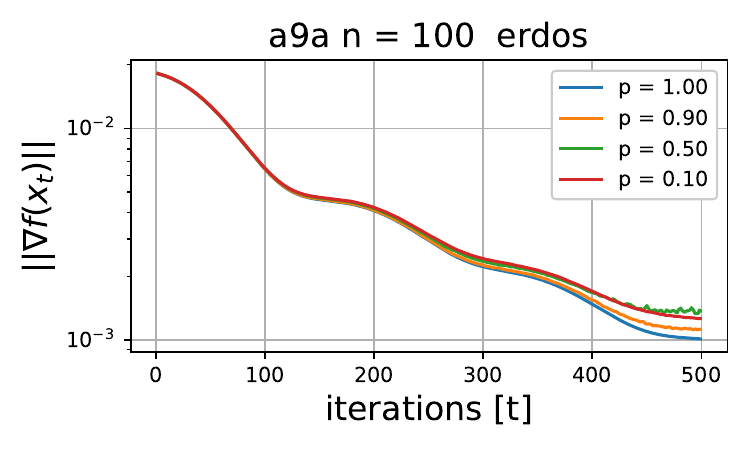}
\includegraphics[width=0.245\textwidth]{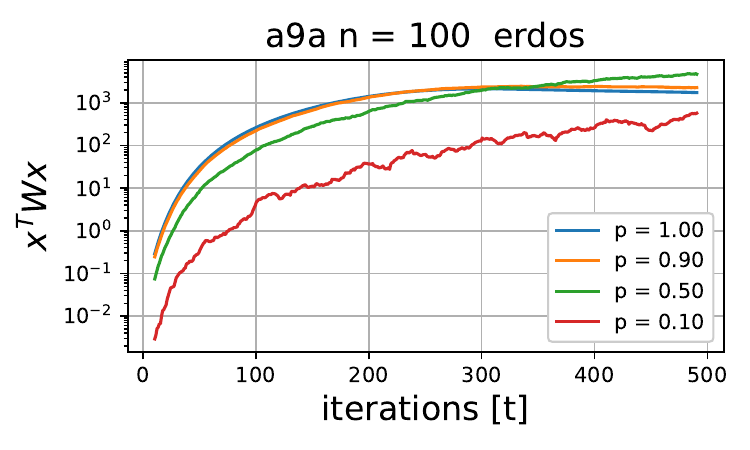}

\includegraphics[width=0.245\textwidth]{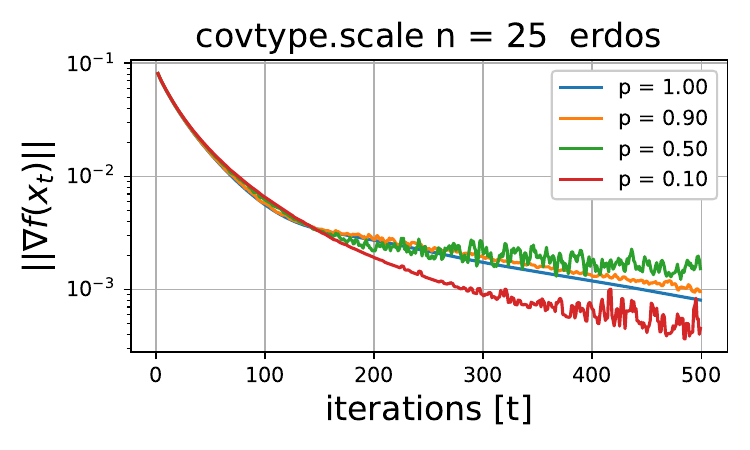}
\includegraphics[width=0.245\textwidth]{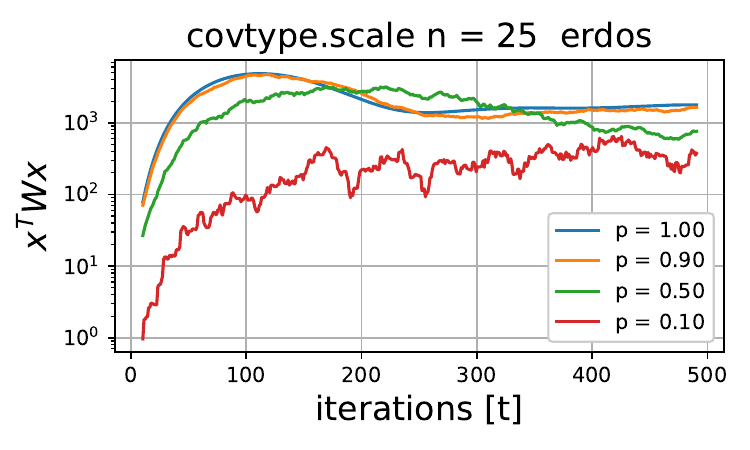}
\includegraphics[width=0.245\textwidth]{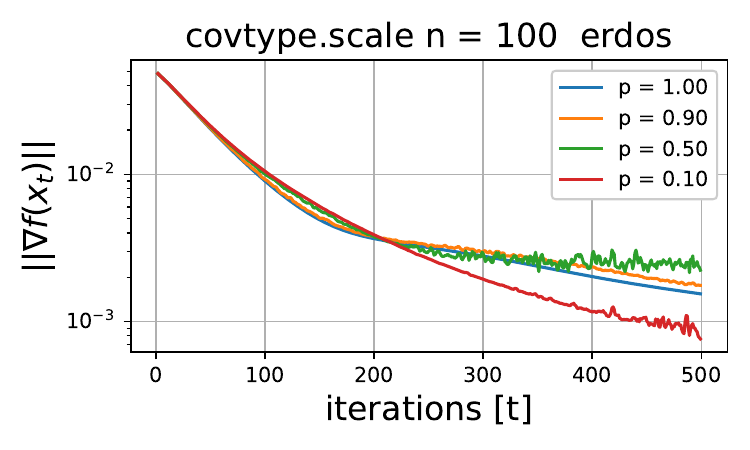}
\includegraphics[width=0.245\textwidth]{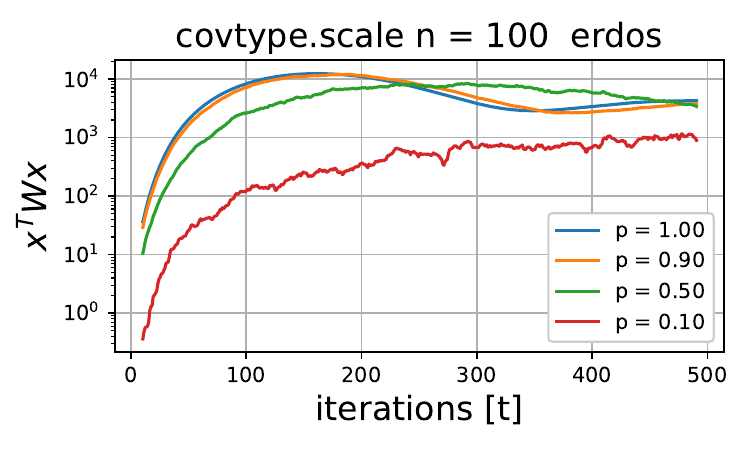}

\includegraphics[width=0.245\textwidth]{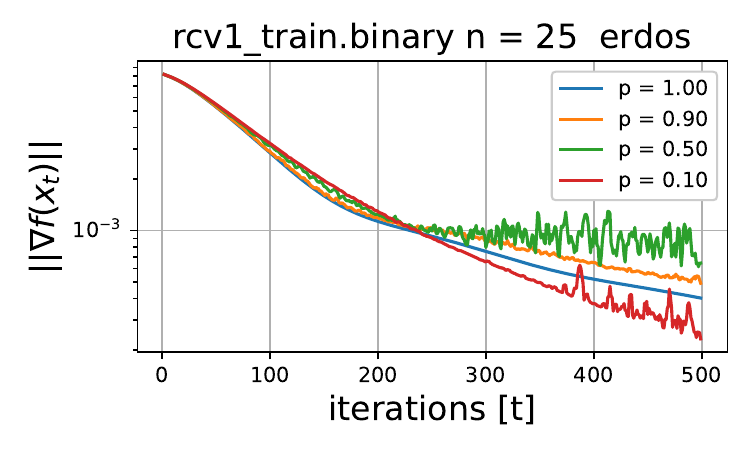}
\includegraphics[width=0.245\textwidth]{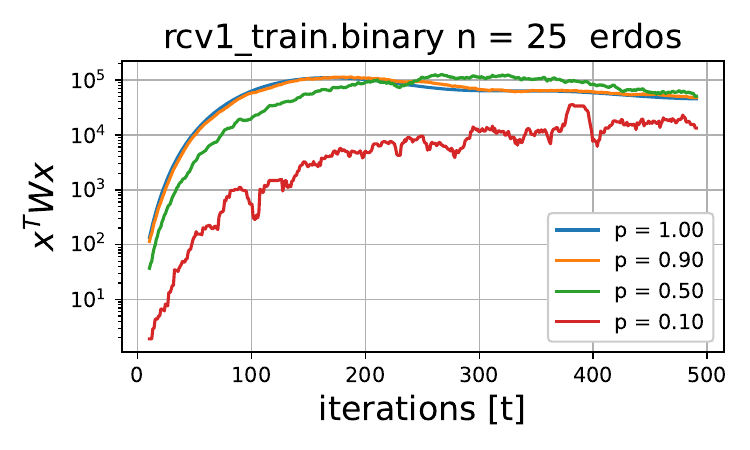}
\includegraphics[width=0.245\textwidth]{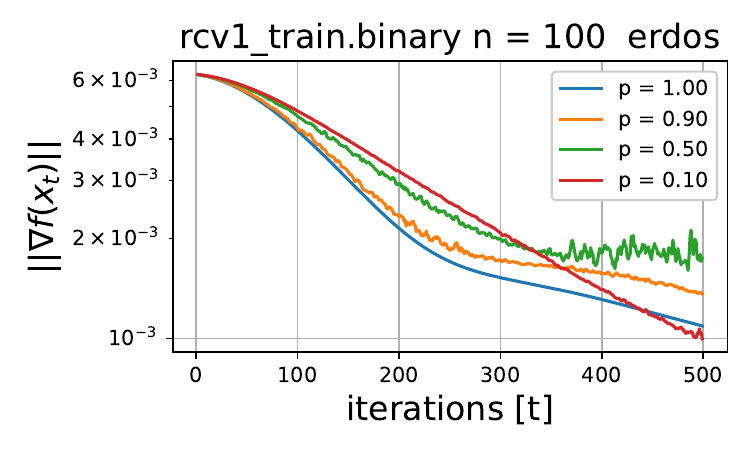}
\includegraphics[width=0.245\textwidth]{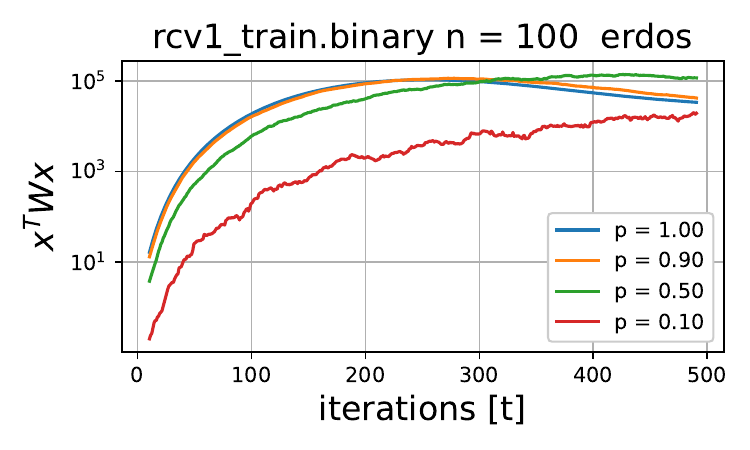}

    \caption{\ab{Evolution of $\|\nabla f(x_t)\|$ and $x_t^T W x_t$ for various datasets, number of devices $n \in \{25, 100\}$
    and different probability $p$ of keeping the device in the communication 
    with erdos network.}}
    \label{fig:droppingCommunication1B}
\end{figure}

\begin{figure}
    \centering

\includegraphics[width=0.245\textwidth]{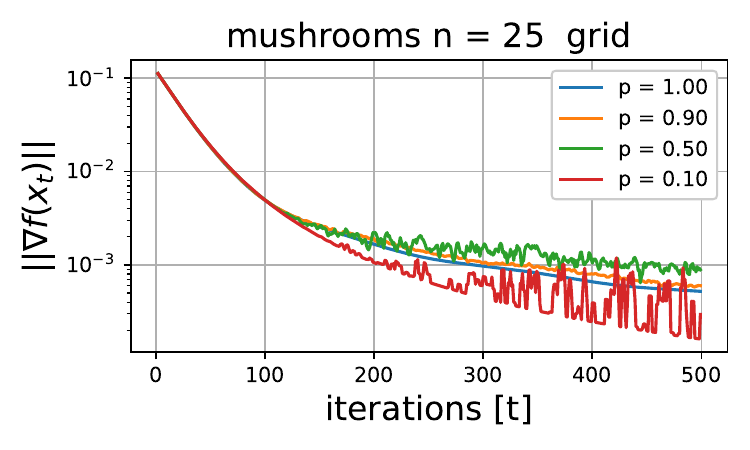}
\includegraphics[width=0.245\textwidth]{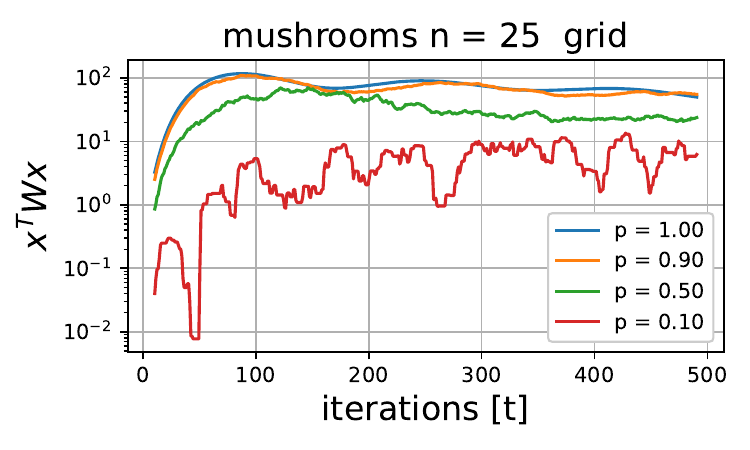}
\includegraphics[width=0.245\textwidth]{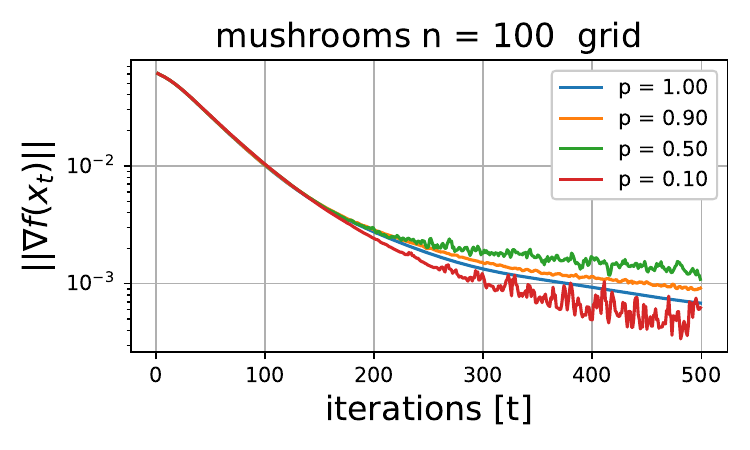}
\includegraphics[width=0.245\textwidth]{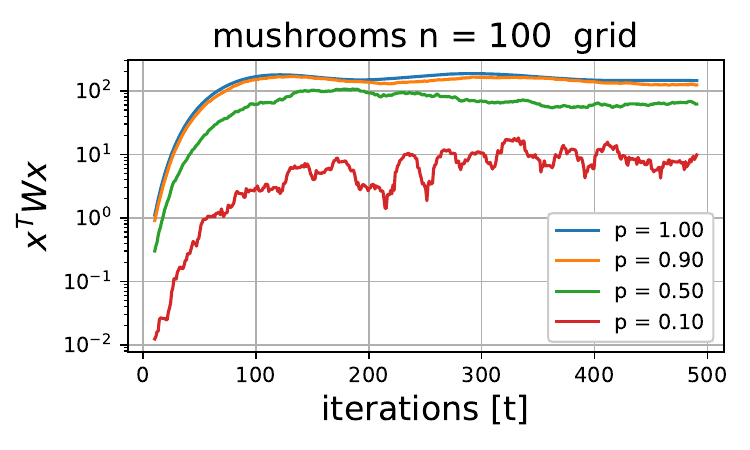}

\includegraphics[width=0.245\textwidth]{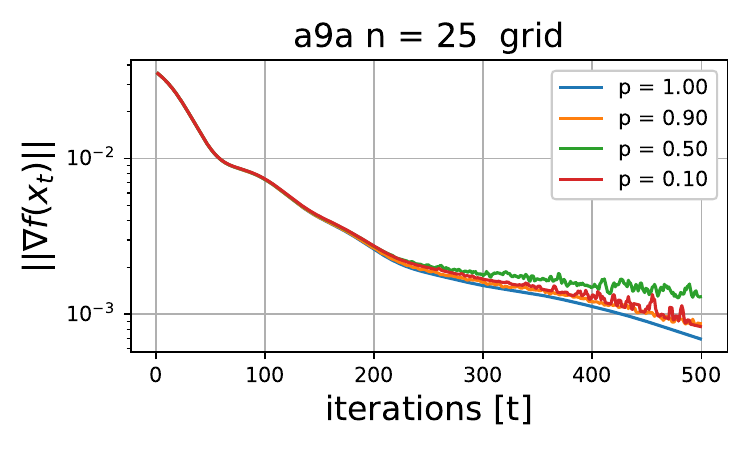}
\includegraphics[width=0.245\textwidth]{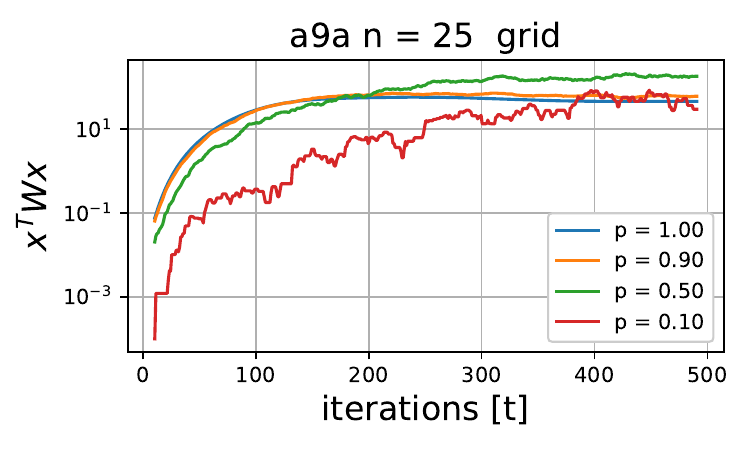}
\includegraphics[width=0.245\textwidth]{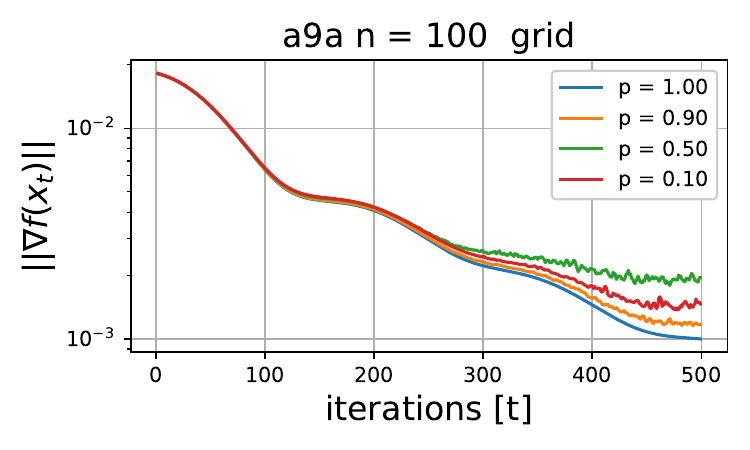}
\includegraphics[width=0.245\textwidth]{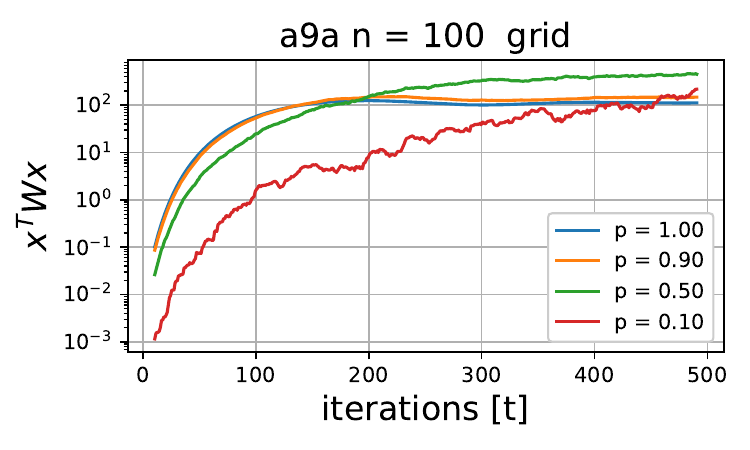}

\includegraphics[width=0.245\textwidth]{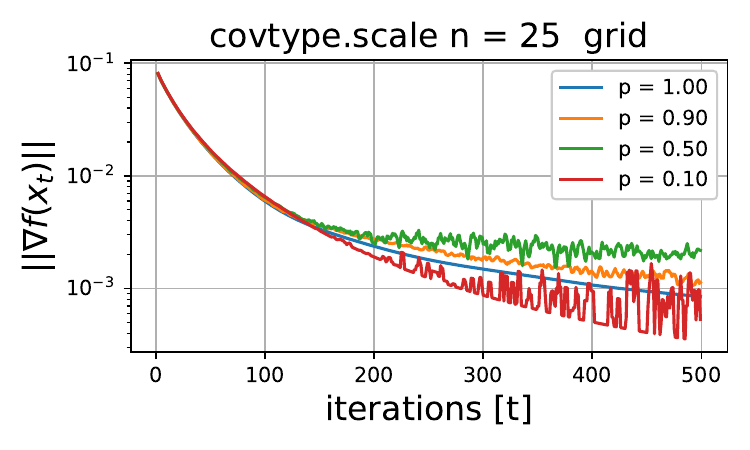}
\includegraphics[width=0.245\textwidth]{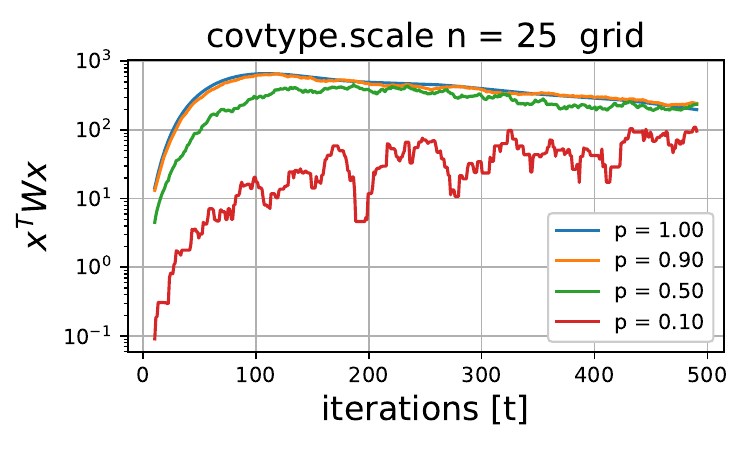}
\includegraphics[width=0.245\textwidth]{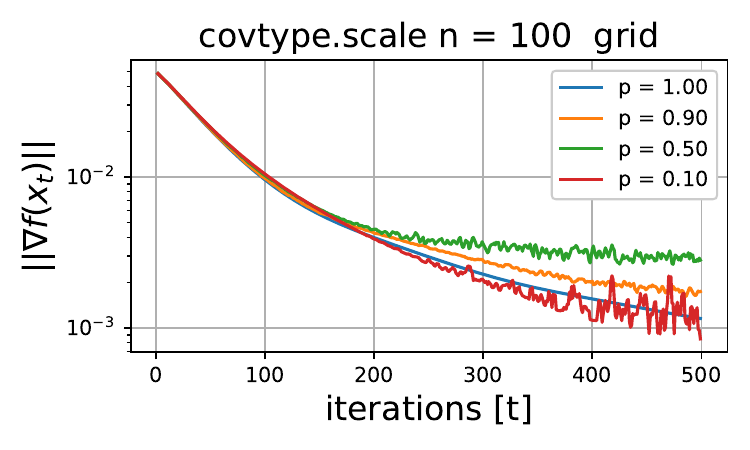}
\includegraphics[width=0.245\textwidth]{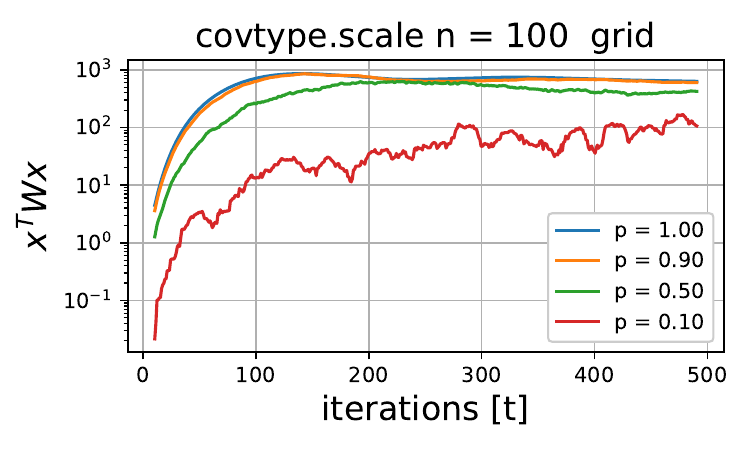}

\includegraphics[width=0.245\textwidth]{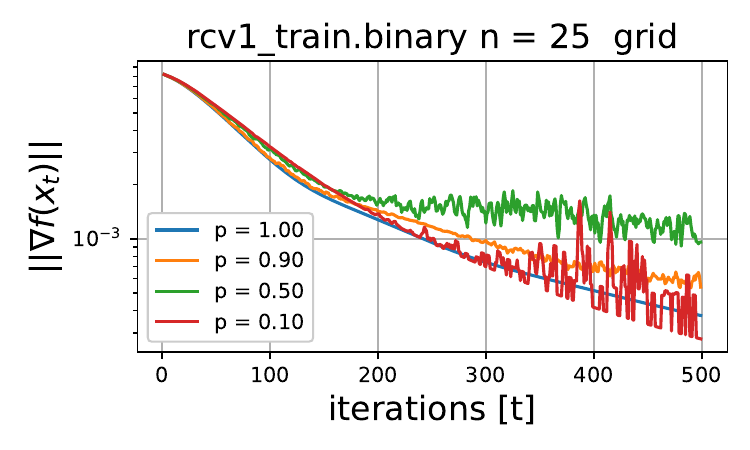}
\includegraphics[width=0.245\textwidth]{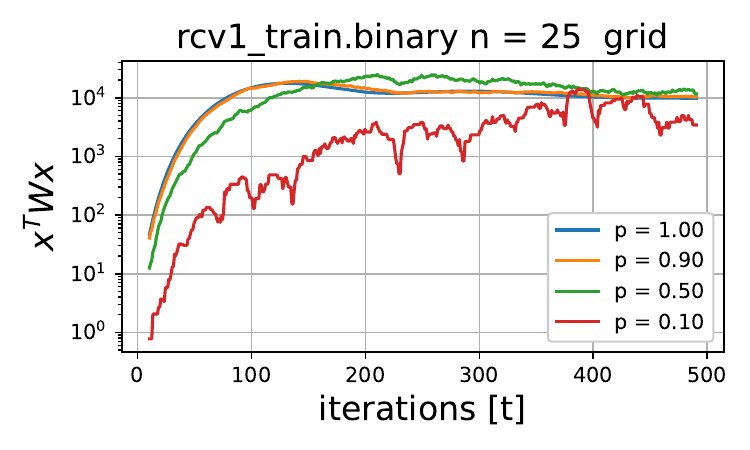}
\includegraphics[width=0.245\textwidth]{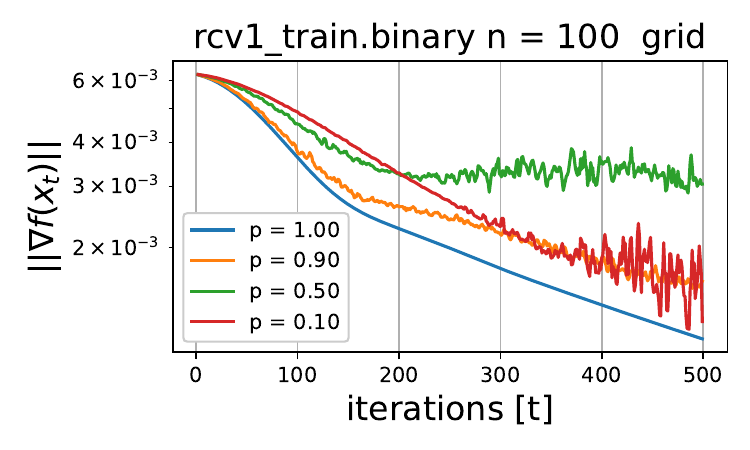}
\includegraphics[width=0.245\textwidth]{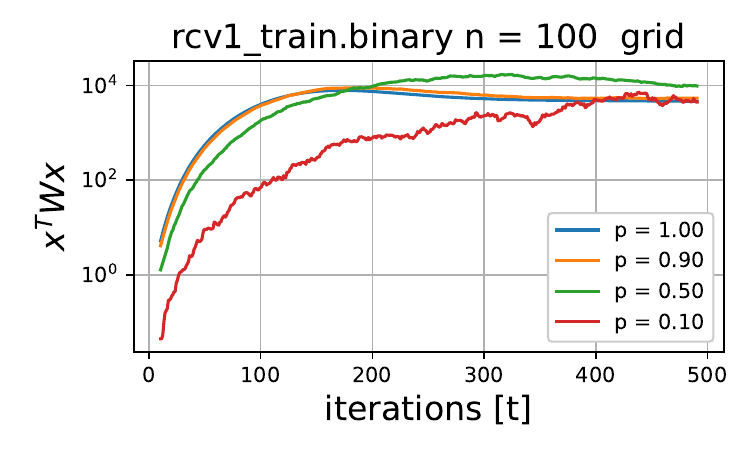}

    \caption{\ab{Evolution of $\|\nabla f(x_t)\|$ and $x_t^T W x_t$ for various datasets, number of devices $n \in \{25, 100\}$
    and different probability $p$ of keeping the device in the communication 
    with grid network.}}
    \label{fig:droppingCommunication2B}
\end{figure}

\ab{
\subsection{Partial worker participation}

One of the challenges of the FL setting is the fact that not all devices can always participate in all the communications \cite{kairouz2021advances,wang2021field}.
To simulate such a scenario, we conducted the following two experiments:
\begin{enumerate}
    \item {\bf Randomly dropping communication edge(s).}
For each iteration, and each communication edge $e$, the edge is kept with probability $p_e$, or dropped with probability $1-p_e$.
The result is that the gossip matrix $W$ is randomly modified at each iteration.
In Figures~\ref{fig:droppingCommunication1} and \ref{fig:droppingCommunication2}
we demonstrate empirically that keeping some communication edges with probability $p_e \in \{1.0, 0.9, 0.5, 0.1\}$ only mildly affects the convergence.
 \item {\bf Randomly dropping the device(s) from communication.}
In this case, a subset of devices is randomly selected. In particular, at each iteration, a device is kept with probability $p_d$, and excluded/dropped with probability $1-p_d$. As before, the effect is that the gossip matrix $W$ is randomly modified at each iteration. In Figures~\ref{fig:droppingCommunication1B} and \ref{fig:droppingCommunication2B} we demonstrate empirically that keeping only some devices with probability $p_d \in \{1.0, 0.9, 0.5, 0.1\}$ only mildly affects the convergence.
\end{enumerate}

}

\begin{figure}
    \centering
\includegraphics[width=0.325\textwidth]{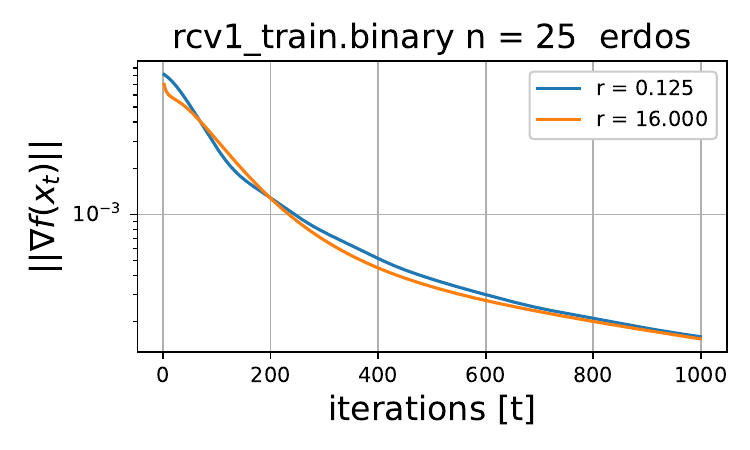}
\includegraphics[width=0.325\textwidth]{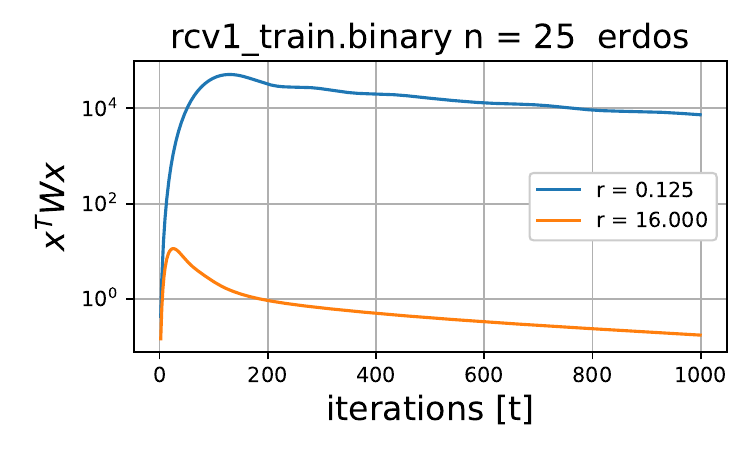}
\includegraphics[width=0.325\textwidth]{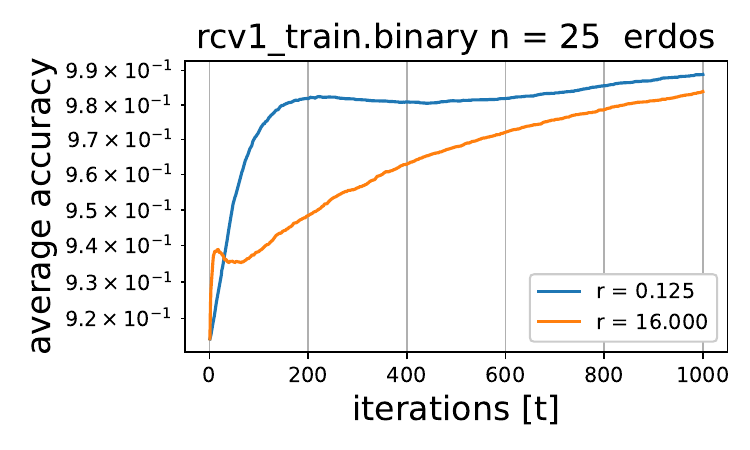}

\includegraphics[width=0.325\textwidth]{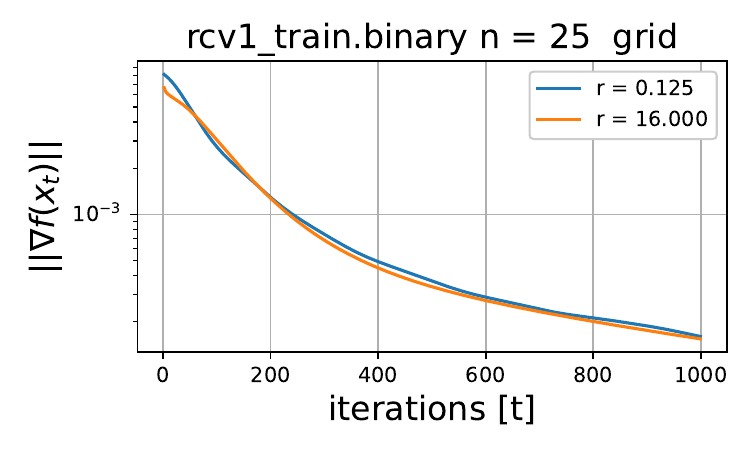}
\includegraphics[width=0.325\textwidth]{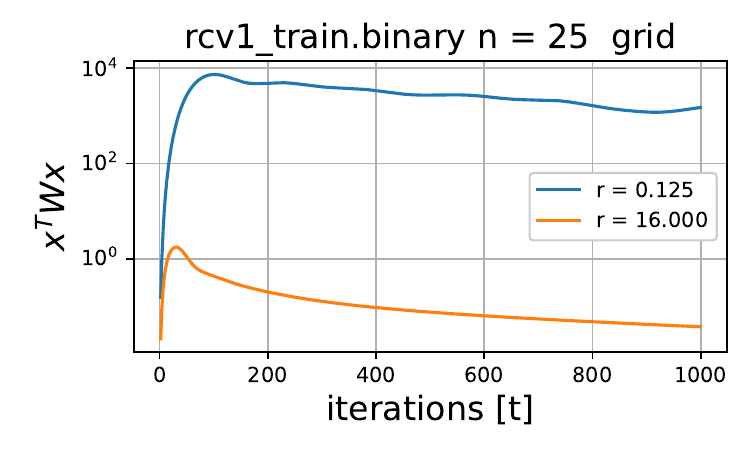}
\includegraphics[width=0.325\textwidth]{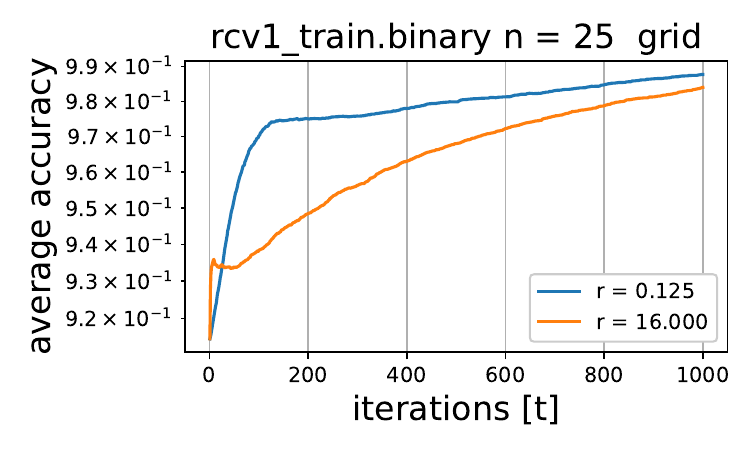}

    \caption{\ab{
    Evolution of $\|\nabla f(x_t)\|$, $x_t^T W x_t$ and
    average accuracy of local models
    for rcv1 datasets
    and erdos and grid network.
    We compare with two levels of regularization: low with
    $\lambda = 0.125 \cdot  \frac{L}{\lambda_{\max}(W)} $
    and high with $\lambda = 16 \cdot \frac{L}{\lambda_{\max}(W)}$.
    Note that, although both regularization values give comparable $\|\nabla f(x_t)\|$, the average accuracy for the case with smaller penalization is better.
    }}
    \label{fig:benefitOfpersonalization}
\end{figure}

\ab{
\subsection{The benefit of personalized training}
In Section~\ref{sec:variousSmallr} we discussed the case when 
$\lambda \leq \frac{L}{\lambda_{\max}(W)}$ 
 that allows for more personalization of local models. 
Note that, as discussed in  Section~\ref{sec:analys},
we use Algorithm~\ref{meta_algorithm}
with different settings for $h_1(x)$ and $h_2(x)$ depending on the value of $\lambda$.
In Figure~\ref{fig:benefitOfpersonalization} we investigate the behavior of Algorithm~\ref{meta_algorithm} for $\lambda = r \frac{L}{\lambda_{\max}(W)}$ with $r \in \{ 0.125, 16 \}$. Note that a larger value of $\lambda$ ($r$) corresponds to larger penalization if the model deviates from the mean ($x^T W x$); consequently, this allows less personalization.
}

\section{Conclusion}
In this work we studied the problem of decentralized personalized federated learning. Problem \eqref{rewrite_main_problem} used a penalty term that was based upon the specific network structure, which was more appropriate than a `deviation from the average' penalty in the decentralized setting. We presented lower bounds on the local communication and computation costs, and we presented algorithms that achieved these lower bounds. Numerical experiments demonstrated the benefits of this approach.

\ab{Interesting issues for further research are those related to the more practical features arising in a federated learning setup, including asynchronous and delayed transmissions, and compression of information to reduce communication cost, among others. It would also be interesting to perform numerical experiments using the Leaf framework (\texttt{https://leaf.cmu.edu}).}

\section*{Acknowledgements}
 This work was partially conducted while A. Sadiev, A. Beznosikov, D.Dvinskikh  were visiting research assistants and A. Gasnikov was a visiting scholar in Mohamed bin Zayed University of Artificial Intelligence (MBZUAI). 
 
 This work was supported by a grant for research centers in the field of artificial intelligence, provided by the Analytical Center for the Government of the Russian Federation in accordance with the subsidy agreement (agreement identifier 000000D730321P5Q0002) and the agreement with the Moscow Institute of Physics and Technology dated November 1, 2021 No. 70-2021-00138.

%  We consider the following  functions $f_j$.
% \begin{equation}
%     f_j(x) = \log \left(1 + \exp\left(-y_j \left\la a_j, x\right\ra\right)\right),
% \end{equation}
% where $y_j$ takes only $1$ and $-1$ value with equal probability, $a_k$ denotes the $j$th column of the matrix $A$, where $A$ is generated randomly. 
% The experiment are performed on the cycle graph which has
% the following Laplacian matrix 
% %\todo{let's make it smaller}
% % \begin{eqnarray*}
% %     W = \left(
% %     \begin{array}{cccccccc}
% %     1&-1 & & & & & &  \\
% %     -1&2 &-1 & & & & &  \\
% %     &-1 &2 & -1 & & & & \\
% %     & & -1&2 &-1 & & & \\
% %     & & &-1 &2 &-1 & & \\
% %     & & & & &\ldots & & \\
% %     & & & & &-1 &2 &-1 \\
% %     & & & & & &-1 &1 \\
% %     \end{array}
% %     \right).
% % \end{eqnarray*}
% \begin{eqnarray}\label{U}
%     U = 
%     \begin{pmatrix}
%         1 & -1 & & &   \\
%         -1 & 2 & -1 &  &  \\
%         & \ddots&\ddots & \ddots&  \\
%         & &   -1 & 2 & -1\\
%         & &    & -1 & 1\\
%     \end{pmatrix},
% \end{eqnarray}
% and we take $\hat W = U$. It can be shown that this $\hat W$ satisfies Definition~\ref{gossipmatrix}, i.e., it is a gossip matrix.
 
%\clearpage
 
\bibliographystyle{plain}
\bibliography{ref}

\begin{thebibliography}{10}

\bibitem{allen2017katyusha}
Zeyuan Allen-Zhu.
\newblock Katyusha: The first direct acceleration of stochastic gradient
  methods.
\newblock {\em The Journal of Machine Learning Research}, 18(1):8194--8244,
  2017.

\bibitem{anderson1985eigenvalues}
William~N Anderson~Jr and Thomas~D Morley.
\newblock Eigenvalues of the laplacian of a graph.
\newblock {\em Linear and multilinear algebra}, 18(2):141--145, 1985.

\bibitem{beznosikov2021decentralized}
Aleksandr Beznosikov, Vadim Sushko, Abdurakhmon Sadiev, and Alexander Gasnikov.
\newblock Decentralized personalized federated min-max problems.
\newblock {\em arXiv preprint arXiv:2106.07289}, 2021.

\bibitem{boyd2006randomized}
Stephen Boyd, Arpita Ghosh, Balaji Prabhakar, and Devavrat Shah.
\newblock Randomized gossip algorithms.
\newblock {\em IEEE transactions on information theory}, 52(6):2508--2530,
  2006.

\bibitem{chang2011libsvm}
Chih-Chung Chang.
\newblock " libsvm: a library for support vector machines," acm transactions on
  intelligent systems and technology, 2: 27: 1--27: 27, 2011.
\newblock {\em http://www. csie. ntu. edu. tw/\~{} cjlin/libsvm}, 2, 2011.

\bibitem{gasnikov2020AM}
Alexander Gasnikov, Darina Dvinskikh, Pavel Dvurechensky, Dmitry Kamzolov,
  Vladislav Matykhin, Dmitry Pasechnyk, Nazarii Tupitsa, and Alexei Chernov.
\newblock Accelerated meta-algorithm for convex optimization.
\newblock {\em arXiv preprint arXiv:2004.08691}, 2020.

\bibitem{gorbunov2019optimal}
Eduard Gorbunov, Darina Dvinskikh, and Alexander Gasnikov.
\newblock Optimal decentralized distributed algorithms for stochastic convex
  optimization.
\newblock {\em arXiv preprint arXiv:1911.07363}, 2019.

\bibitem{hanzely2020lower}
Filip Hanzely, Slavom{\'\i}r Hanzely, Samuel Horv{\'a}th, and Peter
  Richt{\'a}rik.
\newblock Lower bounds and optimal algorithms for personalized federated
  learning.
\newblock {\em arXiv preprint arXiv:2010.02372}, 2020.

\bibitem{L_Katyusha}
Filip Hanzely, Dmitry Kovalev, and Peter Richtarik.
\newblock Variance reduced coordinate descent with acceleration: New method
  with a surprising application to finite-sum problems.
\newblock {\em arXiv preprint arXiv:2002.04670}, Feb 2020.

\bibitem{hanzely2020federated}
Filip Hanzely and Peter Richt{\'a}rik.
\newblock Federated learning of a mixture of global and local models.
\newblock {\em arXiv preprint arXiv:2002.05516}, 2020.

\bibitem{hanzely2021personalized}
Filip Hanzely, Boxin Zhao, and Mladen Kolar.
\newblock Personalized federated learning: A unified framework and universal
  optimization techniques, 2021.

\bibitem{Hendrikx2020}
Hadrien Hendrikx, Francis Bach, and Laurent Massoulie.
\newblock An optimal algorithm for decentralized finite sum optimization.
\newblock {\em arXiv preprint arXiv:2005.10675}, 2020.

\bibitem{NIPS2013_ac1dd209}
Rie Johnson and Tong Zhang.
\newblock Accelerating stochastic gradient descent using predictive variance
  reduction.
\newblock In C.J. Burges, L.~Bottou, M.~Welling, Z.~Ghahramani, and K.Q.
  Weinberger, editors, {\em Advances in Neural Information Processing Systems},
  volume~26. Curran Associates, Inc., 2013.

\bibitem{kairouz2021advances}
Peter Kairouz, H~Brendan McMahan, Brendan Avent, Aur{\'e}lien Bellet, Mehdi
  Bennis, Arjun~Nitin Bhagoji, Kallista Bonawitz, Zachary Charles, Graham
  Cormode, Rachel Cummings, et~al.
\newblock Advances and open problems in federated learning.
\newblock {\em Foundations and Trends{\textregistered} in Machine Learning},
  14(1--2):1--210, 2021.

\bibitem{konevcny2016federated}
Jakub Kone{\v{c}}n{\`y}, H~Brendan McMahan, Daniel Ramage, and Peter
  Richt{\'a}rik.
\newblock Federated optimization: Distributed machine learning for on-device
  intelligence.
\newblock {\em arXiv preprint arXiv:1610.02527}, 2016.

\bibitem{kulkarni2020survey}
Viraj Kulkarni, Milind Kulkarni, and Aniruddha Pant.
\newblock Survey of personalization techniques for federated learning.
\newblock In {\em 2020 Fourth World Conference on Smart Trends in Systems,
  Security and Sustainability (WorldS4)}, pages 794--797. IEEE, 2020.

\bibitem{li2020decentralized}
Huan Li, Cong Fang, Wotao Yin, and Zhouchen Lin.
\newblock Decentralized accelerated gradient methods with increasing penalty
  parameters.
\newblock {\em IEEE Transactions on Signal Processing}, 68:4855--4870, 2020.

\bibitem{mcmahan2017communication}
Brendan McMahan, Eider Moore, Daniel Ramage, Seth Hampson, and Blaise~Aguera
  y~Arcas.
\newblock Communication-efficient learning of deep networks from decentralized
  data.
\newblock In {\em Artificial Intelligence and Statistics}, pages 1273--1282.
  PMLR, 2017.

\bibitem{nedic2009distributed}
Angelia Nedic and Asuman Ozdaglar.
\newblock Distributed subgradient methods for multi-agent optimization.
\newblock {\em IEEE Transactions on Automatic Control}, 54(1):48--61, 2009.

\bibitem{nesterov2003introductory}
Yurii Nesterov.
\newblock {\em Introductory lectures on convex optimization: A basic course},
  volume~87.
\newblock Springer Science \& Business Media, 2003.

\bibitem{scaman2017optimal}
Kevin Scaman, Francis Bach, S{\'e}bastien Bubeck, Yin~Tat Lee, and Laurent
  Massouli{\'e}.
\newblock Optimal algorithms for smooth and strongly convex distributed
  optimization in networks.
\newblock {\em arXiv preprint arXiv:1702.08704}, 2017.

\bibitem{shalev2014understanding}
Shai Shalev-Shwartz and Shai Ben-David.
\newblock {\em Understanding machine learning: From theory to algorithms}.
\newblock Cambridge university press, 2014.

\bibitem{smith2017federated}
Virginia Smith, Chao-Kai Chiang, Maziar Sanjabi, and Ameet Talwalkar.
\newblock Federated multi-task learning.
\newblock {\em arXiv preprint arXiv:1705.10467}, 2017.

\bibitem{tominin2021accelerated}
Vladislav Tominin, Yaroslav Tominin, Ekaterina Borodich, Dmitry Kovalev,
  Alexander Gasnikov, and Pavel Dvurechensky.
\newblock On accelerated methods for saddle-point problems with composite
  structure.
\newblock {\em arXiv preprint arXiv:2103.09344}, 2021.

\bibitem{wang2021field}
Jianyu Wang, Zachary Charles, Zheng Xu, Gauri Joshi, H~Brendan McMahan, Maruan
  Al-Shedivat, Galen Andrew, Salman Avestimehr, Katharine Daly, Deepesh Data,
  et~al.
\newblock A field guide to federated optimization.
\newblock {\em arXiv preprint arXiv:2107.06917}, 2021.

\bibitem{wang2018distributed}
Weiran Wang, Jialei Wang, Mladen Kolar, and Nathan Srebro.
\newblock Distributed stochastic multi-task learning with graph regularization.
\newblock {\em arXiv preprint arXiv:1802.03830}, 2018.

\end{thebibliography}
\newpage
%%%%%%%%%%%%%%%%%%%%%%%%%%%%%%%%%%%%%%%%%%%%%%%%%%%%%%%%%%%%
\appendix

\section{\ab{Versions of Algorithms \ref{meta_algorithm} and \ref{alg_pfl_sum} for problem \eqref{rewrite_main_problem}}} \label{app:alg_prac}

\begin{algorithm} [h!]
\caption{\texttt{MA} for $\lambda\eigWmax \geq L$.}
\label{meta_algorithm_big_lambda}
\hspace*{0.2cm} {\bf Input:} starting point $x^0_i = x^0 \in \R^{d}$, no. of iterations $K$,  parameter $\gamma=2L>0$, accuracy $\delta > 0$\\
\hspace*{0.3cm}{\bf Initialization:} $A^0 = 0,\quad y^0_i = x^0_i, \quad \tau = \frac{1}{2\gamma}$
	\begin{algorithmic}
		\For{ $k = 0, \ldots, K-1$}
		\State $ a^{k+1} = \displaystyle \frac{\tau+\sqrt{\tau^2+4\tau A^k}}{2}$
		\State $A^{k+1} = A^k+a^{k+1}$
		\State Local update: $w^k_i = \displaystyle \frac{A^k}{A^{k + 1}}y^k_i + \frac{a^{k+1}}{A^{k+1}} x^k_i
		$
		\State Local computation: $u^k_i = \nabla f_i (w_i^k)$
		\State Solve subproblem via gossip communications, i.e. find $\y^{k+1}\in \R^{nd}$, such that $\|\hat \y^{k+1} -  \y^{k+1}\|_2^2 \leq \delta$, where
		\begin{equation*}
		    \hat \y^{k+1} = \argmin_{\y\in \R^{nd}} \left\{ \langle \nabla h_1 (\w^k), \y-\w^k\rangle  + h_2(\y)   +\frac{\gamma}{2}\|\y-\w^k\|_2^{2} \right\}
		\end{equation*}
		\State Compute $z^k_i$ via gossip communication with neighbors: $\z^k =  \lambda W \y^{k+1}$
		\State Local update: $x^{k+1}_i = x^k_i-a^{k+1} (\nabla f_i(y^{k+1}_i) + z^k_i)$
		\EndFor
	\end{algorithmic}
	\hspace*{0.2cm} {\bf Output:} $\{y^{K}_i\}$
\end{algorithm}

\begin{algorithm} [h!]
\caption{\texttt{MA} for $\lambda\eigWmax < L$.}
\label{meta_algorithm_small_lambda}
\hspace*{0.2cm} {\bf Input:} starting point $x^0_i = x^0 \in \R^{d}$, no. of iterations $K$,  parameter $\gamma=2L>0$, accuracy $\delta > 0$\\
\hspace*{0.3cm}{\bf Initialization:} $A^0 = 0,\quad y^0_i = x^0_i, \quad \tau = \frac{1}{2\gamma}$
	\begin{algorithmic}
		\For{ $k = 0, \ldots, K-1$}
		\State $ a^{k+1} = \displaystyle \frac{\tau+\sqrt{\tau^2+4\tau A^k}}{2}$
		\State $A^{k+1} = A^k+a^{k+1}$
		\State Local update: $w^k_i = \displaystyle \frac{A^k}{A^{k + 1}}y^k_i + \frac{a^{k+1}}{A^{k+1}} x^k_i
		$
		\State Compute $u^k_i$ via gossip communication with neighbors: $\u^k =  \lambda W \w^k$
		\State Solve local subproblem, i.e. find $y^{k+1}_i \in \R^{nd}$, such that $\|\hat y^{k+1}_i -  y^{k+1}_i\|_2^2 \leq \delta$, where
		\begin{equation*}
		    y^{k+1}_i = \argmin_{y_i\in \R^{d}} \left\{ \langle u^k_i, y_i-w^k_i\rangle  + f_i(y_i)   +\frac{\gamma}{2}\|y_i-w^k_i\|_2^{2} \right\}
		\end{equation*}
		\State Compute $z^k_i$ via gossip communication with neighbors: $\z^k =  \lambda W \y^{k+1}$
		\State Local update: $x^{k+1}_i = x^k_i-a^{k+1} (\nabla f_i(y^{k+1}_i) + z^k_i)$
		\EndFor
	\end{algorithmic}
	\hspace*{0.2cm} {\bf Output:} $\{y^{K}_i\}$
\end{algorithm}

\begin{algorithm}
	\caption{\texttt{ARDM}}
	\label{alg_pfl_sum_prac}
	\hspace*{0.2cm} {\bf Input:} starting point $x^0_i = x^0 \in \R^{d}$, number of iterations $K$,  parameters $0<\theta_1, \theta_2 <1$, $\eta, \beta, \gamma > 0$, probabilities $p$, $\rho$\\
    \hspace*{0.3cm}{\bf Initialization:} $ y^0_i = z^0_i = u^0_i = x^0_i$ and $\hat g^{0}_i=\lambda y^{0}_i + \nabla f_i(y^{0}_i)$
	\begin{algorithmic}[1]
\For {$k=0,1, 2, \ldots K-1$ }
\State Local update: $x^k_i = \theta_1 z^k_i + \theta_2 u^k_i + (1 - \theta_1 -\theta_2)y^k_i$
\State Generate $\xi^k =  \begin{cases}
 1,&  \text{with probability} ~~ 1 - p \\
0 ,& \text{with probability} ~~ p
\end{cases}$
\If{$\xi^k = 0$} 
\State Compute $b^k_i$ via gossip communication with neighbors:
$\b^{k} = \frac{\lambda}{p} W\x^k$
\State Local update: $g^{k}_i = b^k_i - a^k_i + \hat g^{k}_i$
\Else 
\State Sample indices $m^k_1, \ldots m^k_n$ for each node independently and uniformly from $[M]$
\State Local computation: $g^{k}_i = \frac{1}{1-p}\left(\nabla f_{i,m^k_i}(x^k_i) - \nabla f_{i,m^k_i}(u^k_i)\right) + \hat g^{k}_i$ 
\EndIf
\State Local update: $y^{k+1}_i = x^k_i - \eta g^k_i$
\State Local update: $z^{k+1}_i = \beta z^k_i + (1-\beta)x^k_i + \frac{\gamma}{\eta}(y^{k+1}_i - x^k_i)$
\State Generate $\xi^{k+\frac{1}{2}}=  \begin{cases}
1,&  \text{with prob.} ~~ 1 - \rho \\
0 ,& \text{with prob.} ~~ \rho
\end{cases}$
\If{$\xi^{k+1/2} = 0$}
\State $u^{k+1}_i = y^{k+1}_i$
\State Compute $a^{k+1}_i$ via gossip communication with neighbors: $\a^{k+1} = \frac{\lambda}{p} W\u^{k+1}$
\State Compute $c^k_i$ via gossip communication with neighbors: $\c^k = \lambda W\y^{k+1}$
\State Local update: $\hat g^{k+1}_i = c^k_i + \nabla f_i(y^{k+1}_i)$
\Else
\State $u^{k+1}_i = u^{k}_i$
\State $a^{k+1}_i = a^k_i$
\State $\hat g^{k+1}_i = \hat g^{k}_i$
\EndIf
\EndFor
\end{algorithmic}
\end{algorithm}

\newpage

\section{Proof of Theorem \ref{theorem_lower_bounds}}\label{ap:T1}
In this section, we prove lower convergence bounds of algorithms satisfying Assumption \ref{as:lower} for the problem  \eqref{rewrite_main_problem}. To do this, we need to give an example of `bad' functions that satisfy Assumption \ref{as:func}, and an example of a `bad' arrangement of these functions in some graph with a `bad' matrix $\hat W$ (Definition \ref{gossipmatrix}) with an upper bound of condition number $\chi$. Following \cite{nesterov2003introductory,hanzely2020lower} we consider quadratic functions, and following \cite{scaman2017optimal}, we construct a linear graph. 

Let us start with the network. As the gossip matrix, we take the Laplacian of the linear graph. Then, for our problem \eqref{rewrite_main_problem}, we get that the matrix $W$ has the following form $W = \hat{W}\otimes I_d$, where $\hat{W} = \tfrac1{2} U$, and $U$ is 
\begin{eqnarray}\label{U}
    U = 
    \begin{pmatrix}
        1 & -1 & & &   \\
        -1 & 2 & -1 &  &  \\
        & \ddots&\ddots & \ddots&  \\
        & &   -1 & 2 & -1\\
        & &    & -1 & 1\\
    \end{pmatrix}.
\end{eqnarray}
It is known that the spectrum of the (scaled by $1/2$) Laplacian matrix of the linear graph with $n$ vertices, is $2 \sin^2\left( \tfrac{\pi k}{2n}\right)$ for $k = 0, \ldots, n-1,$ \cite{anderson1985eigenvalues}. Thus, the condition number is $\chi (n) = \tfrac{\sin^2\left( \tfrac{\pi (n-1)}{2n}\right)}{\sin^2\left( \tfrac{\pi}{2n}\right)}$. Since we consider $\chi \geq 3$, one can find $n \geq 3$ such that $\chi (n) \leq \chi <  \chi (n+1)$. Moreover, for $n \geq 3$ we can guarantee that $\lambda_{\max}(n) \geq \tfrac{3}{2}$, $ \tfrac{4}{n^2} \leq \lambda^+_{\min}(n)  \leq \tfrac{5}{n^2}$ and $\chi({n+1}) \leq \tfrac{1}{\sin^2\left( \tfrac{\pi}{2(n+1)}\right)} \leq \tfrac{(n+1)^2}{2}$.
It turns out that if we choose as the `bad' network, a linear graph with $n$ vertices (where $n$ is such that $\chi_{n} \leq \chi <  \chi_{n+1}$), and take the Laplacian of this graph as the gossip matrix, then we satisfy Definition \ref{gossipmatrix} and $\chi$ is an upper bound for the condition number of the gossip matrix. And one can note that $n - 1 > \sqrt{2 \chi} - 2 \geq \tfrac{1}{5} \sqrt{\chi}$ (since $\chi \geq 3$), $1 \leq \tfrac{2}{3} \lambda_{\max}(n)$ and $ \frac{4}{n^2}\leq \lambda^+_{\min}(n) \leq \tfrac{5}{n^2}$.

Now let us move on to the `bad' functions. We choose the dimension of these functions equivalent to $d=2T$ with large enough $T$ (to be defined later). Next, we divide the nodes of the network into three types: the first type includes $\mathcal{V}_1=\left\{1\right\}$, the second type includes $\mathcal{V}_2=\left\{2, n - 1\right\}$,  the third type includes $\mathcal{V}_3=\left\{n\right\}$. Each type of node has its own functions: 
\begin{equation}
    \label{eq:bad_func}
    f_i(x) = 
    \begin{cases}
        \frac{\mu}{2}\|x\|^2 + a x^{(1)} + \frac{c\lambda}{2}\left(\sum^{T-1}_{t = 1}\left(x^{(2t)}- x^{(2t+1)}\right)^2\right) + \frac{b\lambda }{2}\left({x}^{(2T)}\right)^2, & \text{if } i \in \mathcal{V}_1,\\
        \phi \cdot \frac{\mu}{2}\|x\|^2, & \text{if } i \in \mathcal{V}_2,\\ 
        \frac{\mu}{2}\|x\|^2  + \frac{c\lambda}{2}\left(\sum^{T-1}_{t = 0}\left(x^{(2t+1)}- x^{(2t+2)}\right)^2\right), & \text{if } i \in \mathcal{V}_3,
    \end{cases}
\end{equation}
where constants $a, b, c$ will be defined shortly. The parameter $\phi$ takes two values: $1$ or $0$. We will consider both values below, we need $0$ to simplify the mathematical calculations, note that in this case we slightly change the class of problems, since not all functions $f_i$ are strongly convex and we slightly go beyond Assumption~\ref{as:func}.

In the proof we will rely on \cite{hanzely2020lower}. In particular, we will prove similar (but not analogous) lemmas.

Let us introduce the solution of the problem \eqref{rewrite_main_problem} with \eqref{eq:bad_func}. For the first type of node, we denote the solution by $x^*$, for the third type node by $z^*$, and for the second type nodes by $y^*_2$, \ldots, $y^*_{n-1}$.
Using this notation we write down the optimality conditions for \eqref{rewrite_main_problem}. First write down for $x^*$:
\begin{equation}
    \label{x_y_3}
    \left( \frac{\mu}{\lambda} + \frac{1}{2} \right) (x^*)^{(1)} + \frac{a}{\lambda} - \frac{1}{2} (y^*_2)^{(1)} = 0,
\end{equation}
\begin{equation}
\label{x_y_2}
    \left(c + \frac{\mu}{\lambda} + \frac{1}{2}\right){(x^*)}^{(2t)}- c {(x^*)}^{(2t+1)} - \frac{1}{2}{(y^*_2)}^{(2t)} = 0,\quad \text{for } 1\leq t \leq T-1,
\end{equation}
\begin{equation}
\label{x_y_1}
    \left(c + \frac{\mu}{\lambda} + \frac{1}{2}\right){(x^*)}^{(2t+1)}- c {(x^*)}^{(2t)} - \frac{1}{2}{(y^*_2)}^{(2t+1)} = 0,\quad \text{for } 1\leq t \leq T-1,
\end{equation}
\begin{equation}
    \label{x_y_4}
    \left( \frac{\mu}{\lambda} + b + \frac{1}{2} \right) (x^*)^{(2T)} - \frac{1}{2} (y^*_2)^{(2T)} = 0.
\end{equation}
Then for $z^*$:
\begin{equation}
\label{z_y_1}
    \left(c + \frac{\mu}{\lambda} + \frac{1}{2}\right){(z^*)}^{(2t-1)}- c {(z^*)}^{(2t)} - \frac{1}{2}{(y^*_{n-1})}^{(2t-1)} = 0,\quad \text{for } 1\leq t \leq T,
\end{equation}
\begin{equation}
\label{z_y_2}
    \left(c + \frac{\mu}{\lambda} + \frac{1}{2}\right){(z^*)}^{(2t)}- c {(z^*)}^{(2t-1)} - \frac{1}{2}{(y^*_{n-1})}^{(2t)} = 0,\quad \text{for } 1\leq t \leq T,
\end{equation}
Finally for $y^*_2$, \ldots, $y^*_{n-1}$:
\begin{equation}
    \label{yy_x}
    \left(1 + \frac{\phi\mu}{\lambda}\right) (y^*_2)^{(t)} - \frac{1}{2}(y^*_3)^{(t)} - \frac{1}{2}(x^*)^{(t)} = 0, \quad \text{for } 1\leq t \leq 2T,
\end{equation}
\begin{equation}
\label{yy}
    \left(1 + \frac{\phi\mu}{\lambda}\right) (y^*_{i})^{(t)} - \frac{1}{2}(y^*_{i+1})^{(t)} - \frac{1}{2}(y^*_{i-1})^{(t)} = 0, \quad \text{for } 1\leq t \leq 2T,
\end{equation}
\begin{equation}
\label{yy_z}
    \left(1 + \frac{\phi\mu}{\lambda}\right) (y^*_{n-1})^{(t)} - \frac{1}{2}(y^*_{n-2})^{(t)} - \frac{1}{2}(z^*)^{(t)} = 0,\quad \text{for } 1\leq t \leq 2T.
\end{equation}

First, we give a proof of the lemma that indicates a recursive connection of coordinates $x^*$ and $z^*$. Before we introduce new notation:
\begin{equation*}
		w_t =
		\begin{cases}
		\begin{pmatrix}
		{(z^*)}^{(t)}\\
		{(x^*)}^{(t)}
		\end{pmatrix} & \text{if } t \text{ is even} \\
		\begin{pmatrix}
		{(x^*)}^{(t)}\\
		{(z^*)}^{(t)}
		\end{pmatrix} & \text{if } t \text{ is odd}
		\end{cases}.
	\end{equation*}
\begin{lemma}\label{lem:1404}
The sequence $w_t$ satisfies the following recursion relation:
	\begin{equation*}
		w_{t+1}
		=
		Q
		w_{t} \quad \text{with} \quad
		Q = \begin{pmatrix}
	- \frac{B}{2c} & ~~~~   &\frac{1}{c} \left(c+ \frac{\mu}{\lambda} +\frac{1}{2} - \frac{A}{2}\right)   \\
	-\frac{1}{c}\left(c+ \frac{\mu}{\lambda} +\frac{1}{2} - \frac{A}{2}\right) &~~~~&
	\frac{2}{Bc}\left(c+ \frac{\mu}{\lambda} +\frac{1}{2} - \frac{A}{2}\right)^2 - \frac{2c}{B}   \\
	\end{pmatrix},
	\end{equation*}
	where
	\begin{equation*}
		A = \left( 1 - \frac{1}{n-1} \right), \quad B = \frac{1}{n-1}, \quad \text{for} \quad \phi = 0,
	\end{equation*}
	or
	\begin{equation*}
		A = \frac{\omega_2^{n-2} -  \omega_1^{n-2}}{\omega_2^{n-1} - \omega_1^{n-1}}, \quad B = \frac{\omega_2 -  \omega_1}{\omega_2^{n-1} - \omega_1^{n-1}}, \quad \text{for} \quad \phi = 1,
	\end{equation*}
	with $\omega_1 = 1 + \tfrac{\mu}{\lambda} - \sqrt{\tfrac{2\mu}{\lambda} + \tfrac{\mu^2}{\lambda^2}}$ and $\omega_2 = 1 + \tfrac{\mu}{\lambda} + \sqrt{\tfrac{2\mu}{\lambda} + \tfrac{\mu^2}{\lambda^2}}$.
\end{lemma}

\begin{proof} We start from \eqref{yy_x}, \eqref{yy}, \eqref{yy_z}. One can note that we have recursion with two initial conditions:
\begin{equation*}
    (y^*_{i})^{(t)} = \left(2 + \frac{2\phi\mu}{\lambda}\right) (y^*_{i-1})^{(t)} - (y^*_{i-2})^{(t)} \quad \text{with} \quad (y^*_{1})^{(t)} = (x^*)^{(t)}, ~~(y^*_{n})^{(t)} = (z^*)^{(t)}.
\end{equation*}
If $\phi = 0$, the expressions for $(y^*_{i})^{(t)}$ are as follows:
\begin{equation*}
    (y^*_{i})^{(t)} = \left(\frac{i}{n-1} - \frac{1}{n-1} \right) (z^*)^{(t)}
    + \left( \frac{n}{n-1}  - \frac{i}{n-1}\right) (x^*)^{(t)}.
\end{equation*}
In particular, $(y_2^*)^{(t)} = \left( 1 - \tfrac{1}{n-1} \right) (x^*)^{(t)} + \tfrac{1}{n-1} (z^*)^{(t)}$ and $(y_{n-1}^*)^{(t)} = \left( 1 - \tfrac{1}{n-1} \right) (z^*)^{(t)} + \tfrac{1}{n-1} (x^*)^{(t)}$. When $\phi = 1$, the expressions for $(y^*_{i})^{(t)}$ become more complicated:
\begin{equation*}
    (y^*_{i})^{(t)} = C_1 \omega_1^{i-1} + C_2 \omega_2^{i-1} = \frac{\omega_2^{n-i} - \omega_1^{n-i}}{\omega_2^{n-1} - \omega_1^{n-1}} (x^*)^{(t)} + \frac{\omega_2^{i-1} -  \omega_1^{i-1}}{\omega_2^{n-1} - \omega_1^{n-1}} (z^*)^{(t)},
\end{equation*}
with $\omega_1 = 1 + \tfrac{\mu}{\lambda} - \sqrt{\tfrac{2\mu}{\lambda} + \tfrac{\mu^2}{\lambda^2}}$ and $\omega_2 = 1 + \tfrac{\mu}{\lambda} + \sqrt{\tfrac{2\mu}{\lambda} + \tfrac{\mu^2}{\lambda^2}}$.
In particular, $(y_2^*)^{(t)} = \tfrac{\omega_2^{n-2} -  \omega_1^{n-2}}{\omega_2^{n-1} - \omega_1^{n-1}} (x^*)^{(t)} + \tfrac{\omega_2 -  \omega_1}{\omega_2^{n-1} - \omega_1^{n-1}} (z^*)^{(t)}$ and $(y_{n-1}^*)^{(t)} = \frac{\omega_2^{n-2} -  \omega_1^{n-2}}{\omega_2^{n-1} - \omega_1^{n-1}} (z^*)^{(t)} + \tfrac{\omega_2 -  \omega_1}{\omega_2^{n-1} - \omega_1^{n-1}} (x^*)^{(t)}$. In both cases of $\phi$ we have that $(y_2^*)^{(t)} = A 
\cdot (x^*)^{(t)} + B \cdot (z^*)^{(t)}$ and $(y_{n-1}^*)^{(t)} = A \cdot (z^*)^{(t)} + B \cdot (x^*)^{(t)}$ with some $A$ and $B$. We can substitute these $(y_2^*)^{(t)}$ and $(y_{n-1}^*)^{(t)}$ into \eqref{x_y_2}, \eqref{x_y_1}, \eqref{z_y_1}, \eqref{z_y_2} and have:
\begin{equation*}
    \left(c + \frac{\mu}{\lambda} + \frac{1}{2}\right){(x^*)}^{(2t)}- c {(x^*)}^{(2t+1)} - \frac{A}{2} (x^*)^{(2t)} - \frac{B}{2} (z^*)^{(2t)} = 0,\quad \text{for } 1\leq t \leq T-1,
\end{equation*}
\begin{equation*}
    \left(c + \frac{\mu}{\lambda} + \frac{1}{2}\right){(x^*)}^{(2t+1)}- c {(x^*)}^{(2t)} - \frac{A}{2} (x^*)^{(2t+1)} - \frac{B}{2} (z^*)^{(2t+1)} = 0,\quad \text{for } 1\leq t \leq T-1,
\end{equation*}
\begin{equation}
    \label{zx_3}
    \left(c + \frac{\mu}{\lambda} + \frac{1}{2}\right){(z^*)}^{(2t-1)}- c {(z^*)}^{(2t)} - \frac{A}{2}(z^*)^{(2t-1)} - \frac{B}{2} \cdot (x^*)^{(2t-1)} = 0,\quad \text{for } 1\leq t \leq T,
\end{equation}
\begin{equation}
    \label{zx_4}
    \left(c + \frac{\mu}{\lambda} + \frac{1}{2}\right){(z^*)}^{(2t)}- c {(z^*)}^{(2t-1)} - \frac{A}{2}(z^*)^{(2t)} - \frac{B}{2} \cdot (x^*)^{(2t)} = 0,\quad \text{for } 1\leq t \leq T,
\end{equation}
The first two expressions together can be rewritten as follows:
\begin{equation*}
	\begin{pmatrix}
	c &~~~~ & 0 \\
	- c- \frac{\mu}{\lambda}-\frac{1}{2} + \frac{A}{2} &~~~~ & \frac{B}{2}  \\
	\end{pmatrix}
	\begin{pmatrix}
	{(x^*)}^{(2t+1)}\\
	{(z^*)}^{(2t+1)}
	\end{pmatrix}
	=
	\begin{pmatrix}
	c+ \frac{\mu}{\lambda} +\frac{1}{2} - \frac{A}{2} &~~~~   & - \frac{B}{2} \\
	-c &~~~~& 0 \\
	\end{pmatrix}
	\begin{pmatrix}
	{(x^*)}^{(2t)}\\
	{(z^*)}^{(2t)}
	\end{pmatrix},
\end{equation*}
or
\begin{align*}
	\begin{pmatrix}
	{(x^*)}^{(2t+1)}\\
	{(z^*)}^{(2t+1)}
	\end{pmatrix}
	&= \begin{pmatrix}
	c &~~~~ & 0 \\
	- c- \frac{\mu}{\lambda}-\frac{1}{2} + \frac{A}{2} &~~~~ & \frac{B}{2}  \\
	\end{pmatrix}^{-1}
	\begin{pmatrix}
	c+ \frac{\mu}{\lambda} +\frac{1}{2} - \frac{A}{2} &~~~~   & - \frac{B}{2} \\
	-c &~~~~& 0 \\
	\end{pmatrix}
	\begin{pmatrix}
	{(x^*)}^{(2t)}\\
	{(z^*)}^{(2t)}
	\end{pmatrix}\\
	&= \frac{2}{Bc}
	\begin{pmatrix}
	\frac{B}{2} &~~~~ & 0 \\
	 c + \frac{\mu}{\lambda} + \frac{1}{2} - \frac{A}{2} &~~~~ &  c  \\
	\end{pmatrix}
	\begin{pmatrix}
	c+ \frac{\mu}{\lambda} +\frac{1}{2} - \frac{A}{2} &~~~~   & - \frac{B}{2} \\
	-c &~~~~& 0 \\
	\end{pmatrix}
	\begin{pmatrix}
	{(x^*)}^{(2t)}\\
	{(z^*)}^{(2t)}
	\end{pmatrix} \\
	&=
	\begin{pmatrix}
	\frac{1}{c} \left(c+ \frac{\mu}{\lambda} +\frac{1}{2} - \frac{A}{2}\right) &~~~~   & - \frac{B}{2c} \\
	\frac{2}{Bc}\left(c+ \frac{\mu}{\lambda} +\frac{1}{2} - \frac{A}{2}\right)^2 - \frac{2c}{B} &~~~~& -\frac{1}{c}\left(c+ \frac{\mu}{\lambda} +\frac{1}{2} - \frac{A}{2}\right) \\
	\end{pmatrix}
	\begin{pmatrix}
	{(x^*)}^{(2t)}\\
	{(z^*)}^{(2t)}
	\end{pmatrix} \\
	&=
	\begin{pmatrix}
	- \frac{B}{2c} & ~~~~   &\frac{1}{c} \left(c+ \frac{\mu}{\lambda} +\frac{1}{2} - \frac{A}{2}\right)   \\
	-\frac{1}{c}\left(c+ \frac{\mu}{\lambda} +\frac{1}{2} - \frac{A}{2}\right) &~~~~&
	\frac{2}{Bc}\left(c+ \frac{\mu}{\lambda} +\frac{1}{2} - \frac{A}{2}\right)^2 - \frac{2c}{B}   \\
	\end{pmatrix}
	\begin{pmatrix}
	{(z^*)}^{(2t)}\\
	{(x^*)}^{(2t)}
	\end{pmatrix} \\
	&=
	Q \begin{pmatrix}
	{(z^*)}^{(2t)}\\
	{(x^*)}^{(2t)}
	\end{pmatrix}.
\end{align*}
Similarly, from \eqref{zx_3} and \eqref{zx_4} one can get that \begin{align*}
	\begin{pmatrix}
	{(z^*)}^{(2t)}\\
	{(x^*)}^{(2t)}
	\end{pmatrix}
	&=
	Q \begin{pmatrix}
	{(x^*)}^{(2t-1)}\\
	{(z^*)}^{(2t-1)}
	\end{pmatrix}.
\end{align*}
Using the definition of $w_t$ completes the proof. 
\EndProof
\end{proof}

Then, we follow the idea from \cite{hanzely2020lower}. 
From the proof of the previous lemma we know that $(y_2^*)^{(t)} = A 
\cdot (x^*)^{(t)} + B \cdot (z^*)^{(t)}$. Then, substituting  $(y_2^*)^{(1)}$ and $(y_2^*)^{(2T)}$ into \eqref{x_y_3} and \eqref{x_y_4}, we obtain that the value of $w_1$ and $w_{2T}$ depends on the parameters $a$ and $b$.  Hence, by varying the parameters $a$ and $b$, one can obtain that $w_1$, $w_2$, \ldots, $w_{2T}$ are eigenvectors of the matrix $Q$, i.e. $w_2 = Qw_1 = \gamma w_1$ etc. This idea is implemented in the following lemma.

\begin{lemma}\label{lem:2404}
For any $L$, $\mu$, and $\lambda$ ($L \geq 2\mu$, and $\lambda\eigWmin \geq \mu$), there exists a choice of parameters $a$, $b$, $c$ such that $w_1$, $w_2$, \ldots $w_{2T}$ are eigenvectors of matrix $Q$ corresponding to the eigenvalue $\gamma \in (0;1)$, where
\begin{equation*}
    \gamma \geq 1 - \max \left\{ 2\sqrt{\frac{\mu n^2 }{\lambda}}, 3\sqrt{\frac{\mu}{L - \mu}}\right\}.
\end{equation*}
Moreover, the problem \eqref{rewrite_main_problem} $+$ \eqref{eq:bad_func} with these parameters $a$, $b$, $c$ satisfies Assumption \ref{as:func}. 
\end{lemma}
\begin{proof}
First we give the values of $a$, $b$, and $c$:
\begin{equation*}
c = 
\begin{cases}
1, \quad \text{for} \quad \mu + \lambda \leq L,\\
\frac{\mu}{\lambda} \cdot \delta = \frac{\mu}{\lambda} \cdot \frac{L-\mu}{\mu}, \quad \text{for} \quad \mu + \lambda > L,
\end{cases}
\end{equation*}
% \begin{equation}
% \label{eq:a}
% a = \frac{1}{2}\left(A \lambda \alpha + B \lambda - (2 \mu + \lambda) \alpha\right),
% \end{equation}
\begin{equation}
\label{eq:b}
b = \frac{B  \alpha }{2} - \frac{\mu}{\lambda} - \frac{1}{2} + \frac{A}{2} \quad \text{and} \quad \text{any } a,
\end{equation}
where
\begin{equation*}
\resizebox{\linewidth}{!}{
$
\alpha = -\frac{1 - 2A + A^2 + B^2 + 4c - 4Ac + 4\frac{\mu}{\lambda} - 4A\frac{\mu}{\lambda} + 8c\frac{\mu}{\lambda} + 4\frac{\mu^2}{\lambda^2} + \sqrt{\left(-1 + 2A - A^2 + B^2 - 4\frac{\mu}{\lambda} + 4A\frac{\mu}{\lambda} - 4\frac{\mu^2}{\lambda^2}\right)\left(-1 + 2A - A^2 + B^2 - 8c + 8Ac - 16c^2 - 4\frac{\mu}{\lambda} + 4A\frac{\mu}{\lambda} - 16c\frac{\mu}{\lambda} - 4\frac{\mu^2}{\lambda^2}\right)}}{2B\left(-1 + A - 2c - 2\frac{\mu}{\lambda}\right)}.
$
}
\end{equation*}
Let us check that the problem \eqref{rewrite_main_problem} + \eqref{eq:bad_func} satisfies Assumption \ref{as:func}. Note that by the choice of $c$, it suffices to verify that $ 0 \leq b\lambda \leq c\lambda \leq L - \mu$. We make this verification with Mathematica (here and below, when using Mathematica, we replace $\tfrac{\mu}{\lambda}$ with $x$). First, we check these inequalities when $\varphi = 0$ ($A = \tfrac{n-2}{n-1}$ and $B = \tfrac{1}{n-1}$):
\begin{itemize}[label=$\ast$]
\item $b \lambda \leq c \lambda$ (or $\alpha \leq \tfrac{1}{B}\left(2c + 1 - A + 2\tfrac{\mu}{\lambda}\right)$) for $0 < c \leq 1$ and $x = \tfrac{\mu}{\lambda} > 0$, $x = \tfrac{\mu}{\lambda} \leq  \lambda^+_{\min} \leq \tfrac{5}{n^2}$ (since in Theorem \ref{theorem_lower_bounds} we assume that $\tfrac{\mu}{\lambda} \leq  \lambda^+_{\min}$ and above we estimated that $\lambda^+_{\min} \leq \tfrac{5}{n^2}$)
\begin{figure}[H]
    \centering
    \includegraphics[width = 1\linewidth]{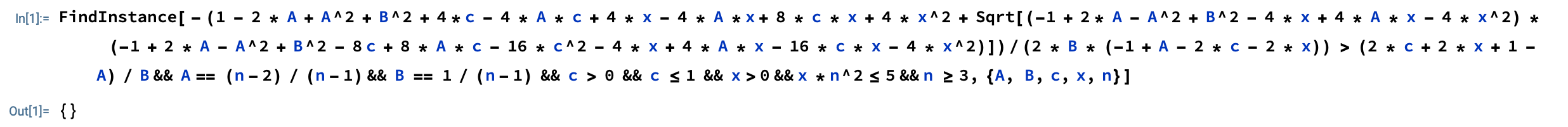}
\end{figure}
\item $ 0 \leq b$ (or $\alpha \geq \tfrac{1}{B}\left(1 - A + 2\tfrac{\mu}{\lambda}\right)$) for $0 < c \leq 1$ and $x = \tfrac{\mu}{\lambda} > 0$, $x = \tfrac{\mu}{\lambda} \leq  \lambda^+_{\min} \leq \tfrac{5}{n^2}$
\begin{figure}[H]
    \centering
    \includegraphics[width = 1\linewidth]{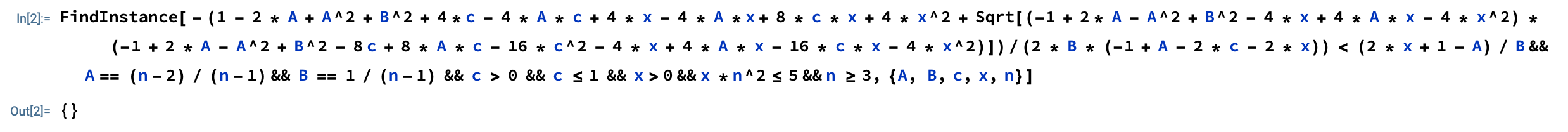}
\end{figure}
\end{itemize}
In the case of $\varphi = 1$, we replace the expressions for $A$ and $B$ from Lemma \ref{lem:1404} by their Taylor approximations:
\begin{figure}[H]
    \centering
    \includegraphics[width = 1\linewidth]{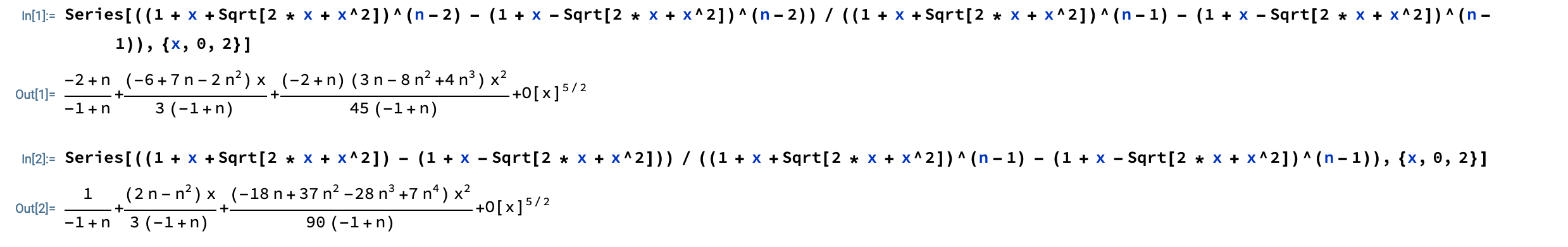}
\end{figure}
\begin{equation}
    \label{AB}
    \begin{split}
        A \approx \frac{n-2}{n-1} - \frac{2n^2 - 7n +6}{3(n-1)}\frac{\mu}{\lambda} + \frac{(4n^3 - 8n^2 + 3n)(n-2)}{45(n-1)}\frac{\mu^2}{\lambda^2}, \\
        B \approx \frac{n-2}{n-1} - \frac{n^2 - 2n}{3(n-1)}\frac{\mu}{\lambda} + \frac{7n^4 - 28n^3+37n^2-18n}{90(n-1)}\frac{\mu^2}{\lambda^2}.
    \end{split}
\end{equation}
Then, we can check inequalities for $b$:
\begin{itemize}[label=$\ast$]
\item $b \lambda \leq c \lambda$ (or $\alpha \leq \tfrac{1}{B}\left(2c + 1 - A + 2\tfrac{\mu}{\lambda}\right)$) for $0 < c \leq 1$ and $x = \tfrac{\mu}{\lambda} > 0$, $x = \tfrac{\mu}{\lambda} \leq  \lambda^+_{\min} \leq \tfrac{5}{n^2}$
\begin{figure}[H]
    \centering
    \includegraphics[width = 1\linewidth]{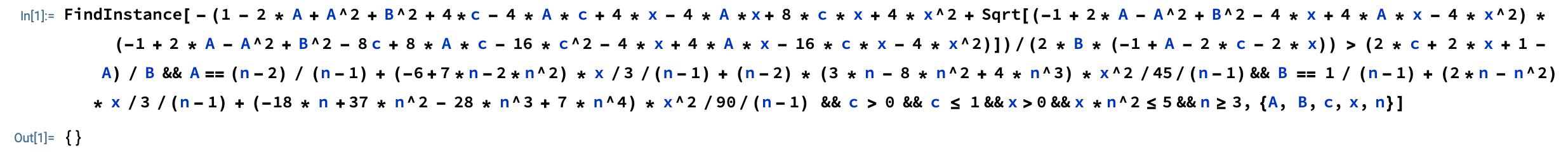}
\end{figure}
\item $ 0 \leq b$ (or $\alpha \geq \tfrac{1}{B}\left(1 - A + 2\tfrac{\mu}{\lambda}\right)$) for $0 < c \leq 1$ and $x = \tfrac{\mu}{\lambda} > 0$, $x = \tfrac{\mu}{\lambda} \leq  \lambda^+_{\min} \leq \tfrac{5}{n^2}$
\begin{figure}[H]
    \centering
    \includegraphics[width = 1\linewidth]{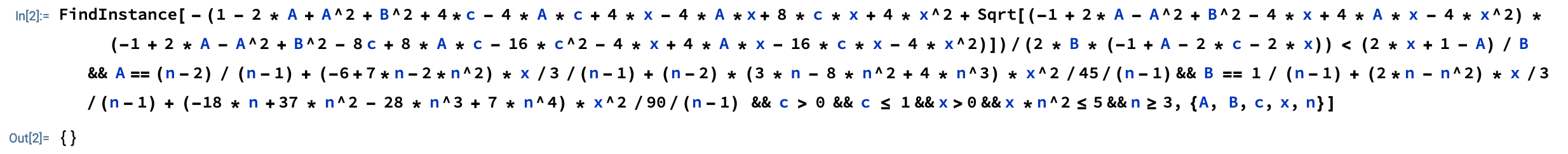}
\end{figure}
\end{itemize}
Next, we turn to eigenvalues and vectors. One can find them:
\begin{figure}[h!]
    \centering
    \includegraphics[width = 1\linewidth]{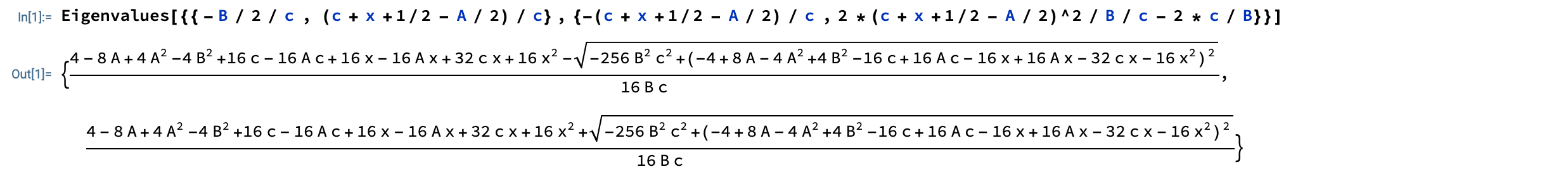}
\end{figure}
\begin{figure}[h!]
    \centering
    \includegraphics[width = 1\linewidth]{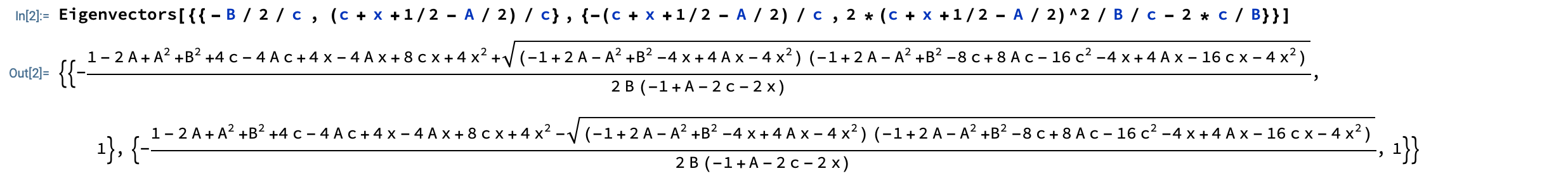}
\end{figure}

We take the smallest eigenvalue
\begin{equation*}
\resizebox{\linewidth}{!}{
$
\gamma = \frac{4 - 8A + 4A^2 - 4B^2 + 16 c - 16Ac + 16 \frac{\mu}{\lambda} - 16 A\frac{\mu}{\lambda} + 32 c\frac{\mu}{\lambda} + 16 \frac{\mu^2}{\lambda^2} - \sqrt{-256 B^2 c^2 + \left(-4 + 8A - 4A^2 + 4 B^2 - 16 c + 16Ac - 16\frac{\mu}{\lambda} + 16 A\frac{\mu}{\lambda} - 32c\frac{\mu}{\lambda} - 16\frac{\mu^2}{\lambda^2}\right)^2}}{16 B c}
$
}
\end{equation*}

and the corresponding eigenvector
\begin{equation*}
v = 
\begin{pmatrix}
\alpha \\
1
\end{pmatrix}.
\end{equation*}
By simply substituting $b$ from expression \eqref{eq:b} and $(y_2^*)^{(t)} = A 
\cdot (x^*)^{(t)} + B \cdot (z^*)^{(t)}$ into equations \eqref{x_y_4}, one can note that $w_{2T}$ is an eigenvector of $Q$. It means that $\gamma w_{2T} = Q w_{2T}$ or $w_{2T} =  \gamma Q^{-1} w_{2T}$. From Lemma \ref{lem:1404} we also have $Q^{-1} w_{2T} = w_{2T-1}$. As the result, $w_{2T} = \gamma w_{2T-1}$, i.e. $w_{2T-1}$ is also an eigenvector of $Q$. Continuing further, we can obtain that all vectors $w_{2T}$, \ldots, $w_1$ are eigenvectors of $Q$. The choice of parameter $a$ does not affect, it only determines the value of $\|w_1\|$.

Finally, we need to make sure that this $\gamma$ satisfies the conditions of the lemma. Let us consider the three cases separately.

1) $\mu + \lambda \leq L$. In this case $c = 1$. We want to verify that $\gamma \in (0;1)$ and $\gamma \geq 1 - 2 \sqrt{\tfrac{\mu n^2}{\lambda}}$. This inequality need to be checked with the constraints: $x = \tfrac{\mu}{\lambda} > 0$, $x = \tfrac{\mu}{\lambda} \leq  \lambda^+_{\min} \leq \tfrac{5}{n^2}$ (since in Theorem \ref{theorem_lower_bounds} we assume that $\tfrac{\mu}{\lambda} \leq  \lambda^+_{\min}$ and above we estimated that $\lambda^+_{\min} \leq \tfrac{5}{n^2}$, when we construct the network).
First, we check these inequalities when $\varphi = 0$:
\begin{itemize}[label=$\ast$]
\item $\gamma > 0$
\begin{figure}[H]
    \centering
    \includegraphics[width = 1\linewidth]{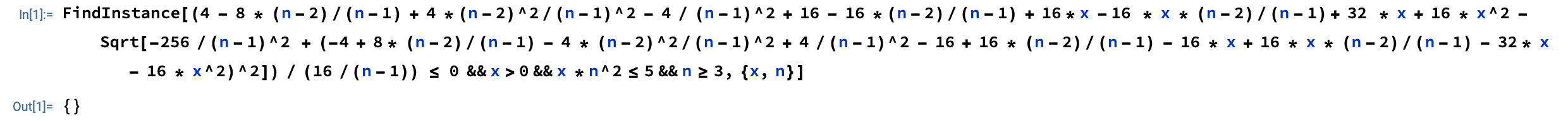}
\end{figure}
\item $\gamma < 1$
\begin{figure}[H]
    \centering
    \includegraphics[width = 1\linewidth]{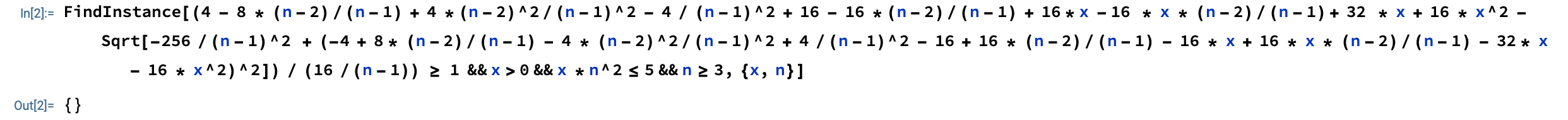}
\end{figure}
\item $\gamma \geq 1 - 2 \sqrt{\tfrac{\mu n^2}{\lambda}}$
\begin{figure}[H]
    \centering
    \includegraphics[width = 1\linewidth]{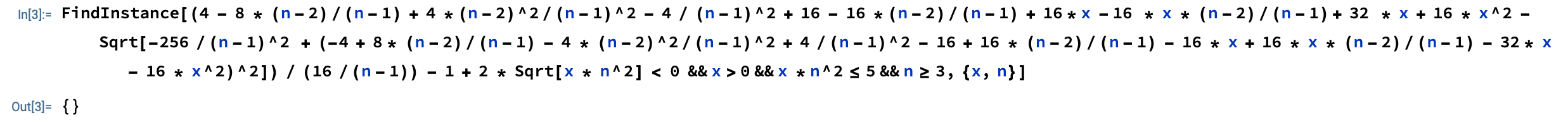}
\end{figure}
\end{itemize}

In the case of $\varphi = 1$, we use \eqref{AB}:
\begin{itemize}[label=$\ast$]
\item $\gamma > 0$
\begin{figure}[H]
    \centering
    \includegraphics[width = 1\linewidth]{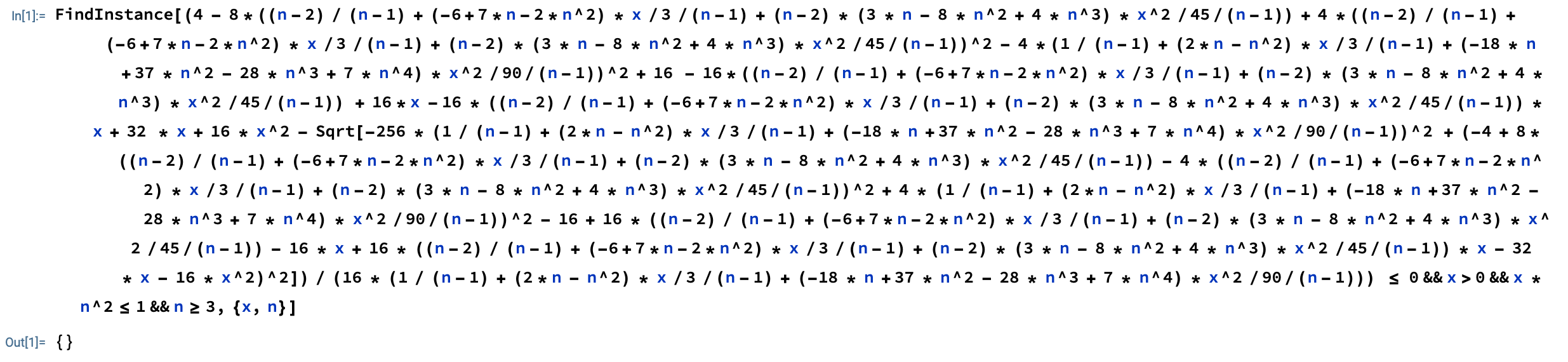}
\end{figure}
\item $\gamma < 1$
\begin{figure}[H]
    \centering
    \includegraphics[width = 1\linewidth]{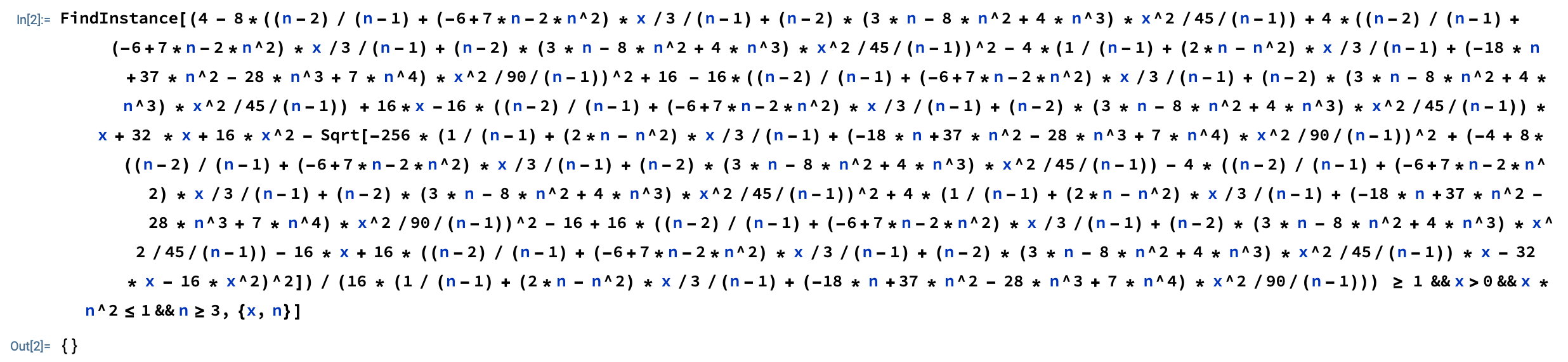}
\end{figure}
\item $\gamma \geq 1 - 2 \sqrt{\tfrac{\mu n^2}{\lambda}}$
\begin{figure}[H]
    \centering
    \includegraphics[width = 1\linewidth]{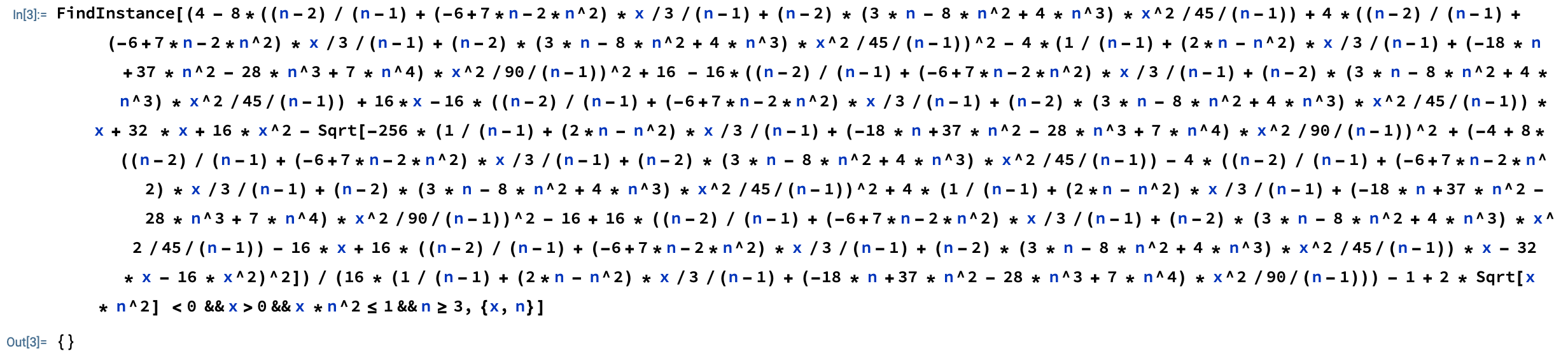}
\end{figure}
\end{itemize}

2) $\mu + \lambda > L \geq \mu + \lambda \lambda_{\min}^+$. In this case $c = \tfrac{L - \mu}{\lambda} = \tfrac{\mu}{\lambda} \cdot \delta = x \delta$. We want to verify that $\gamma \in (0;1)$ and $\gamma \geq 1 - 2 \sqrt{\tfrac{\mu n^2}{\lambda}}$. This inequality need to be checked with the constraints: $\delta \geq 1$ (since in Theorem \ref{theorem_lower_bounds} we assume that $L \geq 2\mu$ and then $\delta = \tfrac{L - \mu}{\mu} \geq 1$), $x = \tfrac{\mu}{\lambda} > 0$, $x = \tfrac{\mu}{\lambda} \leq  \lambda^+_{\min} \leq \tfrac{5}{n^2}$ and $\tfrac{1}{x} > \delta \geq  \tfrac{1}{x} \cdot \lambda^+_{\min} \geq \tfrac{1}{x} \cdot \tfrac{4}{n^2}$ (constraints of the considered case).
First, we check these inequalities when $\varphi = 0$:
\begin{itemize}[label=$\ast$]
\item $\gamma > 0$
\begin{figure}[H]
    \centering
    \includegraphics[width = 1\linewidth]{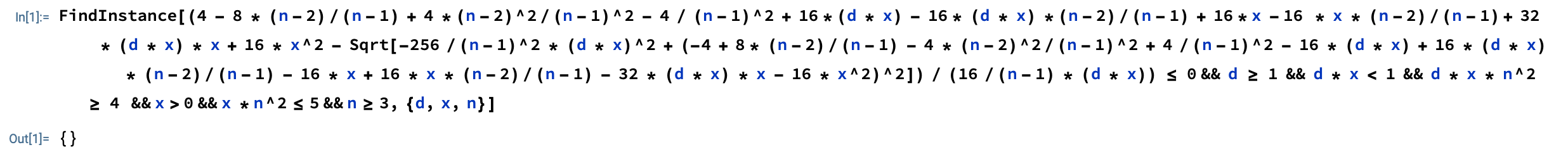}
\end{figure}
\item $\gamma < 1$
\begin{figure}[H]
    \centering
    \includegraphics[width = 1\linewidth]{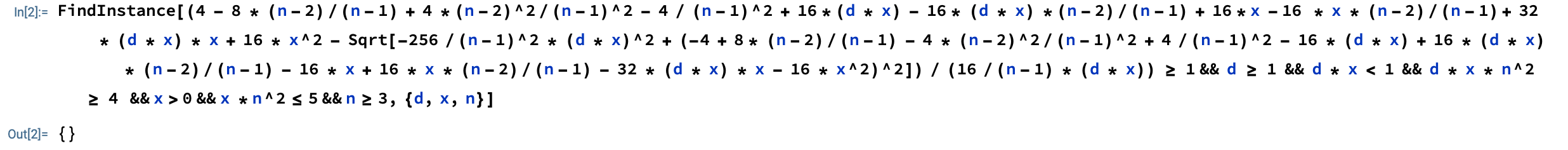}
\end{figure}
\item $\gamma \geq 1 - 2 \sqrt{\tfrac{\mu n^2}{\lambda}}$
\begin{figure}[H]
    \centering
    \includegraphics[width = 1\linewidth]{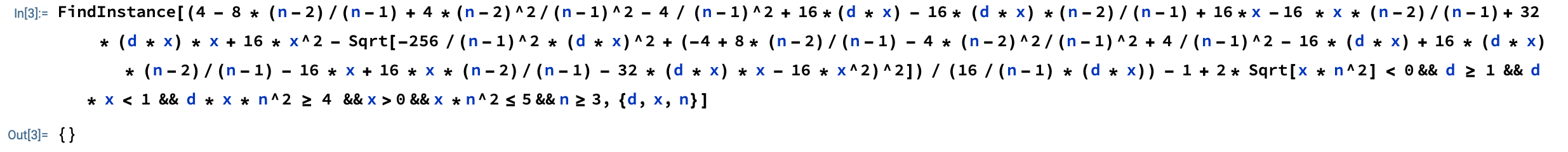}
\end{figure}
\end{itemize}

In the case of $\varphi = 1$, we use \eqref{AB}:
\begin{itemize}[label=$\ast$]
\item $\gamma > 0$
\begin{figure}[H]
    \centering
    \includegraphics[width = 1\linewidth]{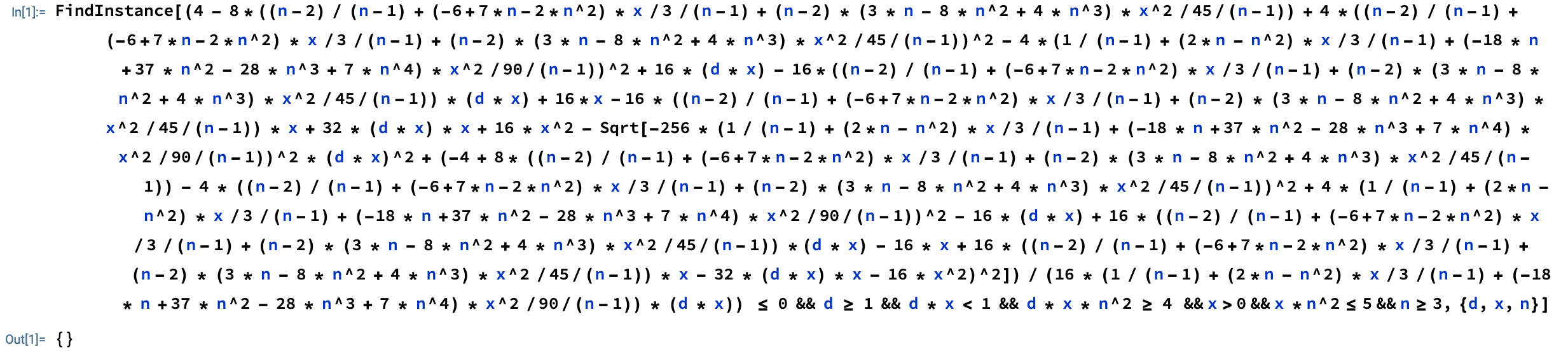}
\end{figure}
\item $\gamma < 1$
\begin{figure}[H]
    \centering
    \includegraphics[width = 1\linewidth]{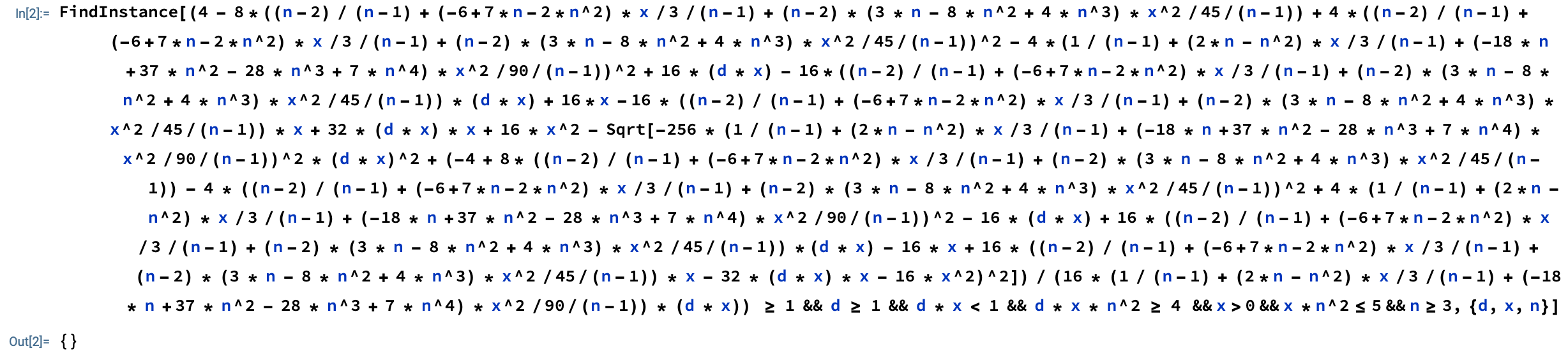}
\end{figure}
\item $\gamma \geq 1 - 2 \sqrt{\tfrac{\mu n^2}{\lambda}}$
\begin{figure}[H]
    \centering
    \includegraphics[width = 1\linewidth]{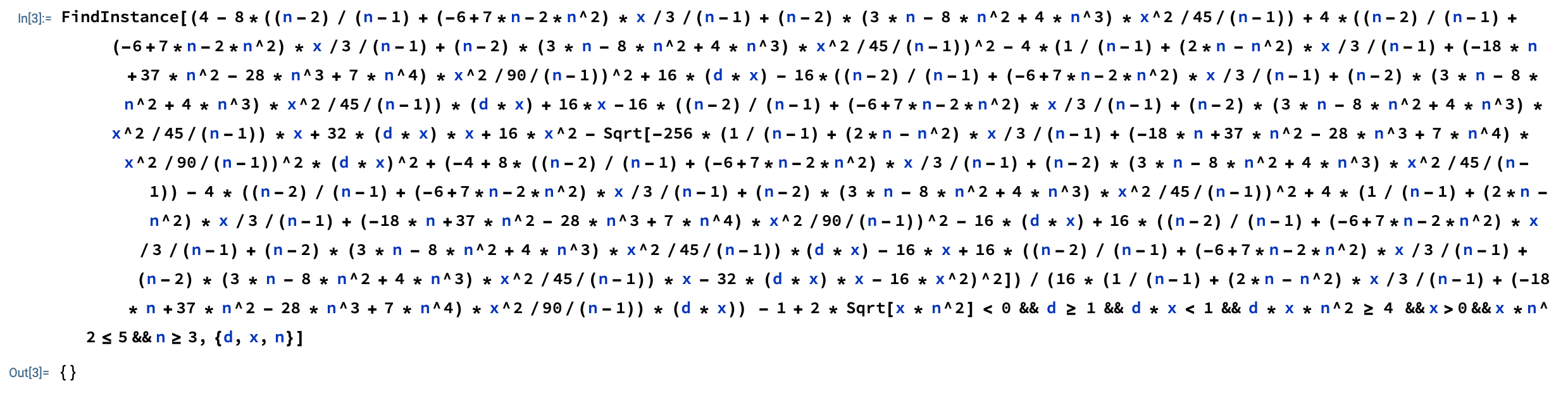}
\end{figure}
\end{itemize}

3) $\mu + \lambda \lambda_{\min}^+  > L$. In this case $c = \tfrac{L - \mu}{\lambda} = \tfrac{1}{x} \cdot \delta = x \delta$. We want to verify that $\gamma \in (0;1)$ and $\gamma \geq 1 - 3 \sqrt{\tfrac{\mu}{L - \mu}}$. This inequality need to be checked with the constraints: $\delta \geq 1$, $x = \tfrac{1}{x} > 0$, $x = \tfrac{1}{x}\tfrac{1}{x} \leq  \lambda^+_{\min} \leq \tfrac{5}{n^2}$ and $\delta < \tfrac{1}{x} \cdot \lambda^+_{\min} \leq \tfrac{1}{x} \cdot \tfrac{5}{n^2}$ (constraints of the considered case).
First, we check these inequalities when $\varphi = 0$:
\begin{itemize}[label=$\ast$]
\item $\gamma > 0$
\begin{figure}[H]
    \centering
    \includegraphics[width = 1\linewidth]{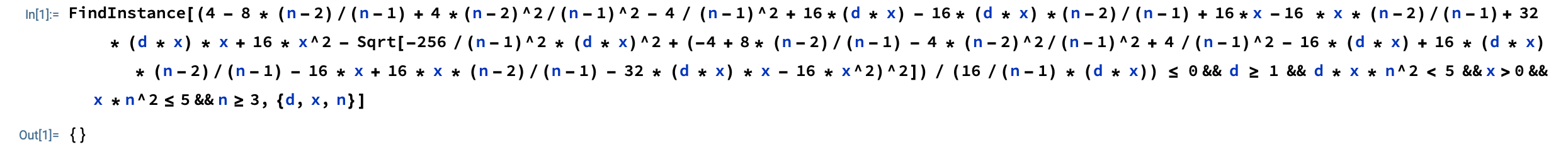}
\end{figure}
\item $\gamma < 1$
\begin{figure}[H]
    \centering
    \includegraphics[width = 1\linewidth]{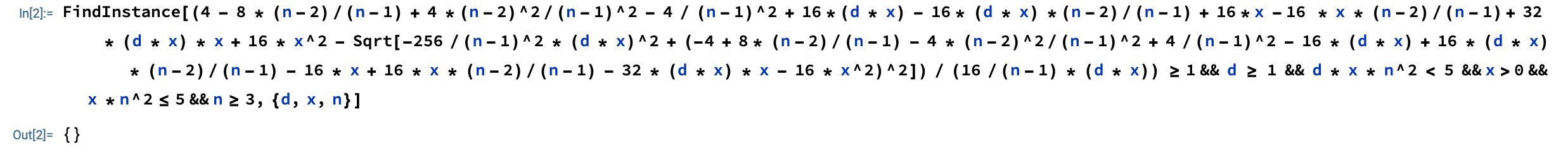}
\end{figure}
\item $\gamma \geq 1 - 3 \sqrt{\tfrac{\mu }{L - \mu}}$
\begin{figure}[H]
    \centering
    \includegraphics[width = 1\linewidth]{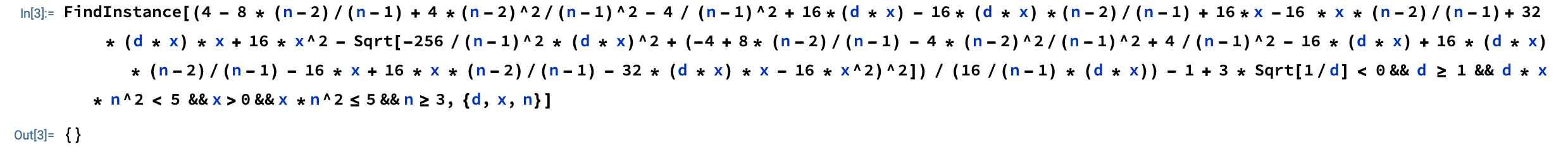}
\end{figure}
\end{itemize}

In the case of $\varphi = 1$, we use \eqref{AB}:
\begin{itemize}[label=$\ast$]
\item $\gamma > 0$
\begin{figure}[H]
    \centering
    \includegraphics[width = 1\linewidth]{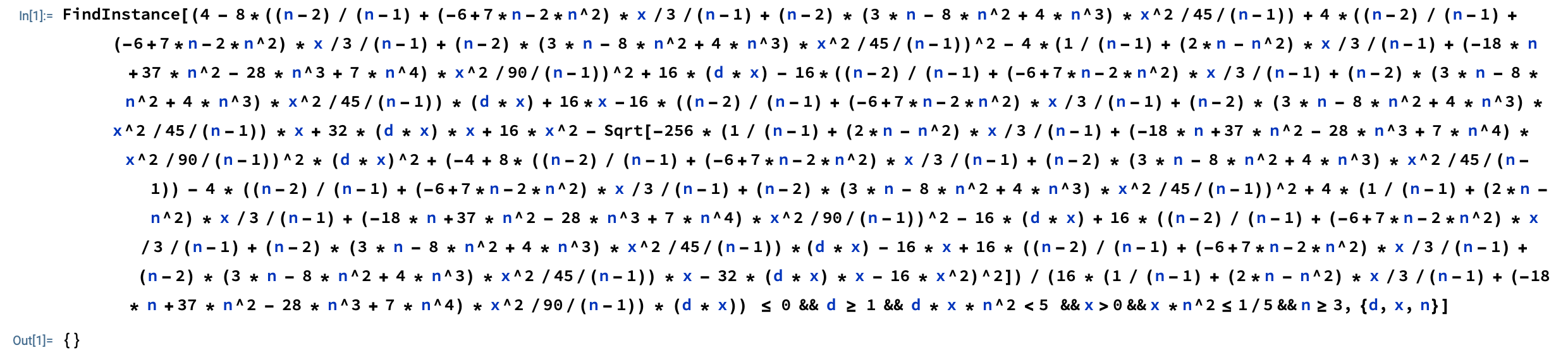}
\end{figure}
\item $\gamma < 1$
\begin{figure}[H]
    \centering
    \includegraphics[width = 1\linewidth]{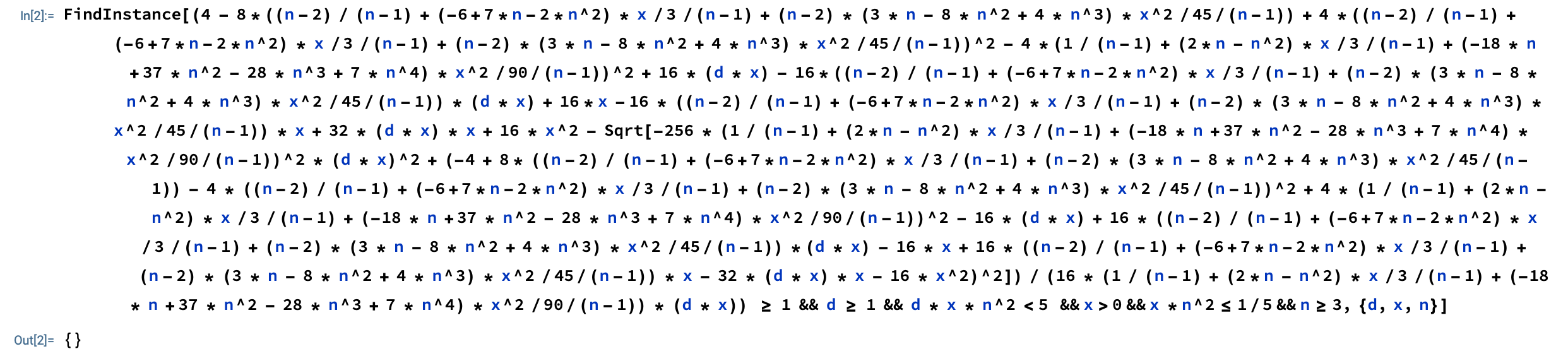}
\end{figure}
\item $\gamma \geq 1 - 3 \sqrt{\tfrac{\mu }{L - \mu}}$
\begin{figure}[H]
    \centering
    \includegraphics[width = 1\linewidth]{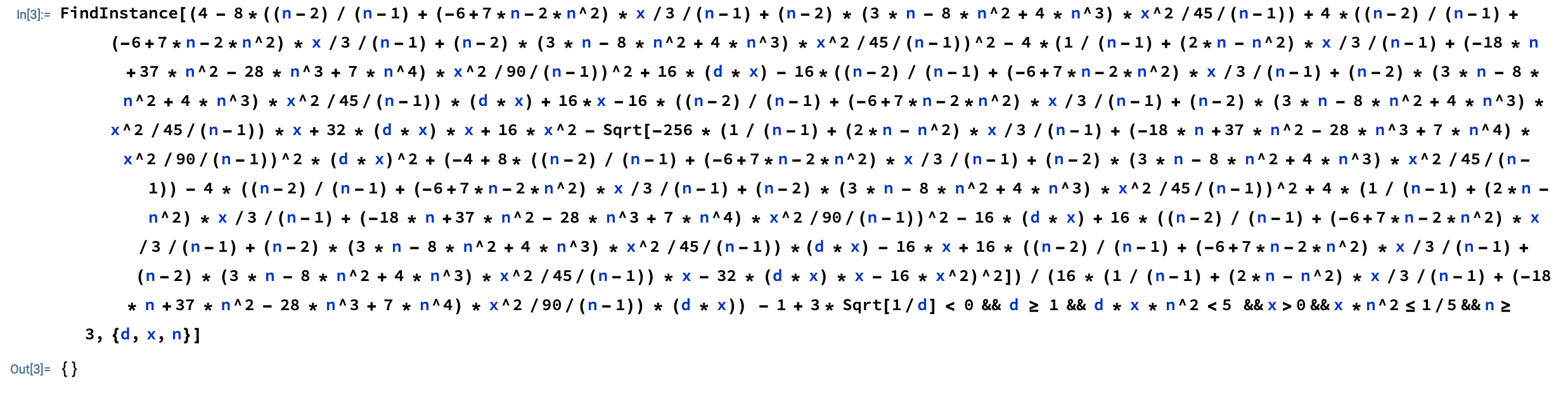}
\end{figure}
\end{itemize}
\EndProof
\end{proof}

The previous Lemmas show what the solution of the problem \eqref{rewrite_main_problem} + \eqref{eq:bad_func} is. Now let us determine how quickly we can approach it.

\begin{lemma} \label{lem_3404}
Let the problem \eqref{rewrite_main_problem} + \eqref{eq:bad_func} be solved by any method that satisfies Assumption \ref{as:lower}. Then after $K$ iterations with $q$ communication rounds,  only the first $\left\lfloor \tfrac{q}{n-1} \right\rfloor$ coordinates of the   global output can be non-zero while  the rest of the $d-\left\lfloor \tfrac{q}{n-1} \right\rfloor$ coordinates are strictly equal to zero.
\end{lemma}

\begin{proof}
We begin introducing some notation for our proof. Let
\begin{align*}
    E_{0} := \{ 0\}, \quad E_{j} := \text{span} \{ e_1, \ldots, e_j\}.
\end{align*}
Note that, if we initialize all $x^0_i = 0$, then we have $\mathcal{M}_{i,0} = E_0$.

Suppose that, at some given time $k$, for some $j$,  $\mathcal{M}_{j,k} = E_l$. Let us analyze how $\mathcal{M}_{j,k}$ can change by performing only local computations. 

We consider the case when $l$ odd (case with even $l$ can be analyzed the same way). After one local update, we have the following: 

1) For node $j \in \mathcal{V}_1$, it holds
\begin{equation}\begin{aligned}
\label{update_lower1}
        \mathcal{M}_{j,k+1} = E_l,
\end{aligned}\end{equation}
because of the block diagonal structure of \eqref{eq:bad_func}. The situation does not change, no matter how many local computations one  does.

2) For node $j \in \mathcal{V}_3$, it holds
\begin{equation*}\begin{aligned}
        \mathcal{M}_{j,k+1} = E_{l+1}, 
\end{aligned}\end{equation*}
It means that, after local computations, one has an update in output and machine on $\mathcal{V}_3$ can progress by one new non-zero coordinate.

This means that we constantly have to transfer progress from the machine from $\mathcal{V}_1$ to the machine from $\mathcal{V}_3$ and back. Initially, all devices have zero coordinates. Further, the machine from $\mathcal{V}_1$ can receive the first nonzero coordinate (but only the first, the second is not), and the rest of the devices are left with all zeros. Next, we pass the first non-zero coordinate to the machine from $\mathcal{V}_3$. To do this, $n-1$ communication rounds are needed. By doing so,  they can make the second coordinate non-zero, and then transfer this progress to the machine from $\mathcal{V}_1$. Then the process continues in the same way. This completes the proof.
\EndProof
\end{proof}

Now we are ready to complete the proof of Theorem \ref{theorem_lower_bounds}. The previous reasoning, as well as Lemmas \ref{lem:1404}, \ref{lem:2404}, and \ref{lem_3404}, gives that we can construct the "bad" problem of type \eqref{rewrite_main_problem} with the "bad" network (satisfying Definition \ref{gossipmatrix}) as well as with the "bad" functions \eqref{eq:bad_func} (satisfying Assumption \ref{as:func}). Moreover, we know that only $\left\lfloor \tfrac{q}{n-1} \right\rfloor$ coordinates in the output can coincide with the solution, and the other coordinates are exactly zero. Then we just have to put $T = \tfrac{1}{2}\left( \max \{ 1, \log_{\gamma} \tfrac{1}{2} \} + \left\lfloor \tfrac{q}{n-1} \right\rfloor \right)$ in the dimension of the problem $d = 2T$, and obtain the following estimate on the outputs from $\mathcal{V}_1$ and $\mathcal{V}_3$:
\begin{align*}
   \frac{\|x^K - x^*\|^2 + \|z^K - z^*\|^2}{\|x^0 - x^*\|^2 + \|z^0 - z^*\|^2} &= \frac{\sum_{i=\left\lfloor \tfrac{q}{n-1} \right\rfloor +1}^{2T} \| w_i\|^2}{\sum_{i=1}^{2T}\| w_i \|^2}
   \geq \frac{\sum_{i=\left\lfloor \tfrac{q}{n-1} \right\rfloor +1}^{2T} \gamma^{i-1} \| w_1\|^2}{\sum_{i=1}^{2T} \gamma^{i-1} \| w_1\|^2}\\
   & = \gamma^{\left\lfloor \tfrac{q}{n-1} \right\rfloor} \frac{\sum_{i=0}^{2T -1 - \left\lfloor \tfrac{q}{n-1} \right\rfloor} \gamma^{i}}{\sum_{i=0}^{2T-1} \gamma^{i}}
   = \gamma^{\left\lfloor \tfrac{q}{n-1} \right\rfloor} \frac{1 - \gamma^{2T - \left\lfloor \tfrac{q}{n-1} \right\rfloor}}{1 - \gamma^{2T}} \\
   & \geq \frac{1}{2} \gamma^{\left\lfloor \tfrac{q}{n-1} \right\rfloor} \geq \frac{1}{2} \left(1 - \max \left\{ 2\sqrt{\frac{\mu n^2 }{\lambda}}, 3\sqrt{\frac{\mu}{L - \mu}}\right\}\right)^{ \tfrac{q}{n-1}}.
\end{align*}
In the other words it means that:
\begin{align*}
    q &= \Omega\left( \min\left\{\sqrt{\frac{\lambda (n-1)^2}{\mu n^2 }}, \sqrt{\frac{(L - \mu)(n-1)^2}{\mu}}\right\}  \log\frac{(\|x^0 - x^*\|^2 + \|z^0 - z^*\|^2)}{\varepsilon}\right) \\
    &= \Omega\left( \min\left\{\sqrt{\frac{\lambda}{\mu }}, \sqrt{\frac{(L - \mu)(n-1)^2}{\mu}}\right\}  \log\frac{(\|x^0 - x^*\|^2 + \|z^0 - z^*\|^2)}{\varepsilon}\right)
\end{align*}
 When constructing the "bad" network, we proved that $n - 1 > \sqrt{2 \chi} - 2 \geq \tfrac{1}{5} \sqrt{\chi}$
$1 \leq \tfrac{2}{3} \lambda_{\max}(n)$. Hence, we get
\begin{align*}
    q &= \Omega\left( \min\left\{\sqrt{\frac{\lambda (n-1)^2}{\mu n^2 }}, \sqrt{\frac{(L - \mu)(n-1)^2}{\mu}}\right\}  \log\frac{(\|x^0 - x^*\|^2 + \|z^0 - z^*\|^2)}{\varepsilon}\right) \\
    &= \Omega\left( \min\left\{\sqrt{\frac{\lambda \lambda_{\max}(W)}{\mu }}, \sqrt{\frac{(L - \mu)\chi}{\mu}}\right\}  \log\frac{(\|x^0 - x^*\|^2 + \|z^0 - z^*\|^2)}{\varepsilon}\right).
\end{align*}
Which is what we needed to prove.
\EndProof

\section{Proof of Theorem \ref{t4}}\label{app:th_ma_LK}

For the following analysis, recall the auxiliary problem \eqref{auxiliary_prob} from Algorithm \ref{meta_algorithm} with $p = 1$, $h_1(\x)$ like sum component, $h_2(\x)$ like $\frac{\lambda}{2}\langle \x, W\x\rangle$, which is restated for convenience:
	\begin{equation*}%\label{auxiliary_prob}
	    \hat \y_{k+1} = \argmin_{\y\in \R^{nd}} \left\{ \langle \nabla h_1 (\w_k), \y-\w_k\rangle  + h_2(\y)   +\frac{\gamma}{2}\|\y-\w_k\|_2^{2} \right\}
	\end{equation*}

Now we look carefully at the auxiliary problem. This problem is $(L+\gamma)$-smooth and $(\mu + \gamma)$ strongly-convex, so we can apply L-Katyusha algorithm from \cite{L_Katyusha}. The complexity of solving problem \eqref{auxiliary_prob} is
\begin{equation*}
    \mathcal{O}\left(\left(M + \sqrt{\frac{M(\gamma + L)}{\gamma + \mu}}\right)\log\frac{1}{\delta}\right),
\end{equation*} 
where $\delta$ denotes the accuracy of the solution to the auxiliary problem \eqref{auxiliary_prob}. 
%\textbf{Case 2} $\lambda \lambda_{\max}(W) < L$
The number of calls of the gradient of $f$ is
\begin{equation}
    N_{W\x} = \mathcal{O}\left(\sqrt{\frac{\lambda \lambda_{\max}(W)}{\mu}}\log\frac{1}{\e}\right)
\end{equation}
while the number of calls of the gradient of $\mathcal{G}$ is
\begin{equation*}
    N_{\nabla f_k} = \mathcal{O}\left(\sqrt{\frac{\lambda \lambda_{\max}(W)}{\mu}}\left( M+ \sqrt{\frac{M(L + \gamma)}{\mu + \gamma}}\right)\log\frac{1}{\e}\log\frac{1}{\delta}\right).
\end{equation*}
Taking $\gamma$ be equal to $\lambda\lambda_{\max}\left(W\right)$, we get
\begin{equation*}
    N_{\nabla f_k} = \mathcal{O}\left(\left( M\sqrt{\frac{\lambda \lambda_{\max}(W)}{\mu}}+ \sqrt{\frac{ML}{\mu }}\right)\log\frac{1}{\e}\log\frac{1}{\delta}\right).
\end{equation*}

Now we consider $\delta$, the accuracy of the auxiliary problem \ref{auxiliary_prob}. According to Theorem \ref{theorem_meta_algorithm}, we can take $\delta$ as
\begin{equation*}
    \delta = \frac{\e\mu}{864^2(L + \lambda\lambda_{\max}(W) + \gamma)^2},
\end{equation*}
because the function $f(\x)$ is $L$-smooth and  $\frac{\lambda}{2}\left\la \x,W\x\right\ra$ is $\lambda\lambda_{\max}\left(W\right)$-smooth. \EndProof

\section{Proof of Theorem \ref{theorem_randomized}}\label{app:rand}
Let us use the additional notation $G (\x^k, \u^k) = \g^k +  \lambda W \u^k + \nabla f(\u^k) $ for short.
Let us consider our problem as a finite sum problem with $r+1$ terms. 
\begin{equation*}
    F(\x) = \frac{1}{r}\sum_{j=1}^r g_j(\x) + g_{r+1}(\x),
\end{equation*}
where $g_j(\x) =\sum\limits_{k=1}^n f_{jk}(x_k)$ and $g_{r+1}(\x) = \frac{\lambda}{2} \langle \x, W\x \rangle$. For such a problem, one can use the results of the convergence of the variance reduction method L-Katyusha (Algorithm 3 from \cite{L_Katyusha}) on which our method is based.
\begin{align*}
\EE&\left[ \|\g^k - \nabla F(\x^k)\|^2 \right] = \\
&= (1-p)\mathbb{E}\left[\left\| \frac{1}{1-p}\left(\nabla f_j(\x^k) - \nabla f_j(\u^k)\right) + \lambda W\u^k + \nabla f(\u^k) - \lambda W\x^k - \nabla f(\x^k)\right\|^2\right] \\
&+ p\mathbb{E}\left[\left\| \frac{\lambda}{p}\left(W\x^k - W\u^k\right) + \lambda W\u^k + \nabla f(\u^k) - \lambda W\x^k - \nabla f(\x^k)\right\|^2\right] \\
&= (1-p)\sum_{i=1}^n\sum_{j=1}^rp_{j} \left[\left\| \frac{1}{1-p}\left(\nabla f_{ij}(\x^k) - \nabla f_{ij}(\u^k)\right)\right.\right. \\&\left.\left.+ \lambda W\u^k + \nabla f_i(\u^k) - \lambda W\x^k - \nabla f_i(\x^k)\right\|^2\right] \\
&+ p\sum_{i=1}^n\left\| \frac{\lambda}{p}\left(W\x^k - W\u^k\right) + \lambda W\u^k + \nabla f_i(\u^k) - \lambda W\x^k - \nabla f_i(\x^k)\right\|^2 \\
&\leq (1-p)\sum_{i=1}^n\sum_{j=1}^rp_{j} \left\| \frac{1}{1-p}\left(\nabla f_{ij}(\x^k) - \nabla f_{ij}(\u^k)\right)\right\|^2 + p\left\| \frac{\lambda}{p}\left(W\x^k - W\u^k\right) \right\|^2\\
&\leq \frac{1}{1-p}\sum_{i=1}^n\sum_{j=1}^rp_{j} \left\|\nabla f_{ij}(\x^k) - \nabla f_{ij}(\u^k)\right\|^2 + \frac{2\lambda\lambda_{\max}(W)}{p}D_{g_{r+1}}(\u^k, \x^k).
\end{align*}
Choose $p_j = \frac{1}{r}$:
\begin{align*}
\EE&\left[ \|\g^k - \nabla F(\x^k)\|^2 \right] \\
&\leq \frac{2L}{1-p}\sum_{j=1}^nD_{f_j}(\u^k, \x^k) + \frac{ 2\lambda \lambda_{\max}(W)}{p} D_{g_{r+1}}(\u^k, \x^k) \\&=\frac{2L}{1-p}D_{f}(\u^k, \x^k) + \frac{2 \lambda \lambda_{\max}(W)}{p} D_{g_{r+1}}(\u^k, \x^k)
\\&\leq\max\left\{\frac{2L}{1-p}, \frac{ 2\lambda \lambda_{\max}(W)}{p} \right\} D_F(\u^k, \x^k)
\end{align*}

Choose $\mathcal{L} = \max\left\{\frac{L}{1-p}, \frac{ \lambda \lambda_{\max}(W)}{p} \right\}$, then,
\begin{align*}
\EE\left[ \|G^k - \nabla F(\x^k)\|^2 \right] \leq 2 \mathcal{L} D_F(\u^k, \x^k).
\end{align*}
Assumption 5.1 from \cite{L_Katyusha} holds. By Proposition 5.1 from \cite{L_Katyusha} iteration complexity of Algorithm \ref{alg_pfl_sum} is
\begin{equation*}
    O\left(\left(\frac{1}{\rho} + \sqrt{\frac{L + \lambda\lambda_{\max}(W)}{\mu}} +\sqrt{\frac{\mathcal{L}}{\rho\mu}} \right)\log\frac{1}{\varepsilon}\right)
\end{equation*}

Note that optimal complexities in Algorithm \ref{alg_pfl_sum} for local computations and communications 
are achieved on \textbf{different sets of $p$ and $\rho$}. Let us get them separately.

\begin{itemize}[label=$\ast$]
    \item The local stochastic gradient complexity of a single
iteration of Algorithm \ref{alg_pfl_sum} is 0 if $\xi^k = 0$, $\xi^{k + \frac{1}{2}} = 1$, $1$ if $\xi^k = 1$, $\xi^{k+\frac{1}{2}} = 1$, $r+1$ if $\xi^k = 1$, $\xi^{k+\frac{1}{2}} = 0$ and $M$ if $\xi^k = 0$, $\xi^{k+ \frac{1}{2}} = 0$. 
\begin{align*}
    \mathcal{ O}& \Bigg(\left((1 - p)(1 - \rho) + (M+1)(1 -p)\rho + M p \rho\right) \cdot \\
    &\quad\quad\quad\cdot\left(\frac{1}{\rho} + \sqrt{\frac{L + \lambda\lambda_{\max}(W)}{\mu}} +\sqrt{\frac{\mathcal{L}}{\rho\mu}} \right)\log\frac{1}{\varepsilon}\Bigg) 
    \\&\quad= \mathcal{ O} \left(\left(1 - p + M\rho \right)\left(\frac{1}{\rho} + \sqrt{\frac{L + \lambda\lambda_{\max}(W)}{\mu}} +\sqrt{\frac{\mathcal{L}}{\rho\mu}} \right)\log\frac{1}{\varepsilon}\right), 
\end{align*}
For $\rho = \frac{1}{M}$, $p = \frac{\lambda \lambda_{\max}(W)}{\bar{L} + \lambda \lambda_{\max}(W)}$ the
total expected local stochastic gradient complexity of Algorithm \ref{alg_pfl_sum} becomes
\begin{align*}
    \mathcal{ O} & \left(\left(1 - p + M\rho \right)\left(\frac{1}{\rho}+ \sqrt{\frac{L + \lambda\lambda_{\max}(W)}{\mu}} +\sqrt{\frac{\mathcal{L}}{\rho\mu}}\right)\log\frac{1}{\varepsilon}\right)
    \\& \leq \mathcal{ O} \left(2\left(M + \sqrt{\frac{L + \lambda\lambda_{\max}(W)}{\mu}} +\sqrt{\frac{M\mathcal{L}}{\mu}}\right)\log\frac{1}{\varepsilon}\right) 
    \\&= \mathcal{ O} \left(\left(M + \sqrt{ \frac{M(L + \lambda \lambda_{\max}(W))}{\mu}}\right)\log\frac{1}{\varepsilon}\right), 
\end{align*}

\item The total communication complexity of Algorithm \ref{alg_pfl_sum} is the sum of communication complexity coming from the full gradient computation (if statement that includes $\xi^{k + \frac{1}{2}}$) and the rest (if statement that includes $\xi^k$). The former requires a communication if $\xi^{k + \frac{1}{2}} = 0$,   the latter if $\xi^k$ is equal to $0$. The expected total communication $\mathcal{ O}\left(\rho + p\right)$ per iteration.  Thus, the total communication complexity is bounded by
\begin{align*}
    \mathcal{ O} \left((p + \rho)\left(\frac{1}{\rho} + \sqrt{\frac{L + \lambda\lambda_{\max}(W)}{\mu}} +\sqrt{\frac{\mathcal{L}}{\rho\mu}} \right)\log\frac{1}{\varepsilon}\right), 
\end{align*}
For $\rho = p$, $p = \frac{\lambda \lambda_{\max}(W)}{\bar{L} + \lambda \lambda_{\max}(W)}$ the total communication complexity of Algorithm \ref{alg_pfl_sum} becomes
\begin{align*}
    \mathcal{ O}& \left(\left(\rho +p\right)\left(\frac{1}{\rho}+\sqrt{\frac{L + \lambda\lambda_{\max}(W)}{\mu}} +\sqrt{\frac{\mathcal{L}}{\rho\mu}}\right)\log{\frac{1}{\varepsilon}}\right)
    \\&= \mathcal{ O} \left(\left(1+\rho\sqrt{\frac{L + \lambda\lambda_{\max}(W)}{\mu}} +\sqrt{\frac{\rho(L + \lambda\lambda_{\max}(W))}{\mu}}\right)\log{\frac{1}{\varepsilon}}\right)
     \\&= \mathcal{ O} \left(\sqrt{\frac{\lambda \lambda_{\max}(W)(L + \lambda\lambda_{\max}(W))}{(L + \lambda\lambda_{\max}(W))\mu}}\log{\frac{1}{\varepsilon}}\right)
     \\&= \mathcal{ O} \left(\sqrt{\frac{\lambda\lambda_{\max}(W) }{\mu}}\log{\frac{1}{\varepsilon}}\right). 
\end{align*}
\end{itemize}
%Let $\xi^k = \xi^{k+1} = 0$, $\xi^{k + \frac{1}{2}} = \xi^{k + 1  + \frac{1}{2}} = 0$.  Let us consider $d$ consecutive values of $\xi^k$, where $d$ is diameter of $W$. If there is only $0$, we have information on each node about all nodes. Then, if the next values of  $\xi^k$ is equal to $0$, we do not need to communicate. 
\EndProof

\end{document}